\documentclass[final,leqno]{siamonline0516}

\usepackage{amssymb,latexsym} 
\usepackage{amsmath}
\usepackage{caption}
\usepackage{subcaption}
\usepackage{textcomp}
\usepackage{algorithm}
\usepackage{siunitx}

\title{Stochastic basis adaptation and spatial domain decomposition for PDEs with random coefficients}

 \author{R. Tipireddy\footnotemark[2]
 \and P. Stinis\footnotemark[2]
 \and A. M. Tartakovsky\footnotemark[1] \footnotemark[2]
}

\begin{document}

\maketitle
\newcommand{\slugmaster}{%
  \slugger{juq}{xxxx}{xx}{x}{x--x}}
\renewcommand{\thefootnote}{\fnsymbol{footnote}}

\footnotetext[1]{Corresponding author.}
\footnotetext[2]{Pacific Northwest National Laboratory, P.O. Box 999, MSIN K7-90, Richland, WA 99352 ({\tt alexandre.tartakovsky@pnnl.gov}).}

\renewcommand{\thefootnote}{\arabic{footnote}}

\begin{abstract}

We present a novel uncertainty quantification approach for high-dimensional stochastic partial differential equations that reduces the computational cost of polynomial chaos methods by decomposing the computational domain into non-overlapping subdomains and adapting the stochastic basis in each subdomain so the local solution has a lower dimensional random space representation. The local solutions are coupled using the
Neumann-Neumann algorithm, where we first estimate the interface solution then evaluate the interior solution in each subdomain using the interface solution as a boundary condition.  The interior solutions in each subdomain are computed independently of each other, which reduces the operation count from $O(N^\alpha)$ to $O(M^\alpha),$ where $N$ is the total number of degrees of freedom, $M$ is the number of degrees of freedom in each subdomain, and the exponent $\alpha>1$  depends on the uncertainty quantification method used. In addition, the  localized  nature of solutions makes the proposed approach  highly parallelizable. We illustrate the accuracy and efficiency of the approach for linear and nonlinear differential equations with random coefficients. 

\end{abstract}


\pagestyle{myheadings}
\thispagestyle{plain}
\markboth{Tipireddy ET AL.}{Stochastic basis adaptation with domain decomposition}


\section{Introduction}
\label{RK:sec:intro}

We propose a novel approach within the framework of polynomial chaos (PC) expansion methods \cite{RK:Ghanem1991, RK:Xiu2002, Babuka2002, Babuka2010} for solving stochastic partial differential equations (PDEs) with a large number of random input parameters. This approach allows us to address some of the computational challenges due to the high dimensionality in PC-based uncertainty quantification (UQ) methods. The dimension, $d,$ of the PC expansion is obtained from the truncated Karhunen-Lo\`eve (KL) expansion  \cite{RK:Loeve1977} of the underlying random field. As $d$ increases, the number of terms in the PC expansion and the stochastic system size increases exponentially, making PC methods computationally intractable \cite{Nouy2007, RK:Doostan2011, Lin2009AWR,Lin2010JSC,Venturi2013JCP, Tipireddy2013}. 

In \cite{Soize2016} data driven approach is presented to construct reduced order models that are statistically consistent with the given data set. Algorithms for a polynomial chaos expansion of multimodal random vector is presented in \cite{Soize2015} in order to reduce the dimensionally of the stochastic system. The effect of geometric transformation such as in manufacturing process on random material properties, its mathematical treatment and applicability of low dimensional probabilistic methods is discussed in \cite{Ghanem2015}.

Naturally, one approach for reducing the number of terms in PC expansion is to reduce $d$. Recently, we have proposed the basis adaptation method for stochastic dimension reduction while solving for a specific quantity of interest (QoI) {\cite{Tipireddy2013, RK:Tipireddy2014}}, i.e, for the scalar solution of a steady-state PDE at a selected spatial location and for a scalar QoI that is a nonlinear function of the solution of the steady-state PDE. In \cite{RK:Tsilifis2016}, the basis adaptation method has been further extended from a scalar- to vector-valued QoI. Here, we propose a method for solving high-dimensional stochastic problems in the entire spatial domain (as opposed to computing QoIs). Our method is based on spatial domain decomposition, which replaces a stochastic PDE with a system of PDEs coupled through a corresponding interfacial problem. A comprehensive review of domain decomposition methods and algorithms for deterministic PDEs is provided in \cite{RK:Toselli2005}. Domain decomposition also has been used for stochastic PDEs. For example, \cite{Xiu2004662} and \cite{Lin2010JCP} apply domain decomposition to reduce the variance of an input parameter field by decomposing the highly heterogeneous domain in mildly heterogeneous subdmains with random boundaries. Pranesh et al. \cite{Pranesh2016} theoretically demonstrate that the random dimensionality of a parameter field decreases with the domain size. Chen et al. \cite{Chen2015} use this property of fields to formulate a smaller-dimensional stochastic problem in each subdomain. 

Our novel idea is to formulate and solve a stochastic problem in each subdomain in its (new) smaller-dimensional random space using the basis adaptation method such that the random variables in each subdomain maintain dependency on random variables in the other subdomains. In \cite{RK:Tipireddy2016}, we obtain low-dimensional  {\it local} representations of the random solution in each subdomain and are able to reconstruct a global solution by {\it stitching together} the local solutions from the subdomains. To obtain a low-dimensional representation in each subdomain, we employ a Hilbert space KL expansion \cite{RK:Doostan2007}. For each of the locally adapted bases, we solve the corresponding equations in the {\it whole} domain, keep the local solution for each subdomain, and stitch together these local solutions to obtain the global solution. This approach makes the numerical implementation of the method straightforward and efficient if the solutions in the whole domain can be obtained reasonably fast.

Unlike the approach followed in \cite{RK:Tipireddy2016}, where the stochastic PDE is solved in the whole computational domain for each local basis, here we solve the stochastic PDE in each subdomain in its adapted local basis. The stochastic PDE in each subdomain is coupled with the PDEs in the adjacent subdomains through continuity conditions at the interface. Because the interface condition from the adjacent subdomains is expressed in a different local stochastic basis, we have to project the interface conditions to represent them in terms of the local basis adapted to that subdomain. The advantages of this new approach include the continuity of solution and its derivatives across the subdomain boundaries and high scalability of resulting algorithm. The latter is achieved by solving independent boundary value problems in each subdomain. We demonstrate our method by solving three different boundary value problems: i) a one-dimensional linear diffusion equation with statistically non-stationary diffusion coefficient (saturation-based Richards equation with discontinuity of the conductivity in physical domain), ii) a one-dimensional nonlinear diffusion equation with statistically non-stationary diffusion coefficient  (pressure based one-dimensional nonlinear Richards equation),  and iii) a two-dimensional linear diffusion equation with point sink at the center of the spatial domain. For the particular type of equations in our numerical examples, there are several methods for implementing domain decomposition and obtaining the solution in each subdomain \cite{RK:Toselli2005}. In this work, we use a Neumann-Neumann (N-N) algorithm (see Section \ref{RK:sec:NN_section} for more details).  We present a detailed computational cost analysis of the proposed method that shows its cost is significantly smaller than that of the traditional PC based collocation method. 

This paper is organized as follows: In Section \ref{RK:sec:spde}, we present traditional PC-based UQ methods for stochastic PDEs and motivate the need for new approaches to reduce the computational cost. In Section \ref{RK:sec:ba_and_dd}, we discuss our methodology for spatial domain decomposition and basis adaptation methods to represent the local solution in a low-dimensional stochastic space. Section \ref{RK:sec:NN_section} describes the boundary value problems for subdomains and also provides a description of the the N-N algorithm for solving the stochastic PDE in low-dimensional random space. In Section \ref{RK:sec:numerical}, we present numerical results for three different examples using the aforementioned method and provide computational cost analysis for the two-dimensional steady-state diffusion equation with random diffusion coefficient. We present conclusions and ideas for future work in Section \ref{conclusions}.

\section{Partial differential equations with random coefficients}
\label{RK:sec:spde}
We propose a domain decomposition basis adaptation method for high-dimensional linear and nonlinear stochastic PDEs. While there is no inherent limitations of the method for a particular type of PDEs,  we focus this work on nonlinear and linear steady-state diffusion equations because of their ubiquity in natural and engineered systems.  Specifically, we are interested in computing the solution $u(x,\omega):D \times \Omega \rightarrow \mathbb{R}$ of the stochastic partial differential equation (SPDE)
\begin{align}\label{RK:eq:spdeop}
 \mathcal{L}(x,u(x,\omega);a(x,\omega)) &= f(x,\omega)  \;\; \rm{in}~D\times \Omega, \nonumber \\
 \mathcal{B}(x,u(x,\omega);a(x,\omega)) &= h(x,\omega)  \;\; \rm{on}~\partial D\times \Omega,
\end{align}
with the random coefficient $a(x,\omega),$ where $D$ is an open subset of $\mathbb{R}^n$ and $\Omega$ is a sample space. 
Here, $\mathcal{L}$ is a differential operator, and $\mathcal{B}$ is a boundary operator. 
 To numerically solve the stochastic PDE (\ref{RK:eq:spdeop}), it is common to discretize the random fields $a(x,\omega)$ and $u(x,\omega)$ in both spatial and stochastic domains ~\cite{Ghanem1999}. We model $a(x,\omega)$ as a log-normal random field, which is $a(x,\omega) = \exp[g(x,\omega)]$, where $g(x,\omega)$ is a Gaussian random field with known mean and covariance function. We approximate $g(x,\omega)$ with a truncated KL expansion~\cite{RK:Loeve1977} as
 \begin{equation}\label{RK:eq:gkl}
 g(x,\omega) \approx g(x, \boldsymbol{\xi}(\omega)) = g_0(x) + \sum_{i=1}^d \sqrt{\lambda_i} g_i(x) \xi_i(\omega),
\end{equation}
where $d$ is the number of random variables in the truncated expansion. The choice of $d$ depends on the decay of the eigenvalues of the covariance function of $g(x, \omega)$ in domain $D$. We choose $d$ such that the contribution from eigenvalues $\{\lambda_{d+1}, \cdots, \lambda_{\infty}\}$ in~\eqref{RK:eq:gkl} toward $ g(x,\omega)$ is quite small and can be neglected. In general, the  smaller the correlation length of the covariance function of $g(x, \omega)$, the larger the dimension $d$. The independent uncorrelated Gaussian random variables $\boldsymbol{\xi} = (\xi_1, \ldots, \xi_d)^T$  have zero mean and unit variance. Here, $g_0(x)$ is the mean of the random field $g(x,\omega)$, and the eigenvalues $\{\lambda_i\}$ and eigenfunctions $\{g_i(x)\}$ are computed as a solution of the eigenvalue problem
 \begin{equation}\label{RK:eq:eig}
 	\int_D C_g(x_1,x_2)  g_i(x_2)dx_2 = \lambda_i g_i(x_1),
\end{equation}
where $C_g(x_1,x_2)$ is the prescribed covariance function of $g(x,\omega)$. Because the covariance function is positive definite, the eigenvalues are positive and non-increasing, and the eigenfunctions $g_i(x)$ are orthonormal, that is:
 \begin{equation}\label{RK:eq:ortho}
 	\int_D g_i(x) g_j(x)dx = \delta_{ij},
\end{equation}
where $\delta_{ij}$ is the Kronecker delta.
The solution $u(x,\omega)$ is approximated with a truncated PC expansion \cite{Cameron1947, Ghanem1999} as
\begin{equation}\label{RK:eq:pce_u}
  u(x,\omega) \approx u(x,\boldsymbol{\xi}(\omega))  =  u_0(x) + \sum_{i=1}^{N_{\xi}} u_i(x) \psi_i(\boldsymbol{\xi}),
\end{equation}
where  $N_{\xi} = \left(\frac{(d+p)!}{d!~p!}-1\right)$ is the number of terms in PC expansion for dimension $d$ and order $p$, $u_0(x)$ is the mean of the solution field, $u_i(x)$ are PC coefficients, and $\{\psi_i(\boldsymbol{\xi})\}$ are multivariate Hermite polynomials. The PC basis, $\{\psi_i(\boldsymbol{\xi})\}$, are also orthonormal, that is:  
\begin{equation}\label{RK:eq:innprod}
\langle \psi_i(\boldsymbol{\xi}), \psi_j(\boldsymbol{\xi}) \rangle \equiv \int_{\Omega} \psi_i(\boldsymbol{\xi}) \psi_j(\boldsymbol{\xi}) d\boldsymbol{\xi} = \delta_{ij},
\end{equation}
where the integration is performed with respect to the Gaussian measure.
Here, $u_i(x)$ are the unknown coefficients to be computed as a solution of~\eqref{RK:eq:spdeop}. 
Once the random coefficient $a(x,\omega)$ and the solution field $u(x,\omega)$ are represented in terms of $\boldsymbol{\xi}$, the stochastic PDE \eqref{RK:eq:spdeop} transforms into 
\begin{align}\label{RK:eq:spdexi}
 \mathcal{L}(x,\boldsymbol{\xi}, u(x,\boldsymbol{\xi});a(x,\boldsymbol{\xi})) &= f(x,\boldsymbol{\xi})  \;\; \rm{in}~D\times \Omega, \nonumber \\
 \mathcal{B}(x,\boldsymbol{\xi}, u(x,\boldsymbol{\xi});a(x,\boldsymbol{\xi})) &= h(x,\boldsymbol{\xi})  \;\; \rm{on}~\partial D\times \Omega.
\end{align}

Intrusive methods, such as stochastic Galerkin~\cite{RK:Ghanem1991}, or non-intrusive methods, such as sparse-grid collocation~\cite{RK:Xiu2002}, can be used to solve the parameterized stochastic PDE \eqref{RK:eq:spdexi}. In the current work, we use a non-intrusive approach, where the PDE is solved at predefined quadrature points (in random space). A popular choice of quadrature points is sparse-grid collocation points based on the Smolyak approximation~\cite{RK:Smolyak1963, RK:Nobile2008}. Unlike the tensor product of one-dimensional quadrature points, the sparse-grid method judiciously chooses products with only a small number of quadrature points. These product rules depend on an integer value called a {\it sparse-grid level}~\cite{RK:Smolyak1963, RK:Nobile2008}. As the order of PC expansion increases, higher sparse-grid levels are required to maintain solution accuracy. In the sparse-grid method, the number of collocation points increases with the stochastic dimension (number of random variables) and  the sparse-grid level. Let $\boldsymbol{\xi}_q$ be the collocation points associated with the random variables $\boldsymbol{\xi}$. Then, the deterministic PDE 
\begin{align}\label{RK:eq:spdexi_q}
 \mathcal{L}(x, u_q(x);a_q(x)) &= f_q(x)  \;\; \rm{in}~D, \nonumber \\
 \mathcal{B}(x, u_q(x);a_q(x)) &= h_q(x)  \;\; \rm{on}~\partial D
\end{align}
is solved for each collocation point ${{\boldsymbol{\xi}}_q}$ ($q = 1,\ldots, Q_{\xi}$), where $u_q(x) = u(x,\boldsymbol{\xi}_q),$  $a_q(x) = a(x,\boldsymbol{\xi}_q),$  $f_q(x) = f(x,\boldsymbol{\xi}_q),$ and  $h_q(x) = h(x,\boldsymbol{\xi}_q)$. The PC coefficient ${u}_i$ is computed through projection,
\begin{equation}\label{RK:eq:ui_xi}
	u_i(x) = \sum_{q=1}^{Q_{\xi}} {u}(x, \boldsymbol{\xi}_q) \psi_i(\boldsymbol{\xi}_q) w_q^{\xi},
\end{equation}
where $w_q^{\xi}$ are weights for the quadrature points. However, with these methods, the computational cost increases exponentially with increasing $d$ and/or $p.$ Here, we propose a novel approach that divides the spatial domain into several subdomains \cite{RK:Tipireddy2016} and employs the basis adaptation method \cite{RK:Tipireddy2014} in each subdomain to reduce the stochastic dimension. In the proposed method, the solutions in subdomains are coupled through interface boundary conditions instead of ``stitching'' together the local solutions as in \cite{RK:Tipireddy2016}. This affords the ability to independently solve the equations in different subdomains, making the approach highly parallelizable. With this approach, the computational savings are two-fold because: i) the reduction in the stochastic dimension in each subdomain and ii) the reduced computational domain for each subdomain. However, there is an additional cost involved in imposing the boundary conditions at the interface, which is discussed in the following sections.

\section{Domain decomposition and local basis adaptation} \label{RK:sec:ba_and_dd}

To reduce the computational cost due to high dimensionality in solving~\eqref{RK:eq:spdexi}, we first decompose the spatial domain $D \subset \mathbb{R}^n$ into a set of non-overlapping subdomains $D_{s} \subset D, s = 1, \ldots, N_D,$ that is:
\begin{equation}\label{RK:eq:domaind}
	{D} = \bigcup_{s=1}^{N_D} {D_s}, \quad D_s \cap D_{s^{\prime}}=\emptyset \; \text{for} \; s \neq s^{\prime}. 
\end{equation}	
Then, we use the basis adaptation in each subdomain $D_s$ to find a low-dimensional stochastic basis (PC) with a new set of random variables $\tilde{\boldsymbol{\eta}}^s = \{\eta_1^s, \ldots, \eta_{r}^s \}, r \ll d.$

In Section \ref{RK:sec:hkle}, we describe the basis adaptation method employed in each subdomain using the Hilbert-Karhunen-Lo\`eve (Hilbert KL) expansion of the Gaussian part of the solution. We assume that the Gaussian part of the solution can be obtained with a small cost by solving the stochastic system in the original dimension only up to the first order.  
\label{RK:sec:stochreduct}

\subsection{Hilbert-Karhunen-Lo\`eve expansion}
\label{RK:sec:hkle}
We combine the basis adaptation method ~\cite{RK:Tipireddy2014} and spatial domain decomposition~\cite{RK:Tipireddy2016} to construct a low-dimensional solution representation in each subdomain $D_s$ using the Hilbert KL expansion.
The KL expansion ~\eqref{RK:eq:gkl} provides optimal representation of a random field in $L_2$ space ~\cite{RK:Loeve1977}. We want an  optimal representation of the solution space. Hence, the solution should satisfy certain regularity and smoothness conditions~\cite{RK:Doostan2007}. Namely, the solution $u(x,\omega)$ should be a subset of $L_2(\Omega)$.  Here, we use the Hilbert KL expansion \cite{RK:Levy1999, RK:Kirby1992, RK:Silverman1996, Berkooz1993539, Christensen1999} of $u(x,\omega),$ where the Gaussian part of the solution $u_g(x,\omega)$ is expanded in terms of a new set of uncorrelated Gaussian random variables, $\{\eta_i\}.$ 

We first compute the Gaussian part (which requires only up to linear terms in the PC expansion) of the solution in the {\it entire} spatial domain $D$ by solving~\eqref{RK:eq:spdexi}, i.e., 
\begin{equation}\label{RK:eq:pce_ug2}
	 u_g(x,\boldsymbol{\xi}(\omega)) = u_0(x) + \sum_{i=1}^{d} u_i(x) \xi_i,
\end{equation}
where, $u_i(x) $ is computed using \eqref{RK:eq:ui_xi}. Note that solving for $u_g$ is computationally less expensive than finding the full solution of \eqref{RK:eq:spdexi} and could be feasible even if the full solution computation is not. 

Consider the Gaussian part of the solution in subdomain $D_s$, $u^s_g(x,\boldsymbol{\xi}(\omega))$, such that
\begin{equation}\label{RK:eq:pce_ug_s}
	u^s_g(x,\boldsymbol{\xi}(\omega)) = u_g(x,\boldsymbol{\xi}(\omega)) \mathbb{I}_{D_s}(x), 
\end{equation}
where $\mathbb{I}_{D_s}(x)$ is the indicator function so that for any set $D_s,$ $\mathbb{I}_{D_s} = 1$ if $x \in D_s,$ and $\mathbb{I}_{D_s} = 0$ if $x \notin D_s$. We construct the covariance function of $u^s_g(x,\boldsymbol{\xi})$ in each subdomain $D_s$ as follows: 
\begin{equation}\label{RK:eq:cov_ug}
	 C^s_{u_g}(x_1,x_2) = \sum_{i=1}^{d} u_i(x_1) u_i(x_2), \quad x_1, x_2 \in {D_s}.
\end{equation}
The Hilbert space KL expansion of $u^s_g(x,\boldsymbol{\xi})$  \cite{RK:Doostan2007}) in subdomain $D_s$ is
\begin{equation}\label{RK:eq:ugkl}
	u^s_g(x,\boldsymbol{\xi}(\omega)) = u^s_0(x) + \sum_{i=1}^{d} \sqrt{\mu^s_i} \phi^s_i(x) \eta^s_i(\omega), \quad x \in {D_s},
\end{equation}
where $(\mu^s_i, \phi^s_i(x))$ are eigenpairs in the Hilbert KL expansion. They can be obtained by solving the eigenvalue problem: 
\begin{equation}\label{RK:eq:cov_cg}
	 \int_{D_s} C^s_{u_g}(x_1,x_2) \phi^s_i(x_1)dx_1 = \mu^s_i \phi^s_i(x_2), \quad i = 1, 2, \ldots, d. 
\end{equation}
The eigenvalue problem \eqref{RK:eq:cov_cg} can be solved numerically by first discretizing the covariance function $C^s_{u_g}(x_1,x_2)$ and eigenfunctions $\phi^s_i(x)$ using a suitable spatial basis, such as finite element basis, to obtain matrix eigenvalue equations. 
We can use \eqref{RK:eq:ugkl} to write the random variable $\eta^s_i$ as
\begin{align}\label{RK:eq:pce_eta}
	 \eta^s_i &= \frac{1}{\sqrt{\mu^s_i}}\int_{D_s} \left(u^s_g(x,\boldsymbol{\xi}) - u^s_0(x) \right) \phi^s_i(x)dx, \nonumber \\
	 		&= \frac{1}{\sqrt{\mu^s_i}}\int_{D_s} \left(u_0(x) + \sum_{j=1}^{d} u_j(x) \xi_j - u^s_0(x) \right) \phi^s_i(x)dx. \quad i = 1, 2, \ldots, d.
\end{align}
Because $u_0(x) = u^s_0(x), x\in D_s$ by construction,  we have 
 \begin{align}\label{RK:eq:eta_xi}
	 \eta^s_i &= \frac{1}{\sqrt{\mu^s_i}}\int_{D_s} \left ( \sum_{j=1}^{d} u_j(x) \xi_j \right) \phi^s_i(x)dx, \quad x \in D_s, i = 1, 2, \ldots, d \nonumber \\
	 		&= \sum_{j=1}^{d} \left ( \frac{1}{\sqrt{\mu^s_i}}\int_{D_s}  u_j(x) \phi^s_i(x)dx  \right)\xi_j, \quad x \in D_s, i = 1, 2, \ldots, d \nonumber \\
	 		&= \sum_{j=1}^{d} a^s_{ij} \xi_j,
\end{align}
where $a^s_{ij} =  \frac{1}{\sqrt{\mu^s_i}}\int_{D_s}  u_j(x) \phi^s_i(x)dx, i,j = 1, \ldots, d.$ $a^s_{ij}$ provides a linear map between $\eta^s_i$ and $\{\xi_j\}$. The resulting normal variables $\{\eta^s_i\}$ can be normalized to get independent, uncorrelated standard Gaussian random variables because the $\xi_j$ are standard Gaussian random variables. In each subdomain, we reformulate the stochastic PDE~\eqref{RK:eq:spdeop} in terms of $\{\eta^s_i\}$ and solve using the non-intrusive method.

\subsection{Dimension reduction in the subdomain $D_s \subset D$}
\label{RK:sec:ba}
Equation \eqref{RK:eq:eta_xi} can be rewritten as 
\begin{equation}\label{RK:eq:eta}
	\boldsymbol{\eta}^s =A_s \boldsymbol{\xi}, \quad A_s {A_s}^T = \boldsymbol{I}, 
\end{equation}
where $A_s = [a^s_{ij}]$ is an isometry in $\mathbb{R}^d$ and $\boldsymbol{\eta}^s = \{\eta^s_1, \ldots, \eta^s_d\}^T$ is a vector of standard normal random variables. 
We can prove that $A_s {A_s}^T = \boldsymbol{I}$ by showing $\sum_{k=1}^d a^s_{ik}a^s_{jk} = \delta_{ij}$, that is,
 \begin{align}\label{RK:eq:aij_aji}
	 \sum_{k=1}^d a^s_{ik} a^s_{jk} &=  \sum_{k=1}^d \left (\frac{1}{\sqrt{\mu^s_i}}\int_{D_s}  u_k(x) \phi^s_i(x) dx \right ) \left (\frac{1}{\sqrt{\mu^s_j}}\int_{D_s}  u_k(y) \phi^s_j(y) dy \right )  \nonumber \\	 
	 		&= \frac{1}{\sqrt{\mu^s_i}\sqrt{\mu^s_j}} \int_{D_s} \int_{D_s} \left ( \sum_{k=1}^d u_k(x) u_k(y)  \right )\phi^s_i(x)   \phi^s_j(y) dx dy. 
\end{align}
 From \eqref{RK:eq:cov_ug}, $\sum_{k=1}^d  u_k(x) u_k(y) =  C_{u_g} (x,y)$. Hence, 
 \begin{align}\label{RK:eq:aij_aji2}
	 \sum_{k=1}^d a^s_{ik} a^s_{jk} & = \frac{1}{\sqrt{\mu^s_i}\sqrt{\mu^s_j}} \int_{D_s} \left ( \int_{D_s}  C_{u_g} (x,y)  \phi^s_i(x)  dx \right )  \phi^s_j(y) dy,
\end{align}
and from \eqref{RK:eq:cov_cg}, $\int_{D_s} \left ( C_{u_g} (x,y) \phi^s_i(x)  dx  \right ) = \mu_i^s \phi_i^s(y),$ hence, 
 \begin{align}\label{RK:eq:aij_aji3}
	 \sum_{k=1}^d a^s_{ik} a^s_{jk} & = \frac{1}{\sqrt{\mu^s_i}\sqrt{\mu^s_j}} \int_{D_s} \mu_i^s \phi^s_i(y)  \phi^s_j(y) dy  \nonumber \\
	 		 &= \delta_{ij}, \quad \text{after normalization}
\end{align}
where, $\delta_{ij}$ is the Kronecker delta.

The mapping described in \eqref{RK:eq:eta} suggests that both  $\boldsymbol{\xi}$ and $\boldsymbol{\eta}^s$ span the same Gaussian Hilbert space, and the solution can be written as
\begin{equation}\label{RK:eq:u_tilde}
	u(x,\boldsymbol{\xi}) = {u}^{A_s}(x,\boldsymbol{\eta}^s(\boldsymbol{\xi}))
	= {u}^{A_s}_0(x) + \sum_{i=1}^{N_{\eta^s}} {u}^{A_s}_i(x) \psi_i(\boldsymbol{\eta}^s), \quad x \in D_s,
\end{equation}
where $N_{\eta^s} = \left(\frac{(d+p)!}{d!~p!}-1\right)$ is the number of terms in the PC expansion for dimension $d$ and order $p.$

For the proposed approach to be computationally efficient, it must represent the solution accurately in the subdomain $D_s$ with a smaller dimension $r<d.$ Because the eigenvalues in the Hilbert KL expansion decay faster than the original KL expansion, fewer terms are sufficient to represent the solution accurately. Two factors contribute to the faster decay of the eigenvalues. The first factor is that the solution generally is smoother than the input random field $a(x,\boldsymbol{\xi}),$ while the second factor stems from the reduction in the domain size. Although the correlation length $L_c$ remains the same for the full domain $D$ and the subdomain $D_s,$ the ratio of the correlation length to the domain size increases for the subdomain, meaning $\frac{L_c}{D_s}>\frac{L_c}{D}.$ As the correlation length relative to domain size increases or as the domain size decreases relative to the correlation length, the eigenvalues in the Hilbert KL expansion decay faster. Thus, a smaller dimension is sufficient. Even when the solution is not smoother relative to the input, the second factor will help reduce the dimension of the adapted basis in a subdomain.

For each subdomain $D_s,$ \eqref{RK:eq:spdexi} can be reformulated with the new set of random variables $\tilde{\boldsymbol{\eta}}^s = \{\eta_1, \ldots, \eta_r \}$ as 

\begin{align}\label{RK:eq:spdeop_sub}
 \mathcal{L}(x,u^s(x,\tilde{\boldsymbol{\eta}}^s);a(x,\boldsymbol{\xi})) &= f(x,\boldsymbol{\xi})  \;\; \rm{in}~D_s\times \Omega, \nonumber \\
 \mathcal{C}(x,u^s(x,\tilde{\boldsymbol{\eta}}^s),u^{s^{\prime}}(x,\tilde{\boldsymbol{\eta}}^s);a(x,\boldsymbol{\xi})) &= 0  \;\; \rm{in}~(\partial D_s \cap \partial D_{s^{\prime}}) \times \Omega, \nonumber \\
 \mathcal{B}(x,u^s(x,\tilde{\boldsymbol{\eta}}^s);a(x,\boldsymbol{\xi})) &= h(x,\boldsymbol{\xi})  \;\; \rm{on}~(\partial D_s \cap \partial D) \times \Omega.
\end{align} 
The coupling operator $\mathcal{C}$ is used to find the interface solution between subdomains $D_s$ and $D_{s^{\prime}},$ while $\cal{B}$ is, as before, the boundary operator.

\subsection{Subdomain solution using non-intrusive methods}\label{RK:subsec:nonintrusive}
In a non-intrusive method, such as the sparse-grid method, the number of collocation points increases with the stochastic dimension (number of random variables) and  the sparse-grid level. Let $\tilde{\boldsymbol{\eta}}_q^s$ be the collocation points associated with the random variables $\tilde{\boldsymbol{\eta}}^s$. Then, the corresponding points ${\boldsymbol{\xi}}_q^s$ are obtained as 
\begin{equation}\label{RK:eq:xi_Aeta}
	{\boldsymbol{\xi}}^s_q =   [A_s^{-1}]_{r} \tilde{\boldsymbol{\eta}}^s_q,
\end{equation}
where the $d \times r$ matrix $[A_s^{-1}]_r$ consists of the first $r$ columns of $A_s^{-1}$. Finally, the deterministic PDE 
\begin{align}\label{RK:eq:spdeop_sub_quad}
 \mathcal{L}\left(x,u^s(x,\boldsymbol{\xi}_q^s);a(x,\boldsymbol{\xi}_q^s)\right) &= f(x,\boldsymbol{\xi}_q^s)  \;\; \rm{in}~D_s\times \Omega, \nonumber \\
 \mathcal{C}\left(x,u^s(x,\boldsymbol{\xi}_q^s),u^{s^{\prime}}(x,\boldsymbol{\xi}_q^s);a(x,\boldsymbol{\xi}_q^s)\right) &= 0  \;\; \rm{in}~(\partial D_s \cap \partial D_{s^{\prime}}) \times \Omega, \nonumber \\
 \mathcal{B}\left(x,u^s(x,\boldsymbol{\xi}_q^s);a(x,\boldsymbol{\xi}_q^s)\right) &= h(x,\boldsymbol{\xi}_q^s)  \;\; \rm{on}~(\partial D_s \cap \partial D) \times \Omega
\end{align} 
is solved for each collocation point ${{\boldsymbol{\xi}}_q^s}$ ($q = 1,\ldots, Q_{\eta^s}$), and the PC coefficient $\tilde{u}^{A_s}_i$ is computed by projection
\begin{equation}\label{RK:eq:ui_eta}
	\tilde{u}^{A_s}_i(x) = \sum_{q=1}^{Q_{\eta^s}} {u}(x, {\boldsymbol{\xi}}^s_q) \psi_i(\tilde{\boldsymbol{\eta}}^s_q) w_q^{\eta^s},
\end{equation}
where $w_q^{\eta^s}$ are weights for the quadrature points. 

\subsection{The error of the low-dimensional solution, $\tilde{u}^{A_s}(x,\tilde{\boldsymbol{\eta}}^s)$} 
\label{RK:sec:baerror}
In the basis adaptation method, if the dimension of the new basis in $\boldsymbol{\eta}^s$ and that of the original basis in $\boldsymbol{\xi}$ are equal, there is no dimension reduction, and the accuracy of the solution $\tilde{u}^{A_s}(x,\boldsymbol{\eta}^s)$ in terms of $\boldsymbol{\eta}^s$ is the same as that of the solution $u(x,\boldsymbol{\xi}).$ This means: 
\begin{align}\label{RK:eq:error}
	\epsilon(x, \omega) = u(x,\boldsymbol{\xi}(\omega)) - \tilde{u}^{A_s}(x,\boldsymbol{\eta}^s(\omega)) = 0, \quad \text{if } r=d.
\end{align} 

For  $r \ll d$,
\begin{align}\label{RK:eq:error2}
	u(x,\boldsymbol{\xi}(\omega)) &= \tilde{u}^{A_s}_0(x) + \sum_{\mathcal{I}_1} \tilde{u}_{\mathcal{I}_1} \psi_{\mathcal{I}_1}(\tilde{\boldsymbol{\eta}}^s)  + \sum_{\mathcal{I}_{12}} \tilde{u}_{\mathcal{I}_{12}} \psi_{\mathcal{I}_{12}}(\tilde{\boldsymbol{\eta}}^s, \hat{\boldsymbol{\eta}}^s) + \sum_{\mathcal{I}_{2}} \tilde{u}_{\mathcal{I}_{2}} \psi_{\mathcal{I}_{2}}(\hat{\boldsymbol{\eta}}^s), \nonumber \\
	u(x,\boldsymbol{\xi}(\omega)) &- \left ( \tilde{u}^{A_s}_0(x) + \sum_{\mathcal{I}_1} \tilde{u}_{\mathcal{I}_1} \psi_{\mathcal{I}_1}(\tilde{\boldsymbol{\eta}}^s)\right) =   \sum_{\mathcal{I}_{12}} \tilde{u}_{\mathcal{I}_{12}} \psi_{\mathcal{I}_{12}}(\tilde{\boldsymbol{\eta}}^s, \hat{\boldsymbol{\eta}}^s) + \sum_{\mathcal{I}_2} \tilde{u}_{\mathcal{I}_2} \psi_{\mathcal{I}_2}(\hat{\boldsymbol{\eta}}^s), \nonumber \\	
	u(x,\boldsymbol{\xi}(\omega)) &- \tilde{u}^{A_s}(x,\tilde{\boldsymbol{\eta}}^s(\omega)) =   \sum_{\mathcal{I}_{12}} \tilde{u}_{\mathcal{I}_{12}} \psi_{\mathcal{I}_{12}}(\tilde{\boldsymbol{\eta}}^s, \hat{\boldsymbol{\eta}}^s) + \sum_{\mathcal{I}_2} \tilde{u}_{\mathcal{I}_2} \psi_{\mathcal{I}_2}(\hat{\boldsymbol{\eta}}^s), \nonumber \\
	\epsilon(x, \omega) &=  \sum_{\mathcal{I}_{12}} \tilde{u}_{\mathcal{I}_{12}} \psi_{\mathcal{I}_{12}}(\tilde{\boldsymbol{\eta}}^s, \hat{\boldsymbol{\eta}}^s) + \sum_{\mathcal{I}_2} \tilde{u}_{\mathcal{I}_2} \psi_{\mathcal{I}_2}(\hat{\boldsymbol{\eta}}^s), 
\end{align} 
where $\tilde{\boldsymbol{\eta}}^s = \{\eta^s_1, \ldots,\eta^s_r \};$ $\hat{\boldsymbol{\eta}}^s = \{\eta^s_{r+1}, \ldots,\eta^s_d \};$ and multi-indices $\mathcal{I}_1$ and $\mathcal{I}_2$ correspond to PC expansion terms in $\tilde{\boldsymbol{\eta}}^s$ and $\hat{\boldsymbol{\eta}}^s,$ respectively; and multi-index  $\mathcal{I}_{12}$ corresponds to mixed PC terms in $\tilde{\boldsymbol{\eta}}^s$ and $\hat{\boldsymbol{\eta}}^s.$

\section{Solution of individual boundary value problems for the subdomains} \label{RK:sec:NN_section}

We solve separate boundary value problems for each subdomain. In other words, we are interested in solving the boundary value problem defined in Eq. \eqref{RK:eq:spdeop_sub} as a collection of boundary value problems: one in each subdomain $D_s,$ where $s=1,\ldots,N_D.$ The proposed basis adaptation is applicable to any type of coupling relation between different subdomains. For linear problems, the interface conditions can be imposed through the superposition principle. However, an iterative method is usually required for a nonlinear problem. Here, we show numerical results for three different examples. We discuss in detail one iterative algorithm and provide detailed computational cost analysis for one example.

Various types of iterative algorithms, including N-N, Dirichlet-Neumann (D-N), and Dirichlet-Dirichlet (D-D) (also known as finite element tearing and interconnect, or FETI) \cite{RK:Toselli2005} can be used to solve \eqref{RK:eq:spdeop_sub} in non-overlapping subdomains. In the current work, we employ the N-N algorithm. Let $\Gamma = \left\{( \cup_{s=1}^{N_D} \partial D_s )\setminus\partial D \right\}$ be the interface corresponding to all of the subdomains and $\Gamma_s = \Gamma \cap \partial D_s$ be the interface for subdomain $D_s.$ Given the interface solution $u^0_{\Gamma},$ the N-N algorithm works as follows: 

\begin{itemize}
	\item Solve in each subdomain $D_s$ a Dirichlet problem with data $u^0_{\Gamma_s}$, obtaining solution $u^{1/2}_s$ in $D_s$, where $u^0_{\Gamma_s} = u^0_{\Gamma}\mathbb{I}_{D_s}.$
	\item Solve in each subdomain $D_s$ with Neumann data on $\Gamma_s,$  chosen as the difference of normal derivatives of $u^{1/2}_s$ obtained from solution of Dirichlet problems.
	\item Update initial $u^0_{\Gamma}$ to obtain $u^1_{\Gamma}$ using Neumann solutions on $\Gamma.$
\end{itemize}

\begin{equation} \label{RK:dirichlet}
(Dirichlet)
\begin{cases}
\mathcal{L}(x,u^{s,n+1/2}(x,\tilde{\boldsymbol{\eta}}^s);a(x,{\boldsymbol{\xi}})) &= f(x,{\boldsymbol{\xi}})  \;\; \rm{in}~D_s\times \Omega, \nonumber \\
 \mathcal{B}(x,u^{s,n+1/2}(x,\tilde{\boldsymbol{\eta}}^s);a(x,{\boldsymbol{\xi}})) &= h(x,{\boldsymbol{\xi}})  \;\; \rm{on}~\partial D_s \times \Omega,
\end{cases}
\end{equation}

\begin{equation}\label{RK:neumann}
(Neumann)
\begin{cases}
\mathcal{L}(x,u^{s,n+1}(x,\tilde{\boldsymbol{\eta}}^s);a(x,{\boldsymbol{\xi}}) &= 0  \;\; \rm{in}~D_s\times \Omega, \nonumber \\
 \mathcal{C}(x,\nabla u^{s,n+1}, \nabla u^{s^{\prime},n+1/2};a(x,{\boldsymbol{\xi}})) &= 0  \;\; \rm{in}~(\partial D_s \cap \partial D_{s^{\prime}}) \times \Omega, \nonumber \\
 \mathcal{B}(x,u^{s,n+1}(x,\tilde{\boldsymbol{\eta}}^s);a(x,{\boldsymbol{\xi}})) &= h(x,{\boldsymbol{\xi}})  \;\; \rm{on}~(\partial D_s \cap \partial D)\times \Omega,
\end{cases}
\end{equation}
$u_{\Gamma}^{n+1} = u_{\Gamma}^{n}-\theta \sum_s^{N_D}u_{\Gamma}^{s,n+1}$ on $\Gamma$ for some $\theta \in (0, \theta_{\max}).$ The N-N algorithm is equivalent to a preconditioned Richardson method for the Schur complement system \cite{RK:Toselli2005}. 
In the following section, we use a non-iterative version of the N-N algorithm \cite{RK:Toselli2005}. 

\subsection{Non-iterative Neumann-Neumann algorithm: Deterministic case} \label{RK:sec:NN}
Let a deterministic PDE on the domain $D$ be discretized to obtain the system of algebraic equations $K_f(u)= f.$ These nonlinear equations are solved using an iterative scheme, such as Newton's method, where they are linearized to obtain a system of linear equations, $K u = f.$  For each subdomain $D_s,$ we can write the linear equations in matrix form as

\begin{equation}\label{RK:eqn:lineq}
\begin{pmatrix}
	K^s_{II} & K^s_{I\Gamma} \\
	K^s_{\Gamma I} & K^s_{\Gamma\Gamma}
\end{pmatrix}
\begin{pmatrix}
	u^{s}_{I}  \\
	u^{s}_{\Gamma} 
\end{pmatrix}
= 
\begin{pmatrix}
	f^s_I  \\
	f^s_{\Gamma}
\end{pmatrix},
\end{equation}
where the subscript $I$ corresponds to interior degrees of freedom and $\Gamma$ corresponds to the interface for each subdomain $D_s.$ To solve the system of equations~\eqref{RK:eqn:lineq} using the N-N algorithm, we choose an initial guess $u_{\Gamma}^0$  for the solution on the interface and solve the following Dirichlet and Neumann problems: 
\begin{align} \label{RK:eqn:neumann}
&(D_i) \quad
K^s_{II} u^{s,n+1/2}_{I} + K^s_{I\Gamma} u^{s,n+1/2}_{\Gamma} = f^s_{\Gamma} \quad \text{ in } D_s, \nonumber \\
&(N_i) \quad
\begin{pmatrix}
	K^s_{II} & K^s_{I\Gamma} \\
	K^s_{\Gamma I} & K^s_{\Gamma\Gamma}
\end{pmatrix}
\begin{pmatrix}
	u^{s,n+1}_{I}  \\
	u^{s,n+1}_{\Gamma} 
\end{pmatrix}
= 
\begin{pmatrix}
	0  \\
	r_{\Gamma}
\end{pmatrix}\quad \text{ in } D_s, \\
	&u^{n+1}_{\Gamma} = u^{n}_{\Gamma} - \theta \sum_s^{N_D} u_{\Gamma}^{s,n+1}, \nonumber
\end{align}
where the residual $r_{\Gamma}$ is defined as 
\begin{equation} \label{RK:eqn:r_gamma}
	r_{\Gamma} = \sum_{s=1}^{N_D} R_s^T  (K^s_{\Gamma I}u^{s,n+1/2}_{I} + K^s_{\Gamma \Gamma} u^n_{\Gamma} -f^s_{\Gamma}).
\end{equation}
Given the interface solution $u_{\Gamma}$ in the whole domain $D,$ $R_s u_{\Gamma}$ is its restriction on the interface of the subdomain $D_s$. Here, $R_s$ is a rectangular matrix whose elements are zeros and ones and that acts as a restriction operator. Its transpose acts as a  prolongation operator. The  Neumann problem in~\eqref{RK:eqn:neumann} solves for {\it all} (interior and interface) degrees of freedom $u^s$ in the subdomain $D_s.$ However, it is efficient to solve initially for the interface degrees $u_{\Gamma}^s$ then for the interior degrees. 

Let $S^s_{\Gamma} = K^s_{\Gamma\Gamma} - K^s_{\Gamma I} (K^s_{II})^{-1}K^s_{I\Gamma} $ be the Schur complement matrix and $g^s_{\Gamma} = f^s_{\Gamma} -K^s_{\Gamma I} (K^s_{II})^{-1}f^s_I $ in subdomain $D_s$. The global Schur complement matrix is computed as
\begin{align} \label{global_schur}
	S_{\Gamma} &= \sum_{s=1}^{N_D} R_s^T S^s_{\Gamma}  R_s, \\
	g_{\Gamma} &= \sum_{s=1}^{N_D} R_s^T g^s_{\Gamma}. \label{global_schur2}
\end{align}
In~\eqref{RK:eqn:neumann}, we can eliminate $ u^{s,n+1/2}_{I}$ and $u^{s,n+1}_{I}$ to obtain
\begin{equation} \label{RK:eqn:r_gamma2}
	r_{\Gamma} = -(g_{\Gamma} - S_{\Gamma} u_{\Gamma}^n), 
\end{equation}
which means that the total difference in fluxes $r_{\Gamma}$ is equal to minus the residual of the Schur complement system. Using block factorization,
\begin{equation} \label{RK:eqn:u_s_gamma}
	u^{s,n+1}_{\Gamma}  = (S^s_{\Gamma})^{-1} r_{\Gamma} = -(S^s_{\Gamma})^{-1}(g_{\Gamma} - S_{\Gamma} u_{\Gamma}^n), 
\end{equation}
the interface solution can be written as
\begin{equation} \label{RK:eqn:u_gamma}
	u^{n+1}_{\Gamma} = u^{n}_{\Gamma} + \theta \left( \sum_s^{N_D} (S^s_{\Gamma})^{-1} \right ) (g_{\Gamma} - S_{\Gamma} u_{\Gamma}^n).
\end{equation}
Solving the Dirichlet problem for the interior solution $u_{I}$ and the Schur complement system in \eqref{RK:eqn:u_s_gamma} and \eqref{RK:eqn:u_gamma} iteratively gives the solution of the linear system of equations $K u = f.$

When the interface degrees of freedom are few, the size of the global Schur complement $S_{\Gamma}$ is not very large. 
It is possible to have a non-iterative scheme by formulating the solution in terms of the global Schur complement matrix \cite{RK:Toselli2005}.  This is done by computing $S_{\Gamma}$ and $g_{\Gamma}$ from \eqref{global_schur}-\eqref{global_schur2} and solving
\begin{equation} \label{schur_solve}
	S_{\Gamma} u_{\Gamma} = g_{\Gamma} \quad \text{on} \quad \Gamma
\end{equation}
to obtain the solution for the interface degrees of freedom on $\Gamma$. Then, we compute the solution for the interior degrees of freedom, $u^{s}_{I}$, as
\begin{equation}\label{eqn:u_s_I}
	u^{s}_{I} = (K^s_{II})^{-1}[f^s_I - K^s_{I\Gamma} R_s u_{\Gamma}], \quad s=1,\ldots,N_D. 
\end{equation}
Note that $S_{\Gamma}$ and $g_{\Gamma}$ depend only on the stiffness matrix $K$ and the right-hand side $f$, not on the solution $u.$ Thus, they do not need to be updated in an iterative fashion. 

In the next section, we extend this method for a stochastic PDE and use stochastic basis adaptation in each subdomain to reduce the computational cost. 

\subsection{Non-iterative Neumann-Neumann algorithm: Stochastic case with basis adaptation}
For the subdomains $D_s, s\in N_D$, we first write a linear system of equations in $\tilde{\boldsymbol{\eta}}^s$ as follows:
\begin{equation}
\begin{pmatrix}
	K^s_{II}(\tilde{\boldsymbol{\eta}}^s) & K^s_{I\Gamma}(\tilde{\boldsymbol{\eta}}^s) \\
	K^s_{\Gamma I} (\tilde{\boldsymbol{\eta}}^s)& K^s_{\Gamma\Gamma}(\tilde{\boldsymbol{\eta}}^s)
\end{pmatrix}
\begin{pmatrix}
	u^{s}_{I}(\tilde{\boldsymbol{\eta}}^s)  \\
	u^{s}_{\Gamma} (\tilde{\boldsymbol{\eta}}^s)
\end{pmatrix}
= 
\begin{pmatrix}
	f^s_I  \\
	f^s_{\Gamma}
\end{pmatrix}.
\end{equation}
Then, the Schur complement $S_{\Gamma}^s(\boldsymbol{\eta}^s)$ and $g^s_{\Gamma}(\boldsymbol{\eta}^s)$ are computed as
\begin{align}\label{Rk:eqn:schur_eta}
	S_{\Gamma}^s(\tilde{\boldsymbol{\eta}}^s) &= K^s_{\Gamma\Gamma}(\tilde{\boldsymbol{\eta}}^s) - K^s_{\Gamma I}(\tilde{\boldsymbol{\eta}}^s) [K^s_{II}(\tilde{\boldsymbol{\eta}}^s)]^{-1}K^s_{I\Gamma}(\tilde{\boldsymbol{\eta}}^s),  \nonumber\\ 
	g^s_{\Gamma}(\tilde{\boldsymbol{\eta}}^s) &= f^s_{\Gamma} -K^s_{\Gamma I}(\tilde{\boldsymbol{\eta}}^s) [K^s_{II}(\tilde{\boldsymbol{\eta}}^s)]^{-1}f^s_I.
\end{align}
For each subdomain $D_s,$ we compute the global Schur complement $S_{\Gamma}(\tilde{\boldsymbol{\eta}}^s)$ and the corresponding interface solution $u_{\Gamma}(\tilde{\boldsymbol{\eta}}^s).$ Then, we compute the solution for the interior points $u^{s}_{I}(\tilde{\boldsymbol{\eta}}^s)$ using \eqref{eqn:u_s_I}.
Algorithm~\ref{RK:alg:NN} describes the implementation of the proposed N-N method. 
\begin{algorithm}[H]
\caption{Neumann-Neumann algorithm} \label{RK:alg:NN} 
  1. Choose collocation points $\tilde{\boldsymbol{\eta}}_q^s$ for $q=1,\ldots, Q_{\eta^s}$ and $s=1, \ldots, N_D$\\
  2. Compute $S_{\Gamma}^s(\tilde{\boldsymbol{\eta}}_q^s)$ and $g_{\Gamma}^s(\tilde{\boldsymbol{\eta}}_q^s)$ using \eqref{Rk:eqn:schur_eta}\\
  3. {\bf for} $ s=1 \ldots N_D$ {\bf do} \\
  4.  $\quad$  {\bf for} $ s'=1 \ldots N_D$, $s'\neq s$ {\bf do} \\  
  5. $\quad$ $\quad$ Project  $S_{\Gamma}^{s'}(\tilde{\boldsymbol{\eta}}^{s'})$ and $g_{\Gamma}^{s'}(\tilde{\boldsymbol{\eta}}^{s'})$ on to $\tilde{\boldsymbol{\eta}}^s$ space using algorithm \ref{RK:alg:projection}\\  
  6. $\quad$ $\quad$ Obtain $S_{\Gamma}^{s'}(\tilde{\boldsymbol{\eta}}^{s})$ and $g_{\Gamma}^{s'}(\tilde{\boldsymbol{\eta}}^{s})$ \\  
  7. $\quad$ $\quad$ Compute  $S_{\Gamma}^{s'}(\tilde{\boldsymbol{\eta}}_q^{s})$ and $g_{\Gamma}^{s'}(\tilde{\boldsymbol{\eta}}_q^{s})$ for $q=1,\ldots, Q_{\eta^s}$\\  
  8. $\quad$ {\bf end for}\\
  9. $\quad$ $S_{\Gamma}(\tilde{\boldsymbol{\eta}}_q^{s}) = \sum_{s'=1}^{N_D} R_{s'}^T S_{\Gamma}^{s'}(\tilde{\boldsymbol{\eta}}_q^{s})  R_{s'}$\\   
  10.$\quad$ $g_{\Gamma}(\tilde{\boldsymbol{\eta}}_q^{s}) = \sum_{s'=1}^{N_D} R_{s'}^T g_{\Gamma}^{s'}(\tilde{\boldsymbol{\eta}}_q^{s})$\\  
  11.$\quad$ $u_{\Gamma}(\tilde{\boldsymbol{\eta}}_q^{s}) =  [S_{\Gamma}(\tilde{\boldsymbol{\eta}}_q^{s}) ]^{-1}g_{\Gamma}(\tilde{\boldsymbol{\eta}}_q^{s}),$ for $q=1,\ldots, Q_{\eta^s}$\\
  12.$\quad$ $u^s_{I}(\tilde{\boldsymbol{\eta}}_q^{s}) =   (K^s_{II}(\tilde{\boldsymbol{\eta}}_q^{s}))^{-1}[f^s_I - K^s_{I\Gamma}(\tilde{\boldsymbol{\eta}}_q^{s}) R_s u_{\Gamma}(\tilde{\boldsymbol{\eta}}_q^{s})]$\\
  13. {\bf end for}.
\end{algorithm}

To compute the solution $u^s(\tilde{\boldsymbol{\eta}}^{s})$ in each subdomain $D_s$ (steps 9-12 in Algorithm~\ref{RK:alg:NN}), we need to project the local Schur complement $S_{\Gamma}^{s'}(\tilde{\boldsymbol{\eta}}^{s'})$ and vector $g_{\Gamma}^{s'}(\tilde{\boldsymbol{\eta}}^{s'})$ on to the PC basis in random variables $\tilde{\boldsymbol{\eta}}^{s}$, for $s'=1 \ldots N_D$, $s'\neq s.$  In Algorithm~\ref{RK:alg:NN}, $N_D \times (N_D-1)$ number of projections are required for $S_{\Gamma}^{s'}(\tilde{\boldsymbol{\eta}}^{s'})$ and $g_{\Gamma}^{s'}(\tilde{\boldsymbol{\eta}}^{s'}).$ As the number of subdomains $N_D$ increases, the computational cost for an interface solution becomes a bottleneck because of the large number ($N_D\times (N_D-1)$) of projections involved. 
Algorithm 2 illustrates the procedure for projecting a function $y(\tilde{\boldsymbol{\eta}}^{s'})$ on the $\tilde{\boldsymbol{\eta}}^{s}$ space and computing the function $y$ at $\tilde{\boldsymbol{\eta}}_q^{s}$.

\begin{algorithm}[H]
\caption{Projection} \label{RK:alg:projection} 
  1. Choose collocation points $\tilde{\boldsymbol{\eta}}_q^s$ for $q=1,\ldots, Q$\\
  2. $\tilde{\boldsymbol{\eta}}_q^{s'} = A_{s'} A_s^{-1}\tilde{\boldsymbol{\eta}}_q^s$, [since $\tilde{\boldsymbol{\eta}}^{s'} = A_{s'}\boldsymbol{\xi}$ and $\tilde{\boldsymbol{\eta}}^{s} = A_{s} \boldsymbol{\xi}]$\\
  3. $y(\tilde{\boldsymbol{\eta}}^{s'}) = y(\tilde{\boldsymbol{\eta}}^{s}) = \sum_{i=1}^{N^s_{\eta}} y_i
 \psi_i(\tilde{\boldsymbol{\eta}}^{s})$\\
  4. $y_i = \langle y(\tilde{\boldsymbol{\eta}}^{s})   \psi_i(\tilde{\boldsymbol{\eta}}^{s}) \rangle = \sum_{q=1}^Q,   y[\tilde{\boldsymbol{\eta}}^{s'}(\tilde{\boldsymbol{\eta}}_q^{s})] \psi_i(\tilde{\boldsymbol{\eta}}_q^{s}) w_q, \quad i = 1, \ldots,N^s_{\eta}$\\
  5. $ y(\tilde{\boldsymbol{\eta}}_q^{s}) = \sum_{i=1}^{N^s_{\eta}} y_i
 \psi_i(\tilde{\boldsymbol{\eta}}_q^{s})$ for $q=1,\ldots, Q$\\
  
\end{algorithm}

\subsection{Computational cost}\label{comp_cost}
In this section, we define metrics for comparing the computational cost of solving stochastic PDE using PC with and without domain decomposition and basis adaptation. Without loss of generality, we assume that the number of floating point operations (flops) for solving a system of linear algebraic equations with $n$ unknown variables is approximately $ (\frac{2}{3}n^3+2n^2),$ as it is for the lower-upper (LU) factorization method. We use this metric (cost of the LU factorization method for solving a linear system) to compare the total computational cost for the two cases.

\subsubsection{Non-intrusive methods}\label{Nonintru_comp_cost}
Let $N^l_{\xi}$ be the number of sparse-grid collocation points for a given dimension $d$ and sparse-grid level $l$ \textcolor{red}{}. For the original $\boldsymbol{\xi}$ random space, we solve a deterministic system of equations for each collocation point in that random space. Hence, the computational cost for solving a $d$-dimensional stochastic system with sparse-grid level $l$ is $\approx N^l_{\xi} \left (\frac{2}{3} n^3+2 n^2\right )$. For each subdomain, we solve for a low-dimensional solution $\tilde{u}_I^{A_s}(x,\boldsymbol{\eta}^s)$ at interior points and $\tilde{u}^{A_s}_{\Gamma}(x,\boldsymbol{\eta}^s)$ at interface points. Using a $d^s_{\eta}$ dimensional space in $\boldsymbol{\eta}^s$ and sparse-grid level $l,$ the total number of collocation points for each subdomain $D_s$ is $N^l_{\eta}.$ Hence, the total computational cost is approximately equal to $N^l_{\eta} \sum_{s=1}^{N_D} \left(\frac{2}{3} (n^s_I)^3+2  (n^s_I)^2 +  \frac{2}{3} (n_{\Gamma})^3+2 (n_{\Gamma})^2 \right).$

\subsubsection{Computational cost of projection algorithm}\label{projection_cost}
In Algorithm \ref{RK:alg:NN}, the Schur complement $S^{s'}(\tilde{\boldsymbol{\eta}}^{s'})$ and the vector $g_{\Gamma}^{s'}(\tilde{\boldsymbol{\eta}}^{s'})$ for $s'=1,\ldots,N_D$ are projected onto the space spanned by the random variables $\tilde{\boldsymbol{\eta}}^{s}.$ The projection method is implemented as Algorithm \ref{RK:alg:projection}. Here, we describe the computational cost of Algorithm 2 in terms of flops. We assume the addition or subtraction of two vectors of size $m$ needs $m$ flops, and multiplication of a vector of size $m$ with a scalar requires $m$ flops. If $N^l_{\eta_s}$ is the number of collocation points in $\tilde{\boldsymbol{\eta}}^{s}$, $N^s_{\eta}$ is the total number of PC coefficients, $n_{\Gamma} \times n_{\Gamma}$ is the size of the matrix $S^{s'}(\tilde{\boldsymbol{\eta}}^{s'})$, and $n_{\Gamma} \times 1$ is the size of the vector $g^{s'}(\tilde{\boldsymbol{\eta}}^{s'})$, then the computational cost of the projection of the Schur complement $S^{s'}(\tilde{\boldsymbol{\eta}}^{s'})$ is $\{(N_D-1) N^l_{\eta_s} N^s_{\eta} (2n_{\Gamma}^2-1)+N_D N^l_{\eta_s} n_{\Gamma}^2\},$ while that of the vector $g_{\Gamma}^{s'}(\tilde{\boldsymbol{\eta}}^{s'})$ is $\{(N_D-1) N^l_{\eta_s} N^s_{\eta} (2n_{\Gamma}-1)+N_D N^l_{\eta_s} n_{\Gamma}\}.$ Hence, to obtain the total cost of these projections in Algorithm~\ref{RK:alg:NN}, we need to multiply this cost with the total number of subdomains $N_D, $ that is:
\begin{align} \label{comp_cost_estimate}
	\text{Cost in flops} &= N_D \{(N_D-1) N^l_{\eta_s} N^s_{\eta} (2n_{\Gamma}^2-1)+N_D N^l_{\eta_s} n_{\Gamma}^2\} \nonumber \\
	&  + N_D \{(N_D-1) N^l_{\eta_s} N^s_{\eta} (2n_{\Gamma}-1)+N_D N^l_{\eta_s} n_{\Gamma}\}.
\end{align}

\section{Numerical examples} \label{RK:sec:numerical}
Here, we apply our methodology to three examples: i) a saturation-based one-dimensional linear Richards equation, ii) a pressure-based one-dimensional nonlinear Richards equations, and iii) a two-dimensional steady-state diffusion equation with a point sink at the center of the computational domain.

\subsection{One-dimensional linear Richards equation} \label{RK:sec:1d-lin-rich}

In this section, we compute the solution of the stochastic one-dimensional saturation-based Richards equation describing vertical unsaturated flow in heterogeneous porous media \cite{Tartakovsky2008}: 
\begin{align} \label{RK:eqn:rich_lin}
	 \frac{d}{dz}\left [ K(z, \theta) \frac{d\psi}{d\theta} \frac{d\theta}{d z} \right ] +  \frac{d K(z, \theta)}{d z} = 0, \quad 0 < z < L. 
\end{align}
This equation is subject to the Dirichlet boundary condition at the bottom of the domain and  Neumann boundary condition at the top:

\begin{align} \label{RK:eqn:rich_lin_BC}
	 \theta(0) = \Theta_0 \quad \text{and}  \quad \left ( K(z, \theta) \frac{d\psi}{d\theta} \frac{d\theta}{dz} \right ) _{z=L} = -q.
\end{align}
We assume the Gardner-Russo exponential model \cite{Tartakovsky2008} for unsaturated hydraulic conductivity as a function of pressure head, $\psi(z)$, and the saturation $\theta$ as

\begin{align} \label{RK:eqn:rich_lin_K}
	K &= K_s \exp(\alpha \psi) \\
	\theta &= \theta_s \exp(\alpha \psi) ,
\end{align}
where $K_s$ is the saturated hydraulic conductivity, $\theta_s$ is the complete saturation, and $\alpha$ is the Gardner parameter. In this model, we introduce uncertainty through the saturated hydraulic conductivity by treating it as a random field with log-normal distribution. For the numerical experiment we use a two-layer soil that has different values for $\alpha$ and $K_s$ in each layer. The solution is computed using a non-intrusive stochastic collocation method with dimension $15$ with sparse-grid level $5$ that requires $39941$ collocation points (deterministic solutions). 

\begin{figure}[t!]
    \centering
    \begin{subfigure}[t]{0.45\textwidth}
        \centering
        \includegraphics[scale=.3]{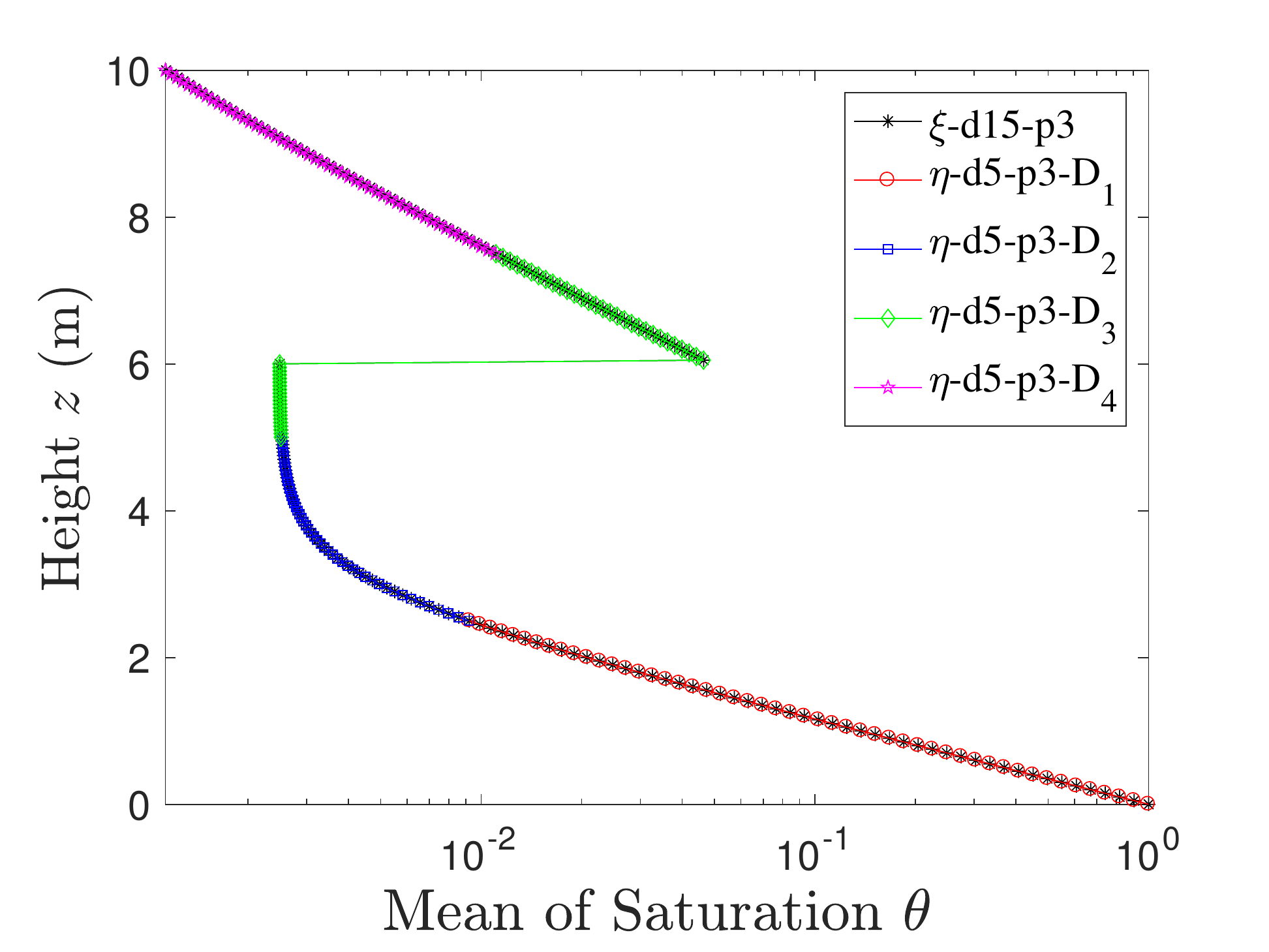}
        \caption{Mean} \label{RK:fig:theta_mean_xid15_etad5}
    \end{subfigure}        
    \begin{subfigure}[t]{0.45\textwidth}
        \centering
        \includegraphics[scale=.3]{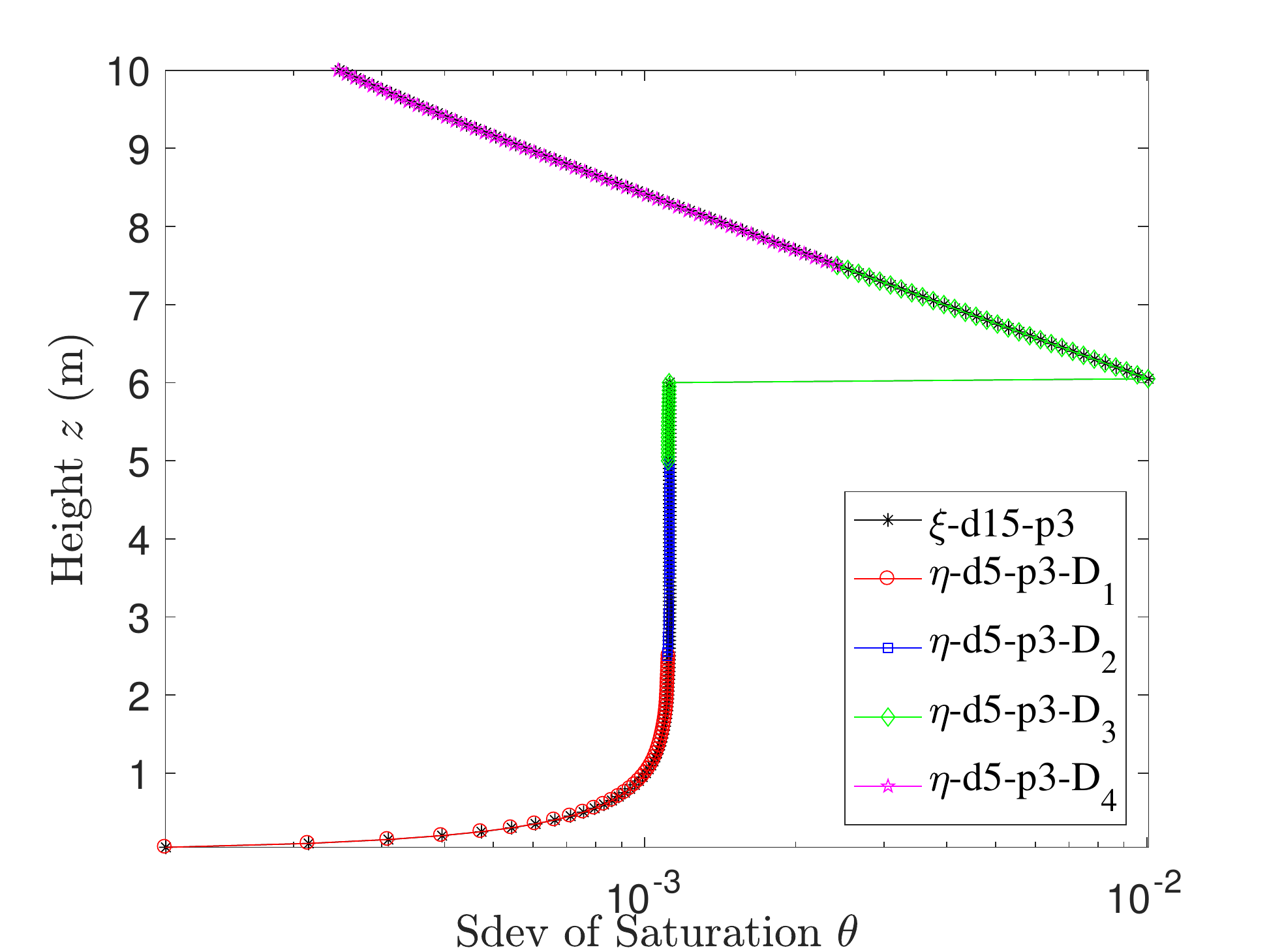}
        \caption{Sdev} \label{RK:fig:theta_sdev_xid15_etad5}
    \end{subfigure}    
       \caption{(a) Mean and (b) standard deviation (Sdev) of the saturation computed without domain decomposition and basis adaptation (dimension, d = 15; sparse-grid level, l = 5; order, p = 3) shown in `solid black'; and with basis adaptation (dimension, r = 5; sparse-grid level, l = 5; order, p = 3) and domain decomposition method ($|D|=4$) shown in `red', `blue', `green' and `magenta' in subdomains $D_1, D_2, D_3 $and $D_4$ respectively}. \label{RK:fig:mean_sdev_lin_rich_theta}

\end{figure}

To implement the proposed approach, we decompose the domain into four non-overlapping subdomains and computed the solution in each subdomain in a new set of random variables with dimension $5$ and sparse-grid level $5.$ For this computation, we require $781$ deterministic solutions for each subdomain. As such, for the four subdomains, we need $3124$ deterministic solutions in a smaller computational domain, $L_d = L/4.$ Here, the computational cost is quite small compared to the full solution because the number of deterministic solutions is much smaller ($3124 \ll 39941$), and the size of each deterministic solution is also four times smaller. Our aim here is to show that the method works even when there is a discontinuity of the material property (hydraulic conductivity) in the physical domain. In the first layer ($0< z < 6.0 \text{ m}$), the mean of $K_s$ is assumed to be $1.0 \text{ m} \text{d}^{-1},$ and $\alpha = 2.0 \text{ m}^{-1}.$ In the second layer ($6.0 < z < 10.0 \text{ m}$), the mean of $K_s=10.0 \text{ m} \text{d}^{-1},$ and $\alpha = 1.0 \text{ m}^{-1}.$ 
Fig. \ref{RK:fig:mean_sdev_lin_rich_theta} shows the discontinuity in the mean and standard deviation of the saturation $\theta$.
Results corresponding to different subdomains are plotted with different colors. It is apparent that the subdomain $D_3$ contains the interface of the two soil layers, and the basis adaptation method still works well, even when the material properties are discontinuous in the physical domain. 

\begin{figure}[t!]
    \centering
    \begin{subfigure}[t]{0.45\textwidth}
        \centering
        \includegraphics[scale=.3]{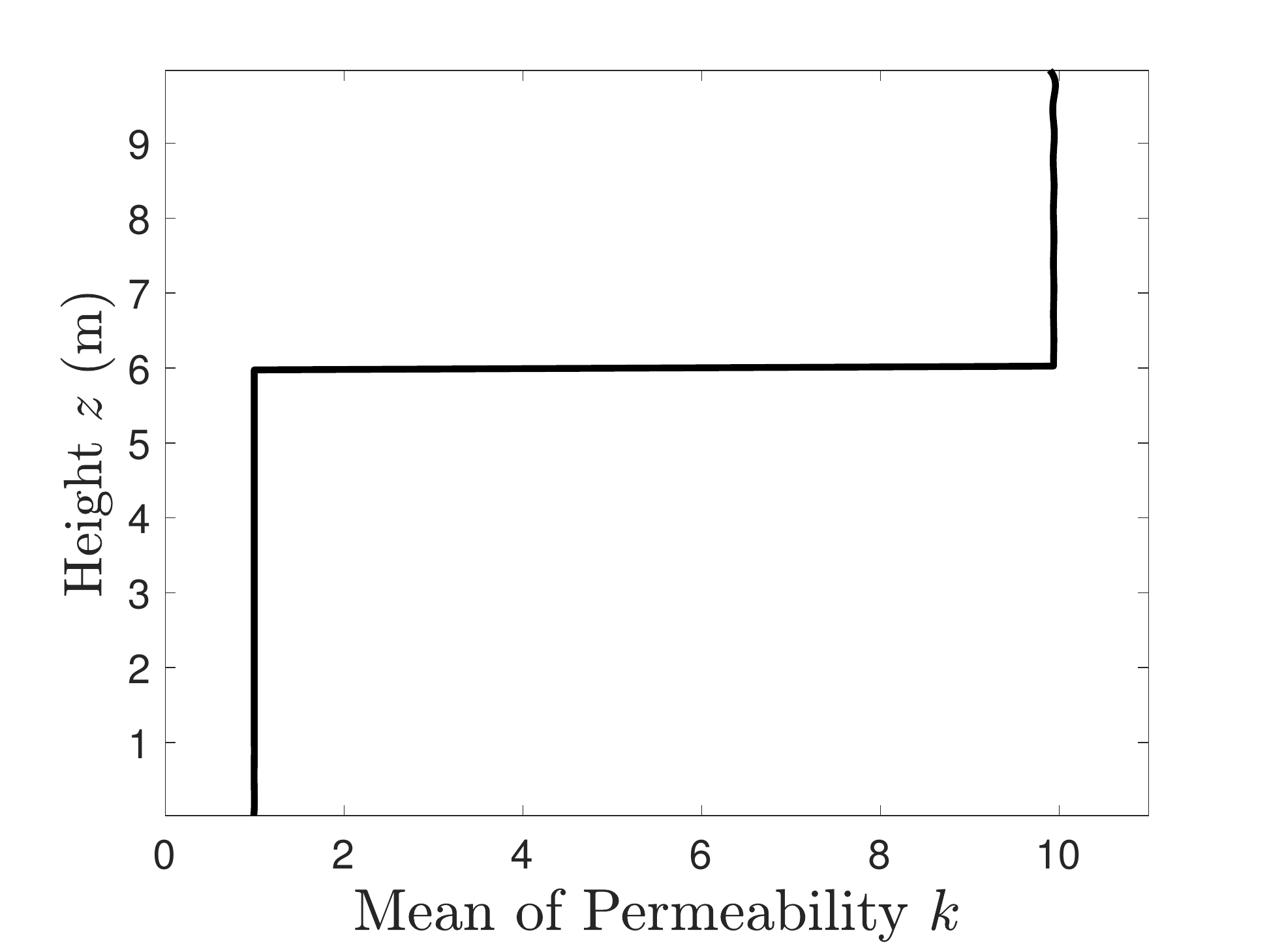}
        \caption{Mean} \label{RK:fig:k_mean_xid15_etad5}
    \end{subfigure}        
    \begin{subfigure}[t]{0.45\textwidth}
        \centering
        \includegraphics[scale=.3]{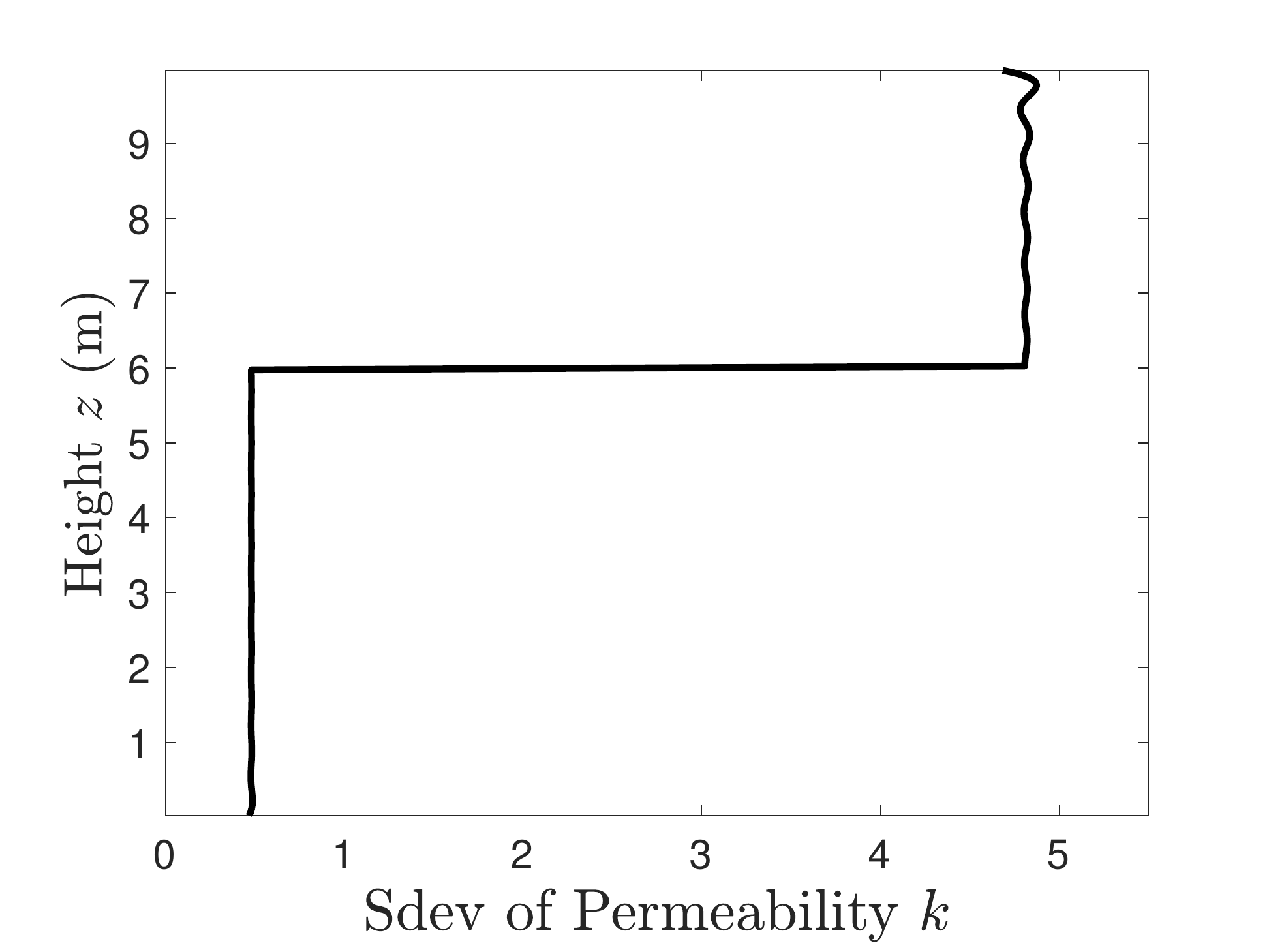}
        \caption{Sdev} \label{RK:fig:k_sdev_xid15_etad5}
    \end{subfigure}    
       \caption{(a) Mean and (b) standard deviation (Sdev) of the hydraulic conductivity $K_s~(\text{m} \text{d}^{-1})$.} \label{RK:fig:mean_sdev_lin_rich_k}
\end{figure}

\subsection{One-dimensional nonlinear Richards equation} \label{RK:sec:1d-nl-rich}
\begin{figure}[t!]
    \centering
    \begin{subfigure}[t]{0.45\textwidth}
        \centering
        \includegraphics[scale=.3]{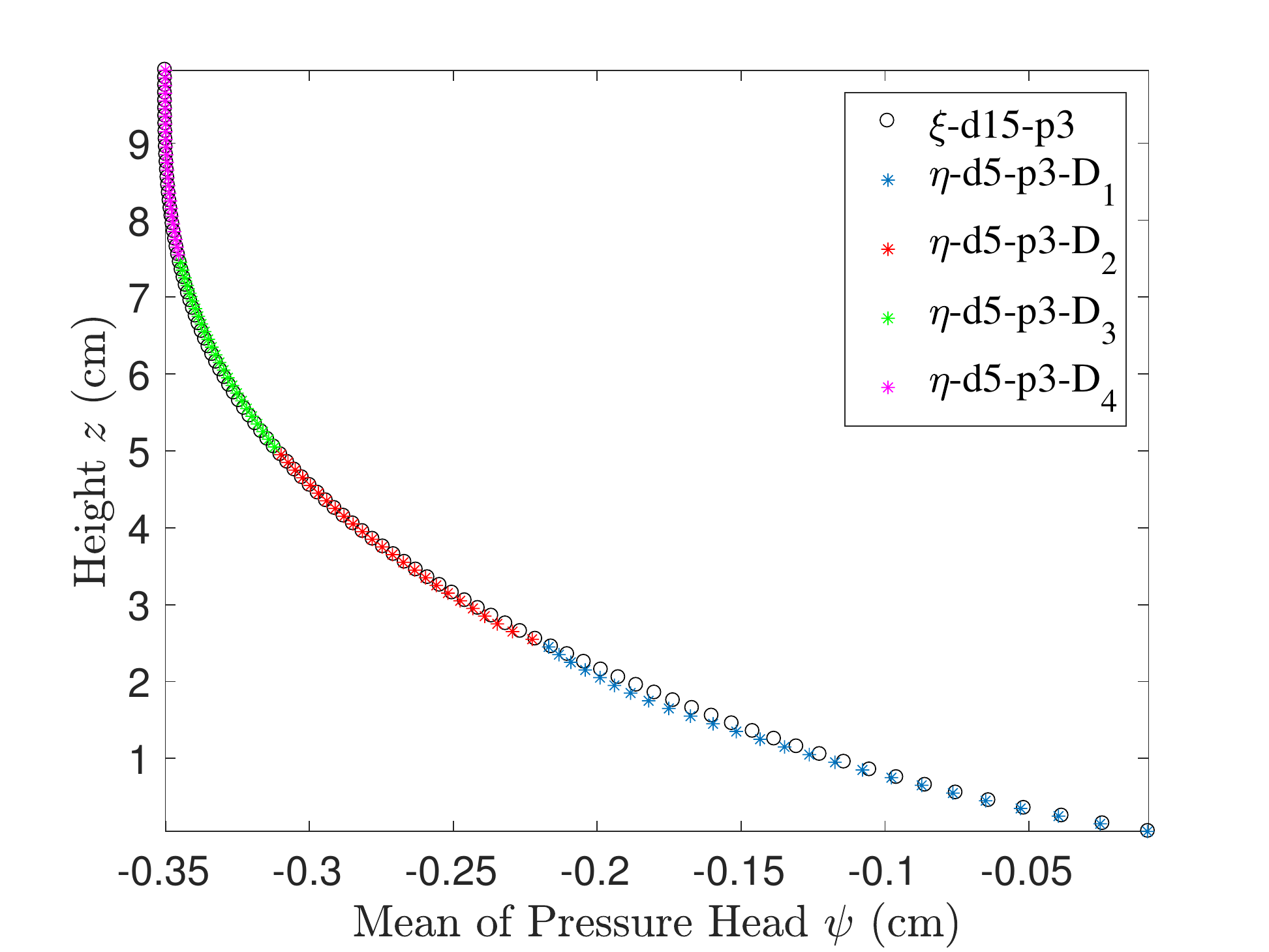}
        \caption{Mean} \label{RK:fig:y_mean_xid15_etad5}
    \end{subfigure}        
    \begin{subfigure}[t]{0.45\textwidth}
        \centering
        \includegraphics[scale=.3]{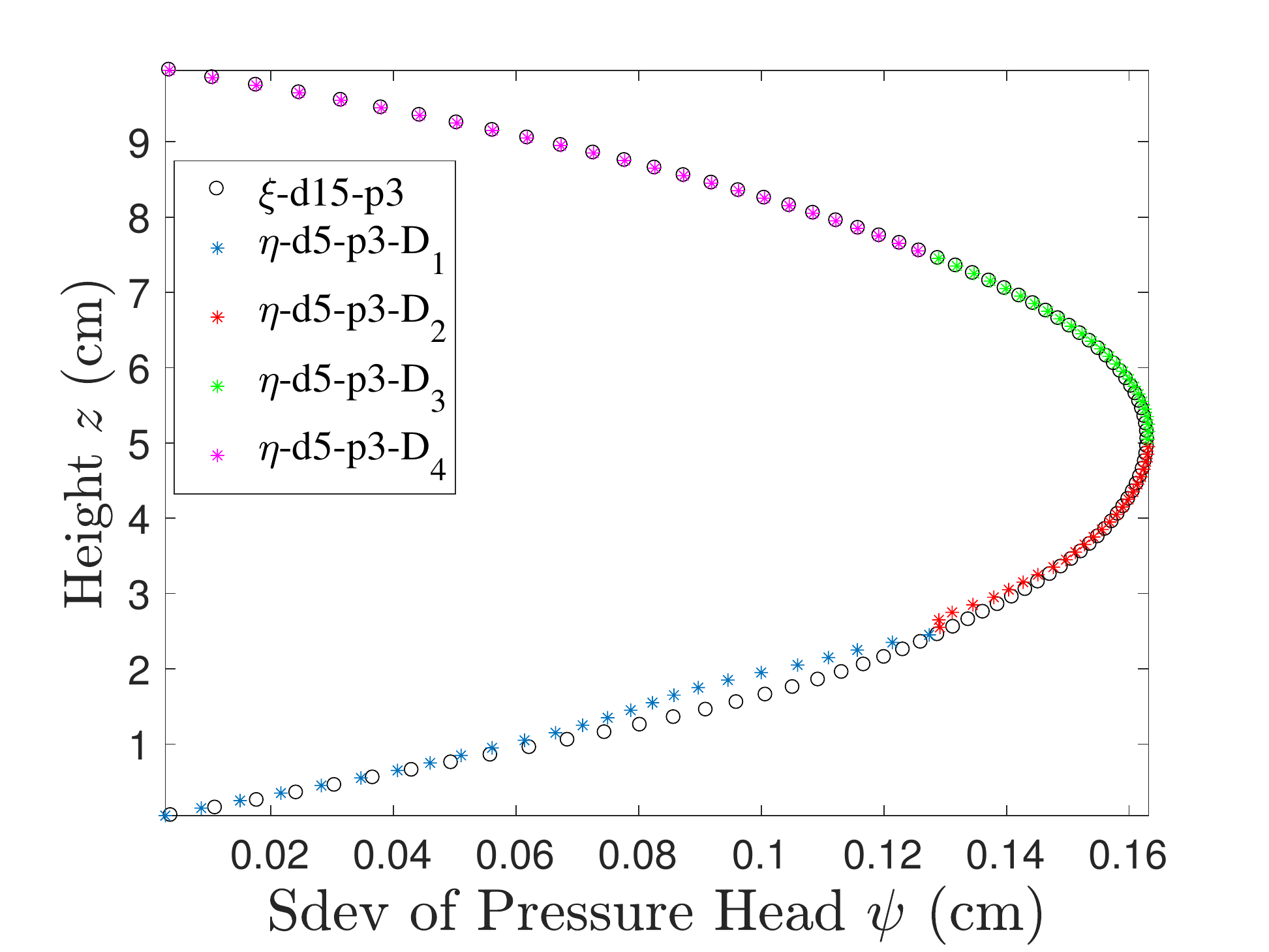}
        \caption{Sdev} \label{RK:fig:y_sdev_xid15_etad5}
    \end{subfigure}    
       \caption{(a) Mean and (b) standard deviation (Sdev) of the pressure head computed without domain decomposition and basis adaptation (dimension, d = 15; sparse-grid level, l = 5; order, p = 3) shown in `solid black'; and with basis adaptation (dimension, r = 5; sparse-grid level, l = 5; order, p = 3) and domain decomposition method ($|D|=4$) shown in `blue', `red', `green' and `magenta' in subdomains $D_1, D_2, D_3 $and $D_4$ respectively} \label{RK:fig:mead_sdev_nl_rich_y}
\end{figure}

Here, we solve a nonlinear pressure-based Richards equation with van Genuchten model \cite{Pan1995} for hydraulic conductivity $K$ and saturation $\theta.$ 
The nonlinear Richards equation subject to the Dirichlet boundary conditions is given by 

\begin{align} \label{RK:eqn:rich_nl}
	& \frac{d}{dz}\left [ K(z, \psi) \left (\frac{d\psi}{dz}+1 \right )  \right ] = 0, \quad \quad 0<z<L \\
	& \psi(0) = \psi_b \quad \text{and}  \quad \psi(L)=\psi_t,
\end{align}
where $\psi$ is the pressure.
The water retention $S_e$ and hydraulic conductivity $K$ are related to the pressure head $\psi$ through 
\begin{equation} \label{RK:eqn:Se}
	S_e = \frac{\theta-\theta_r}{\theta_s-\theta_r} = \left [\frac{1}{1+(\alpha_{vg} |\psi|)^n} \right ]^m
\end{equation}
and
\begin{equation} \label{RK:eqn:K_vg}
	K = K_s \sqrt{S_e} [1-(1-S_e^{1/m})^m]^2,
\end{equation}
where $m=1-1/n;$ $\theta_r$ and $\theta_s$ are residual and saturated water contents, respectively; $K$ is the hydraulic conductivity; $\alpha_{vg}$ and $n$ are van Genuchten model parameters obtained from the experiments. To solve this equation, we decompose the one-dimensional domain into four subdomains, and compute the nonlinear solution using the following approach: we obtain the solution iteratively, starting with an initial guess for the interface conditions then compute the solution in each subdomain.  We assume a log-normal distribution for $K_s(z,\boldsymbol{\xi}).$ For the numerical implementation, we used the soil properties from \cite{Pan1995}: $n=1.3954, \alpha_{vg}=0.0104, \theta_s=0.4686, \theta_r=0.1060$, and $K_s$ is $0.5458$ cm/h. The reference solution is obtained for a stochastic dimension of $15$ and sparse-grid level of $5$, which requires a total of $39941$ deterministic solutions. In the approximate solution with domain decomposition and basis adaptation, we obtained  the solution using reduced dimension $5$ and sparse-grid level \textcolor{green}{of 5} in each subdomain, which requires $781$ deterministic solutions for each subdomain. Hence, we need $4\times781=3124$ deterministic solutions in each step of the iterative algorithm. We stopped after $5$ iterations and hence the total number of deterministic solutions needed is $5\times3124=15620$, which is half the total number required without basis adaptation. However, we note that the size of each subdomain is 1/4th that of the original domain $D.$ Therefore, each deterministic problem in a subdomain $D_s$ is much smaller than that in the  original domain $D.$ We used a coefficient of variation $\sigma_K/\mu_K=0.1$ and the correlation length of $L/4$ where $L=10$ cm. We assumed Dirichlet boundary conditions on both sides with values, $\psi_t=-0.35$cm and $\psi_b=0.0$ cm. Fig.~\ref{RK:fig:mead_sdev_nl_rich_y} shows the (a) mean and (b) standard deviation of the pressure head computed without basis adaptation and with basis adaptation and with basis adaptation and domain decomposition. The low dimensional solution is in very good agreement with the full solution. 

\subsection{Two-dimensional steady-state diffusion} \label{RK:sec:2d-diff}

In this section, we apply the proposed approach to a two-dimensional steady-state diffusion equation with random space-dependent diffusion coefficient defined on the spatial domain $D = [0, 240] \times [0, 60],$ such that $x=(x_1,x_2) \in D$. We use Dirichlet boundary conditions at the boundaries perpendicular to the  $x_1$ direction and Neumann boundary conditions at the other two boundaries. We also consider a sink at the center of the domain $(L_x/2,L_y/2.)$ Let the random coefficient $a(x,\omega):D \times \Omega \rightarrow \mathbb{R}$ be bounded and strictly positive, 
\begin{equation}\label{RK:eq:rf_bound}
	0 < a_l \leq a(x,\omega) \leq a_u < \infty \quad \rm{a.e.} \quad \rm{in} \quad D\times \Omega.
\end{equation}
In the first case, we solve the boundary-value problem:
\begin{align}\label{RK:eq:spde}
	-\nabla \cdot (a(x,\omega) \nabla u(x,\omega))&=f(x,\omega) \;\; \rm{in}~D\times \Omega,  \nonumber \\  
	u(x,\omega)&=50 \;\; \text{on}~x_1 = 0, \nonumber \\
	u(x,\omega)&=25 \;\; \text{on}~x_1 = 240, \nonumber \\
	\vec{n} \cdot \nabla u(x,\omega) & = 0  \;\; \text{on}~x_2 = 0, \nonumber \\
	\vec{n} \cdot \nabla u(x,\omega) & = 0  \;\; \text{on}~x_2 = 60,
\end{align} 
where  $u(x,\omega):D \times \Omega \rightarrow \mathbb{R}.$ Here, we consider a deterministic $f(x,\omega)$ given by,

\begin{equation}\label{RK:eq:fx}
  f(x,\omega)=   
  	\begin{cases}
	  -1.0, & \text{if } (x_1,x_2) = \left(\frac{L_x}{2},\frac{L_y}{2}\right), \\
	  0.0, & \text{otherwise}.
	\end{cases}
\end{equation} 
We assume that the random coefficient $a(x,\omega) = \exp[g(x,\omega)]$ has log-normal distribution with mean $a_0(x) = 5.0$ and standard deviation $\sigma_a = 2.5.$ The Gaussian random field $g(x,\omega)$ has correlation function 
\begin{equation}\label{RK:eq:cov_a}
	C_g(x,y) = \sigma_g^2 \exp \left(-\frac{(x_1-y_1)^2}{l_1^2} -\frac{(x_2-y_2)^2}{l_2^2} \right) \quad \rm{in} \quad D\times \Omega,
\end{equation}
with standard deviation $\sigma_g = \sqrt{\ln \left(1+\frac{\sigma_a}{a_0(x)^2}\right)},$ mean $g_0(x) = \ln \left (\frac{a_0(x)}{(\sqrt{1+\frac{\sigma_a}{a_0(x)^2}})} \right),$ and the correlation lengths $l_1 = 24$ and $l_2 = 20.$
We decompose the spatial domain into $N_D$ subdomains and independently find a low-dimensional solution in each subdomain  using the KL expansion and  basis adaptation methods described in Sections~\ref{RK:sec:hkle} and \ref{RK:sec:ba}. In the first case, we consider the random coefficient that can be represented with $d=10$ random dimensions. To evaluate the accuracy of the solution obtained with the reduced model, we solve the full stochastic system in the domain $D$ with stochastic dimension 10 and sparse-grid level 5 to compute the PC expansion of order 3. This requires $Q_{\xi}=8761$ collocation points.  In the domain decomposition with the basis adaptation approach, we initially solve the full system in the entire domain $D$ with a coarse grid (sparse-grid level 2) and dimension 10, which requires $Q_{\xi}=21$ collocation points. In the second case, we consider a $40$-dimensional system where the reference solution in the full spatial domain is obtained by $100000$ MC simulation samples. In the domain decomposition basis method, we use the sparse-grid level 2, which requires $Q_{\xi}=81$ collocation points, to compute the Gaussian part of the solution. 

\begin{figure}[t!]
    \centering
    \begin{subfigure}[t]{0.75\textwidth}
        \centering
        \includegraphics[scale=.4]{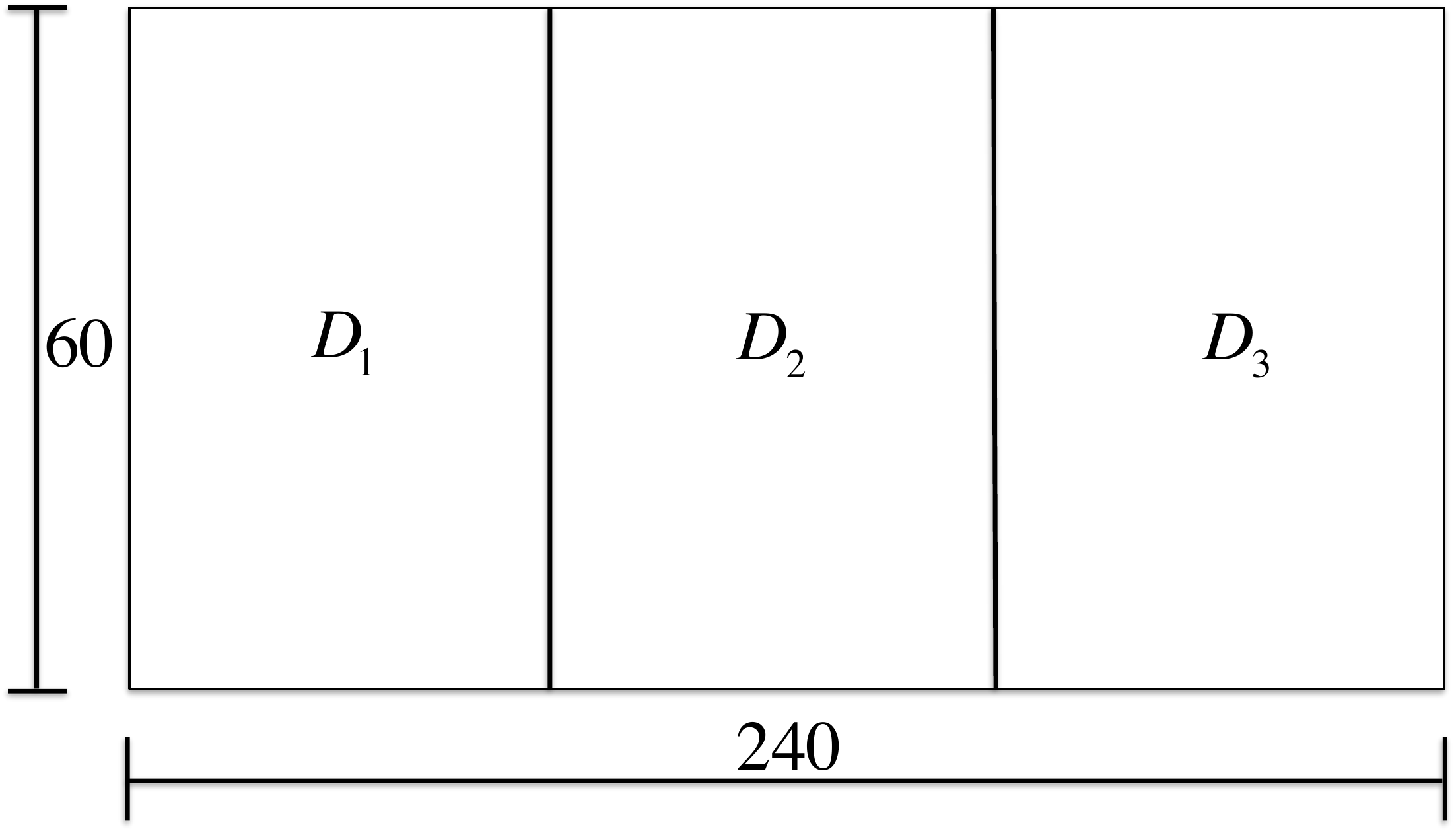}
        \caption{3 subdomains} \label{RK:fig:domain_decomp_3DOM}
    \end{subfigure}        
    \begin{subfigure}[t]{0.75\textwidth}
        \centering
        \includegraphics[scale=.4]{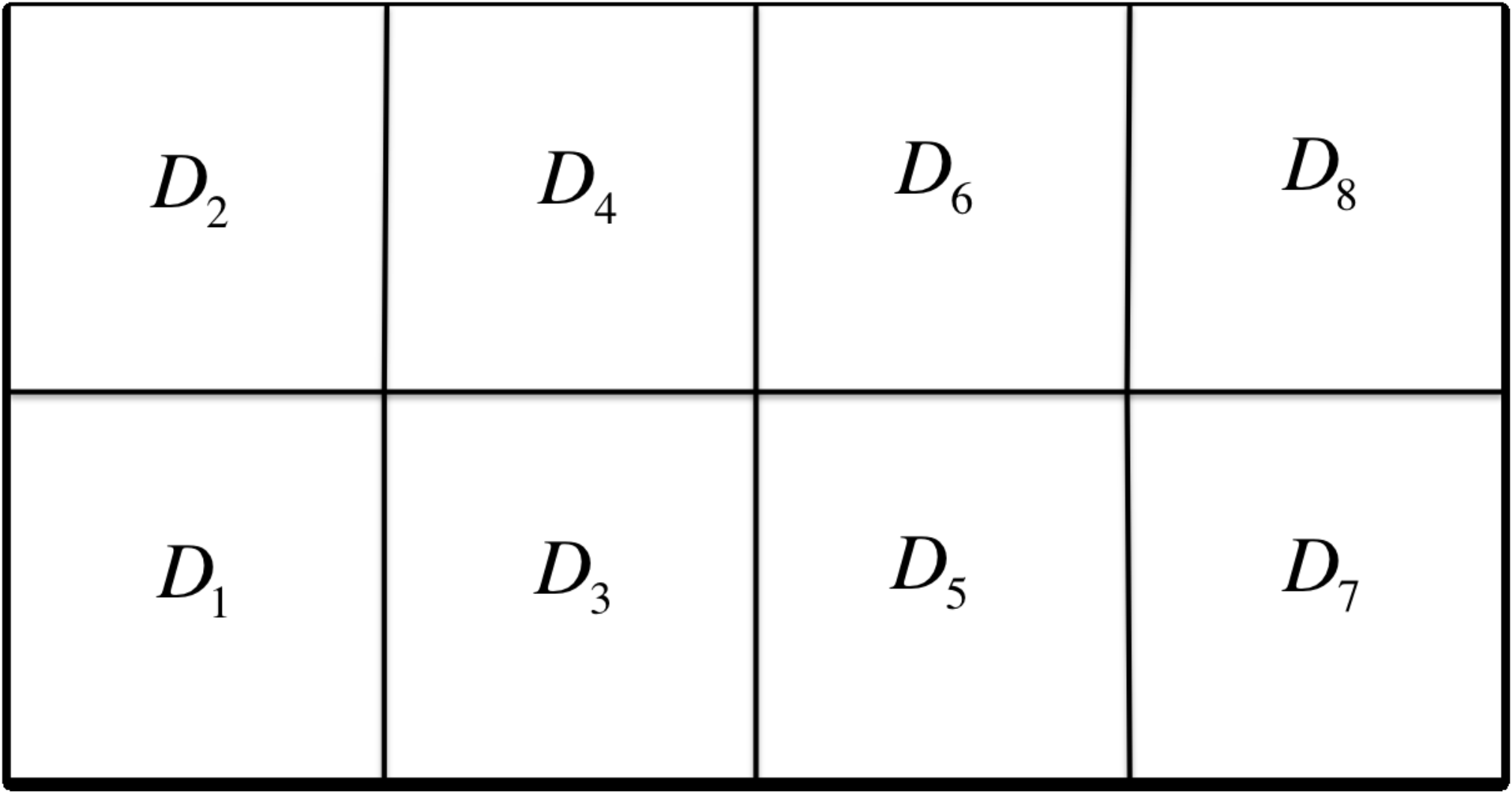}
        \caption{8 subdomains} \label{RK:fig:domain_decomp_8DOM}
    \end{subfigure}    
    \begin{subfigure}[t]{0.75\textwidth}
        \centering
        \includegraphics[scale=.4]{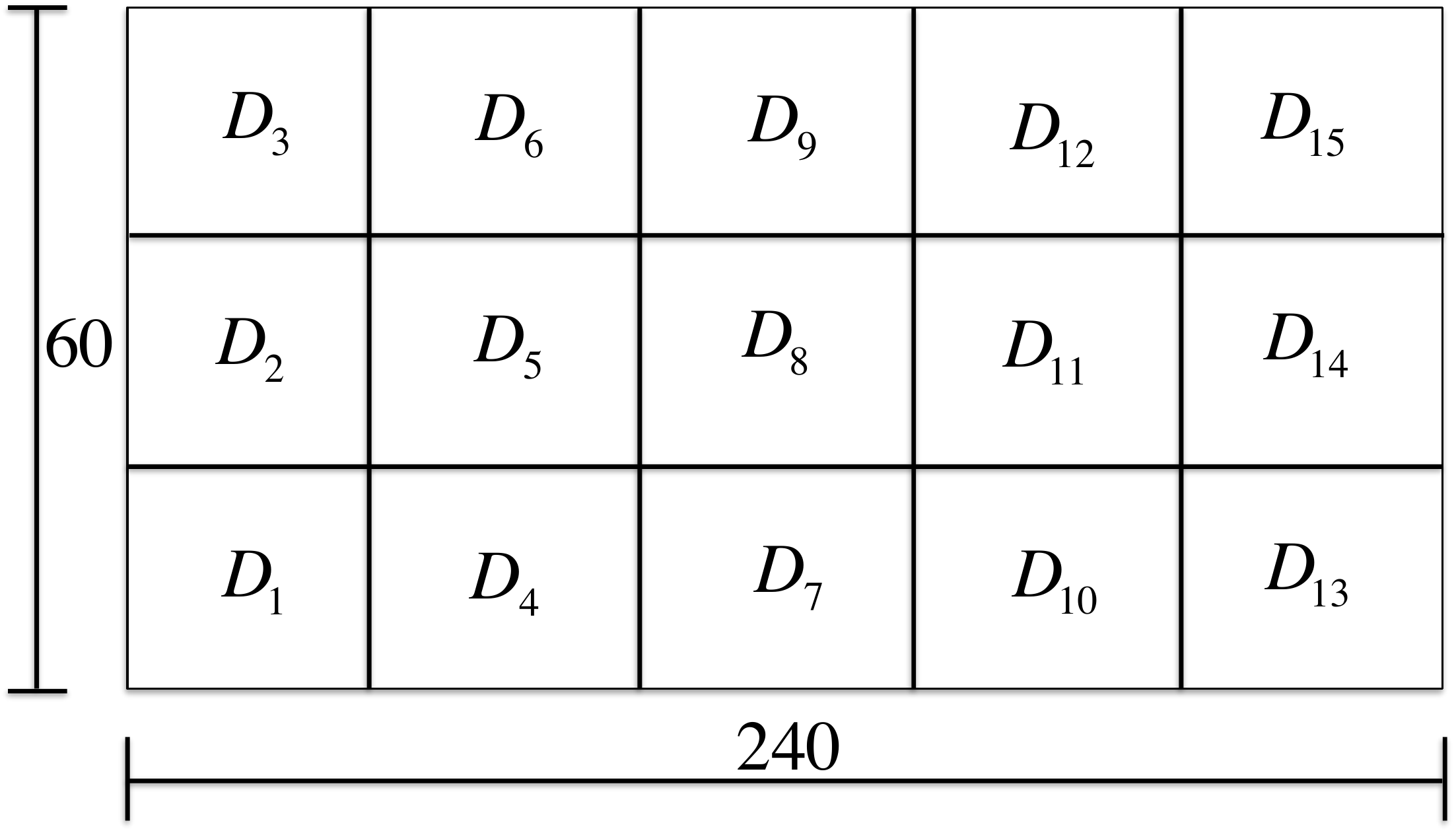}
        \caption{15 subdomains} \label{RK:fig:domain_decomp_15DOM}
    \end{subfigure}  
   
    \caption{Spatial domain decomposed into 3, 8, and 15 non-overlapping subdomains. Domain decomposition into 27 subdomains is not shown here.} \label{RK:fig:domain_decomp_plots}
\end{figure}

\begin{figure}[t!]
    \centering
    \begin{subfigure}[t]{0.45\textwidth}
        \centering
        \includegraphics[height=1.8in]{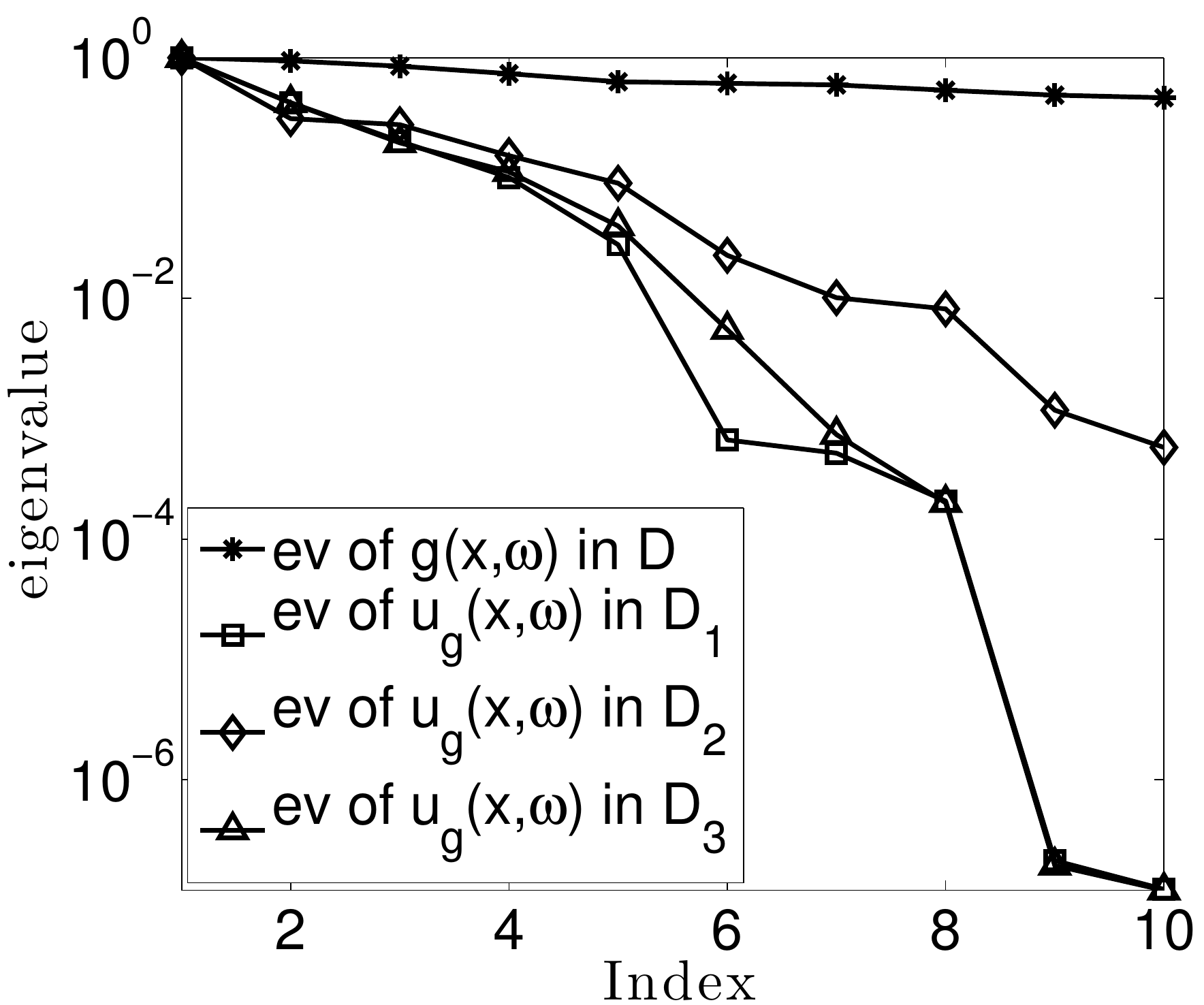}
        \caption{$|D|=3$} \label{RK:fig:eigen_xi_d10_gq2_p1_src0_sink1_bcs50_25_3DOM}
    \end{subfigure}        
    \begin{subfigure}[t]{0.45\textwidth}
        \centering
        \includegraphics[height=1.8in]{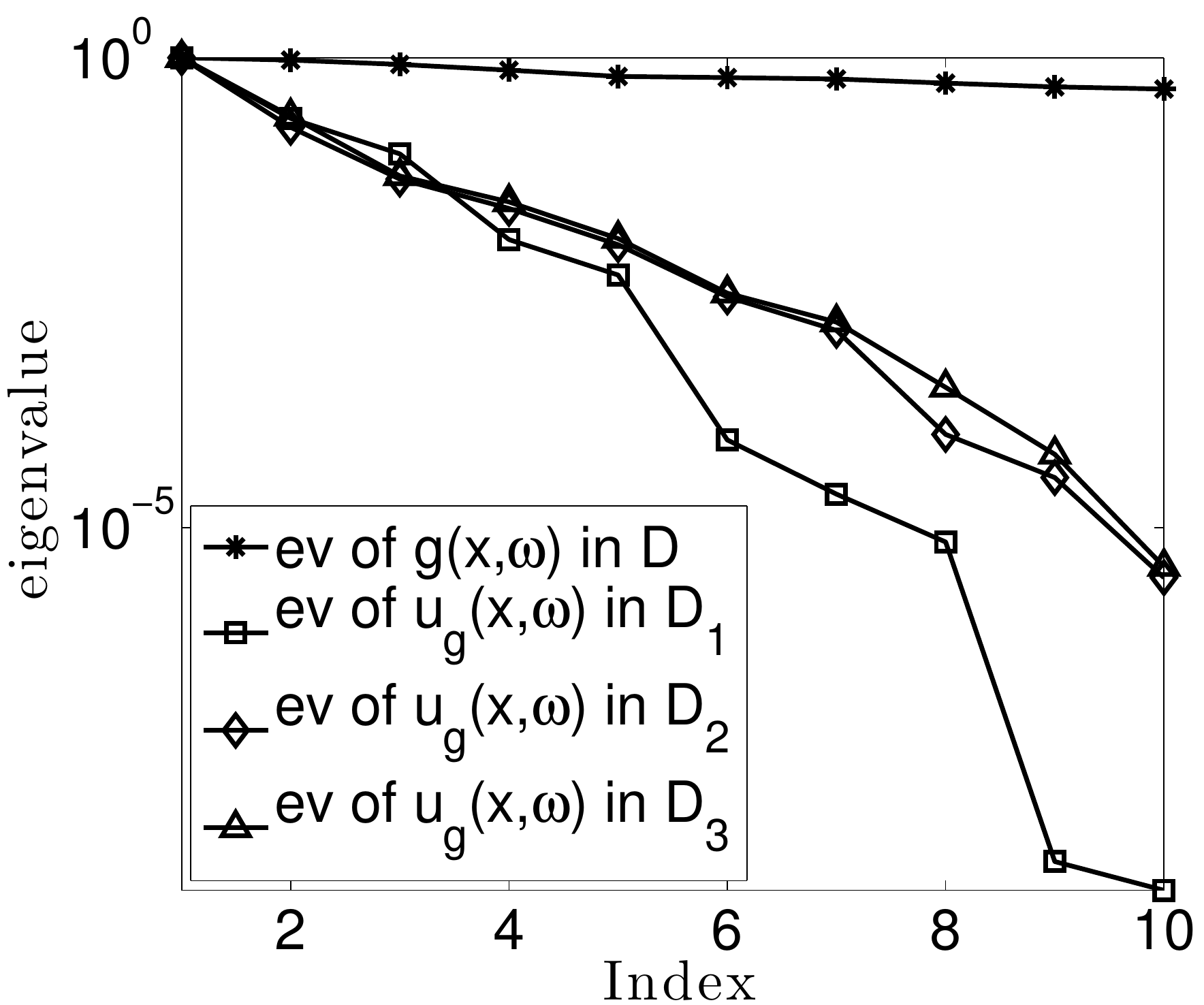}
        \caption{$|D|=8$} \label{RK:fig:eigen_xi_d10_gq2_p1_src0_sink1_bcs50_25_8DOM}
    \end{subfigure}    
    \begin{subfigure}[t]{0.45\textwidth}
        \centering
        \includegraphics[height=1.8in]{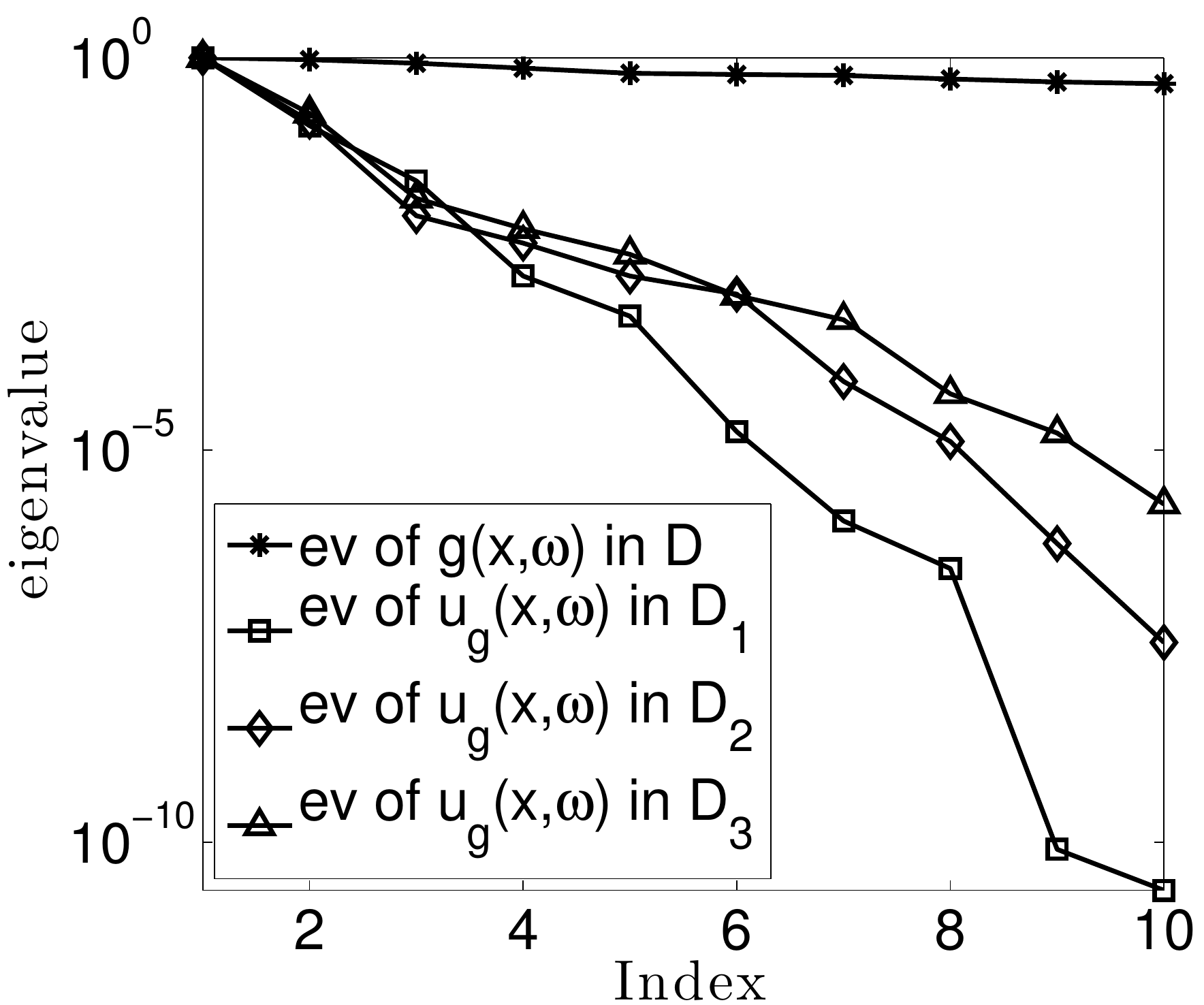}
        \caption{$|D|=15$} \label{RK:fig:eigen_xi_d10_gq2_p1_src0_sink1_bcs50_25_15DOM}
    \end{subfigure}  
    \caption{Decay of the eigenvalues of the covariance function of $g(x,\omega)$ in domain $D$ and the covariance function of the Gaussian solution from subdomains $D_1, D_2,$ and $D_3$ for input random variables $\xi$ with dimension, $d$ = 10, order, $p$ = 1. Spatial domain decomposed into 3, 8,  and 15 non-overlapping subdomains.} \label{RK:fig:eigen_xid10}
\end{figure}

\begin{figure}[t!]
    \centering
    \begin{subfigure}[t]{0.45\textwidth}
        \centering
        \includegraphics[height=1.8in]{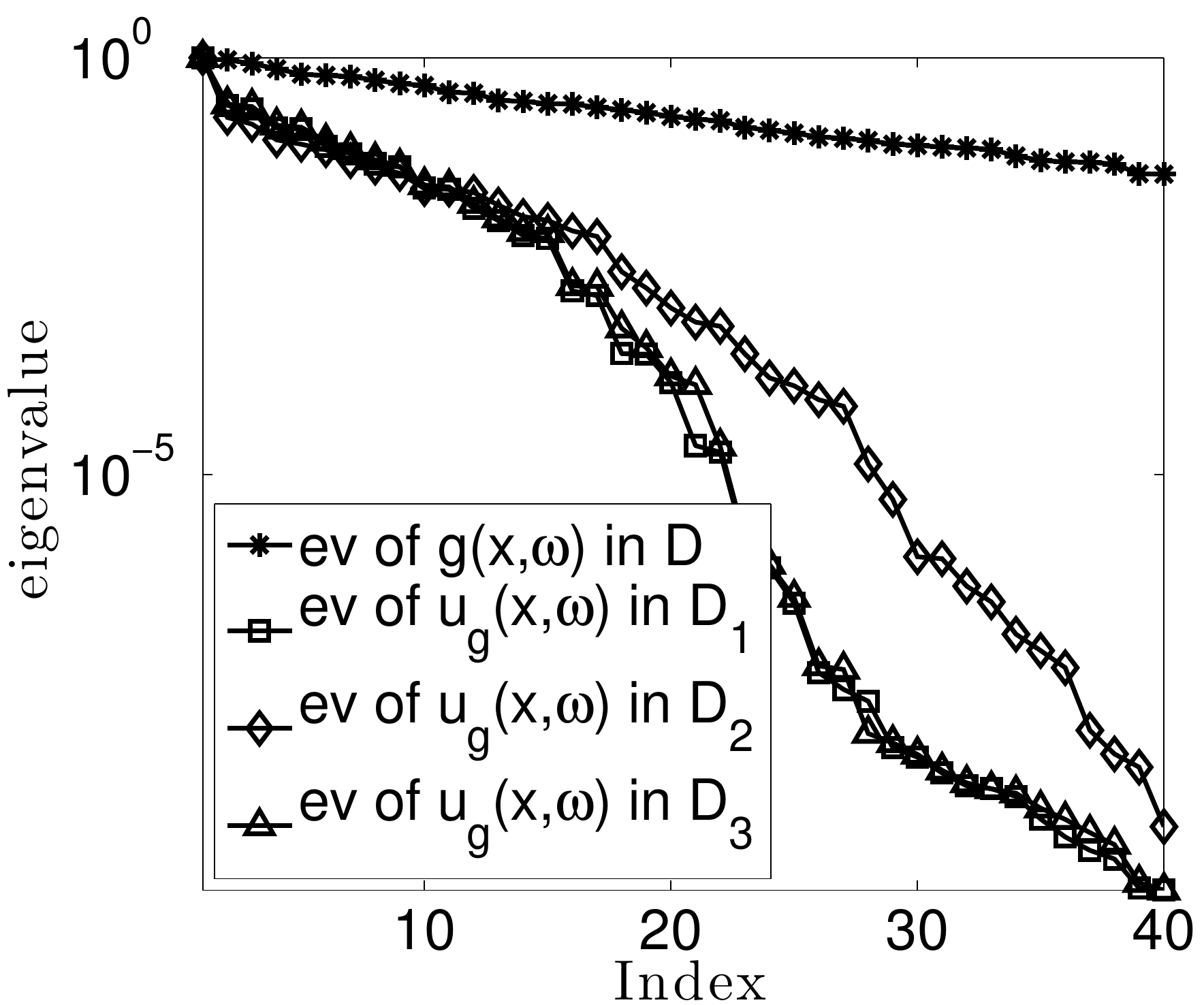}
        \caption{$|D|=3$} \label{RK:fig:eigen_xi_d40_gq2_p1_src0_sink1_mc100000_bcs50_25_3DOM}
    \end{subfigure}        
    \begin{subfigure}[t]{0.45\textwidth}
        \centering
        \includegraphics[height=1.8in]{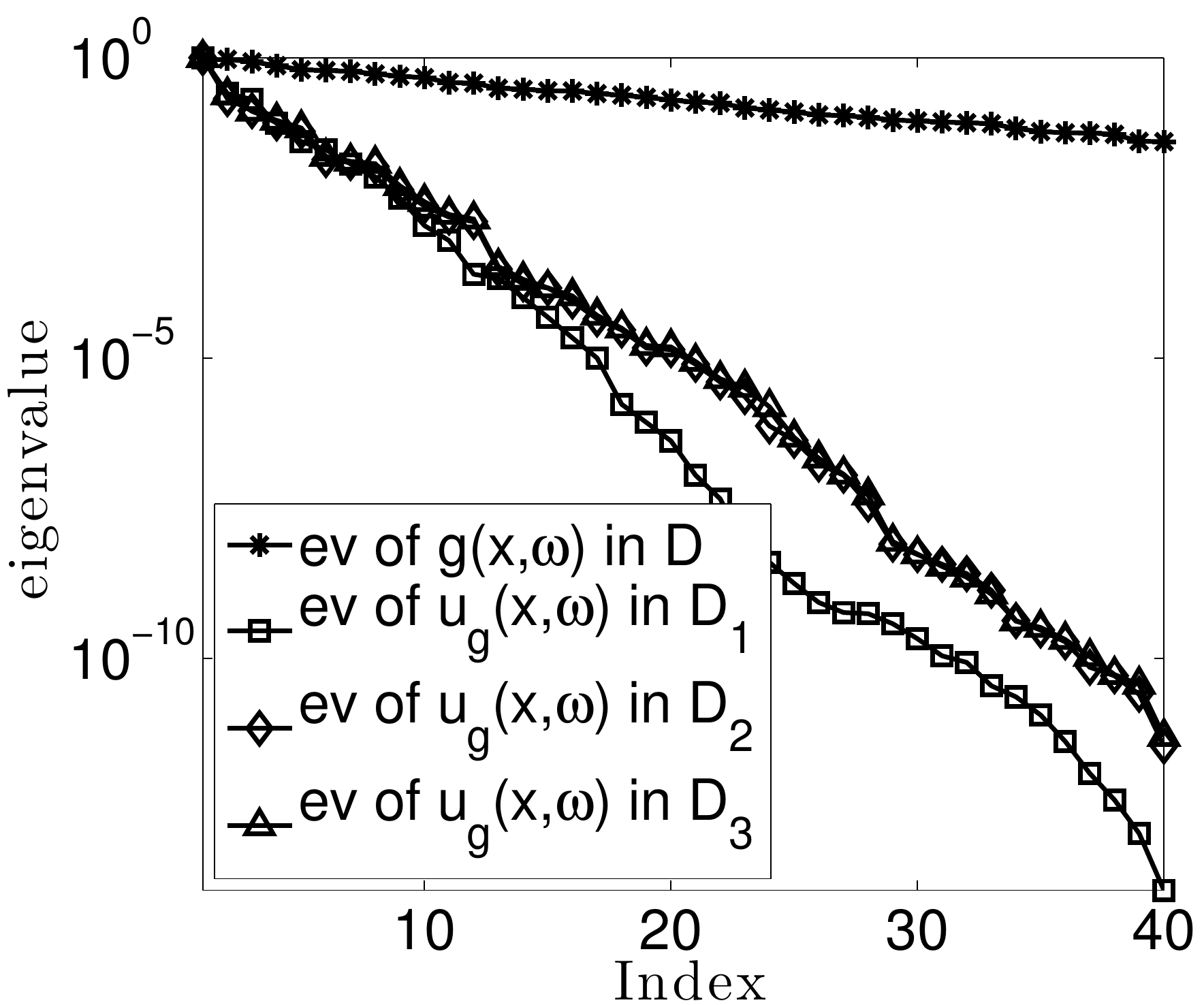}
        \caption{$|D|=8$} \label{RK:fig:eigen_xi_d40_gq2_p1_src0_sink1_mc100000_bcs50_25_8DOM}
    \end{subfigure}    
    \begin{subfigure}[t]{0.45\textwidth}
        \centering
        \includegraphics[height=1.8in]{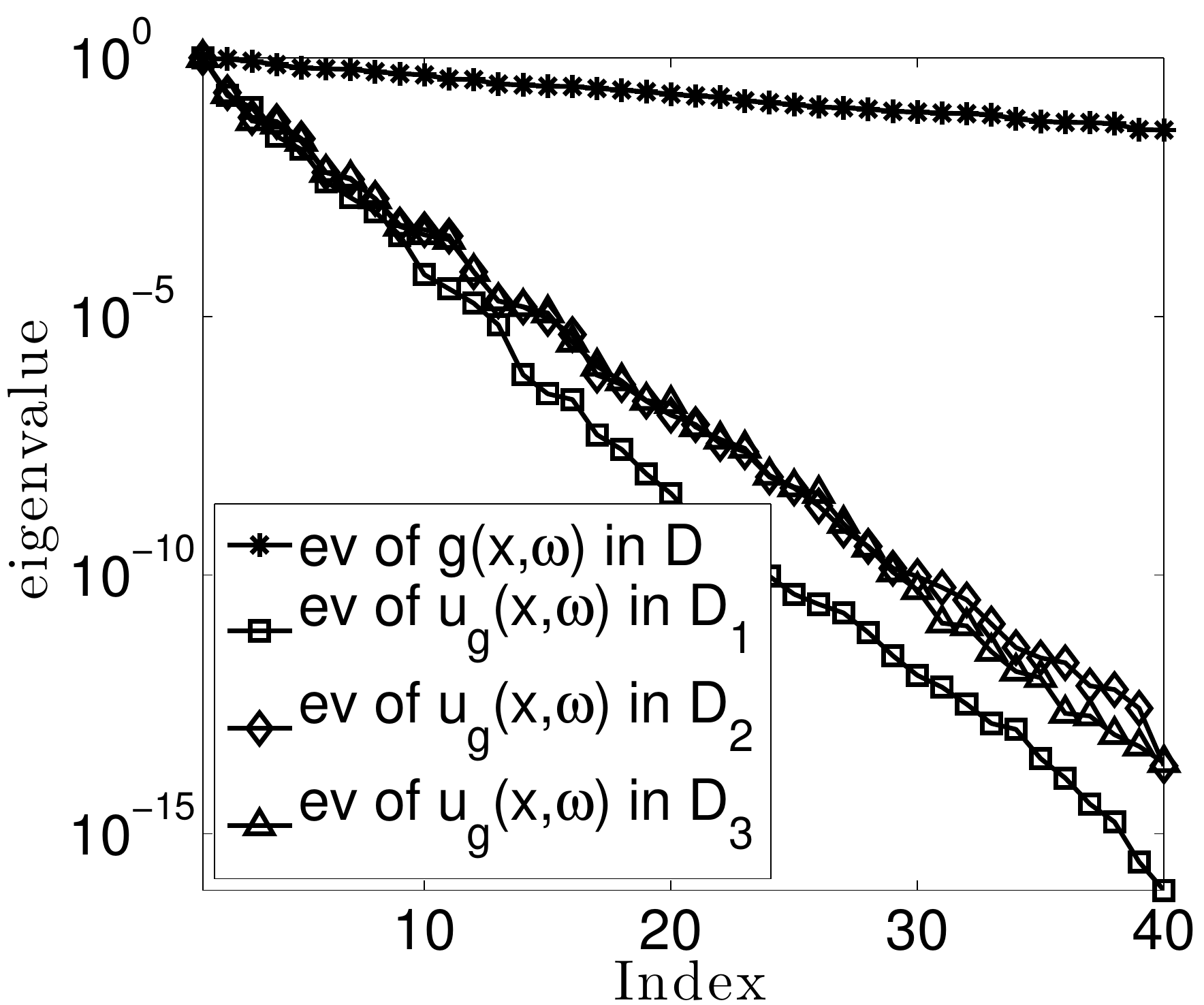}
        \caption{$|D|=15$} \label{RK:fig:eigen_xi_d40_gq2_p1_src0_sink1_mc100000_bcs50_25_15DOM}
    \end{subfigure}  
     \begin{subfigure}[t]{0.45\textwidth}
        \centering
        \includegraphics[height=1.8in]{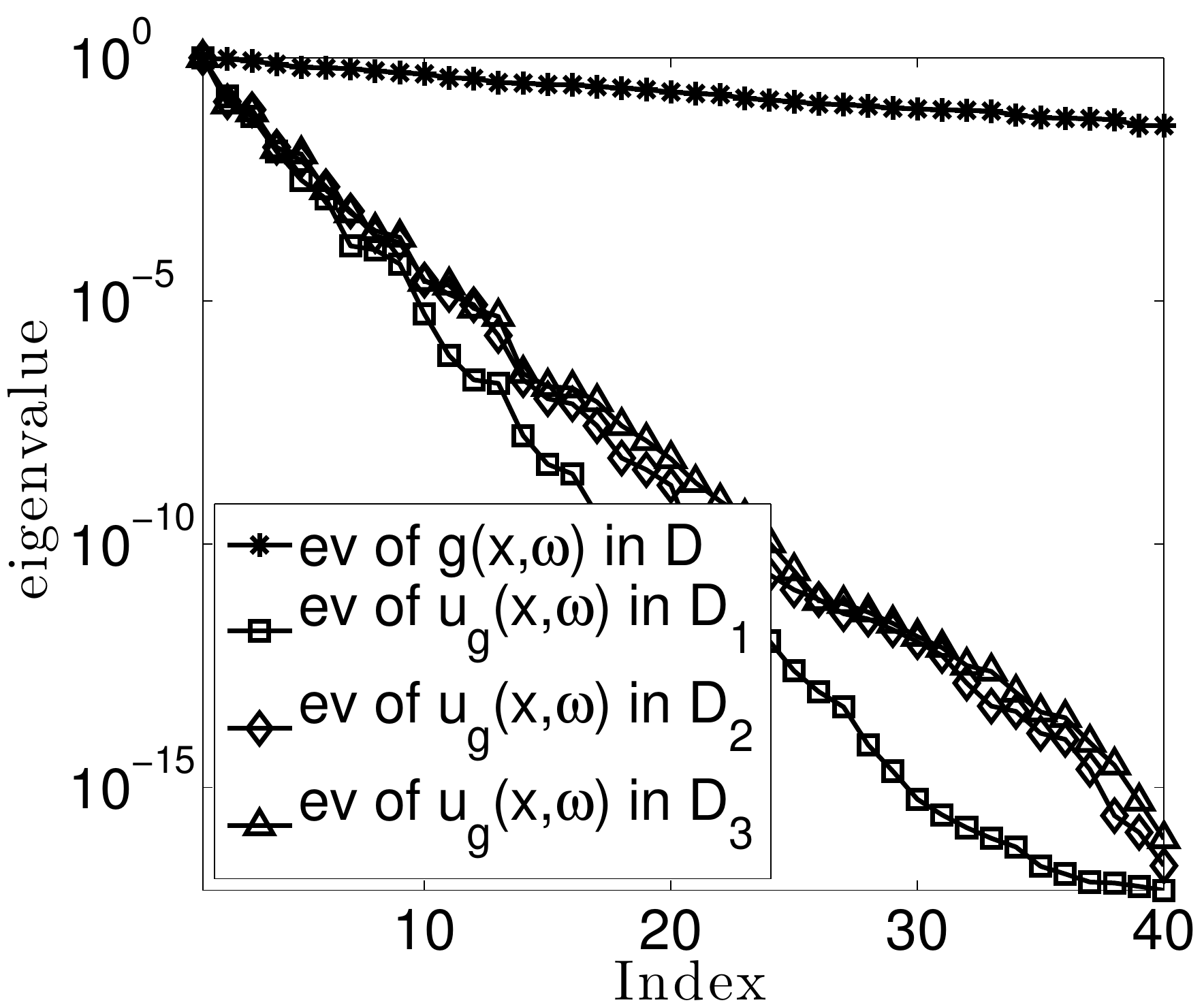}
        \caption{$|D|=27$} \label{RK:fig:eigen_xi_d40_gq2_p1_src0_sink1_mc100000_bcs50_25_27DOM}
    \end{subfigure} 
    
    \caption{Decay of the eigenvalues of the covariance function of $g(x,\omega)$ in domain $D$ and the covariance function of the Gaussian solution from subdomains $D_1, D_2,$ and $D_3$ for input random variables $\xi$ with dimension, $d$ = 40, order, $p$ = 1. Spatial domain decomposed into 3, 8, and 15 non-overlapping subdomains.} \label{RK:fig:eigen_xid40}
\end{figure}

For both cases, we consider the domain decomposition method coupled with basis adaptation, which was discussed in Section~\ref{RK:sec:NN_section}. For both cases ($d=10$ and $d=40$), we solve the system using Algorithm~\ref{RK:alg:NN} for various combinations of numbers of subdomains (3, 8, 15, and 27) and stochastic dimension in each subdomain (3, 4, and 5). Fig. \ref{RK:fig:domain_decomp_plots} shows the domain decomposed into 3, 8, and 15 subdomains. In Fig. \ref{RK:fig:eigen_xid10} and Fig. \ref{RK:fig:eigen_xid40}, we show the decay of eigenvalues of the covariance function of the Gaussian solution in subdomains $D_1, D_2,$ and $D_3$ compared to the input random field $g(x,\omega)$ in the whole domain $D$ for random input dimensions $d=10$ and $d=40,$ respectively. We can observe that the eigenvalues decay faster for smaller domain sizes, requiring smaller stochastic dimensions. 

In all cases we use the non-intrusive method. The computational cost in terms of the flop number in each case is obtained as described in Section~\ref{Nonintru_comp_cost}. The cost of projecting the Schur complement $S^s_{\Gamma}$ and vector $g^s_{\Gamma}$ from one subdomain to another is computed as described in Section \ref{projection_cost}. This cost of projection involves addition of vectors and scalar multiplication of vectors. As such, it is not significant compared to the cost of solving the linear system for the interior points of each subdomain. By comparing the mean and standard deviation with respect to the reference solution in the entire space, we compute the relative error for the solution obtained with domain-decomposition and basis adaptation. We also compare the probability density function (pdf) at an arbitrary spatial point. For the relative error in mean and standard deviation, we use the following norms:
\begin{equation}\label{RK:eq:mu_rel_norm}
	\mu_e = \left \|\frac{\mu_u(x) - \mu_{\tilde{u}}(x)}{\max(\mu_u(x))}\right \|_{2},
\end{equation}
where, $\mu_u(x)$ and $\mu_{\tilde{u}}(x)$ are the mean of the reference and approximate solutions, respectively, and
\begin{equation}\label{RK:eq:sig_rel_norm}
	\sigma_e = \left \|\frac{\sigma_u(x) - \sigma_{\tilde{u}}(x)}{\max(\sigma_u(x))}\right \|_{2},
\end{equation}
where $\sigma_u(x)$ and $\sigma_{\tilde{u}}(x)$ are the standard deviation of the reference and approximate solutions, respectively. 

The computational cost reduction or cost ratio (CR) due to the domain decomposition and basis adaptation is defined as 
\begin{equation}\label{RK:eq:CR}
	\text{CR} =\frac{\text{number of flops for reference solution}}{\text{number of flops for approximate solution}}. 
\end{equation}
In Tables \ref{RK:tab:10d} and \ref{RK:tab:40d}, we compare the relative errors and CRs for the $10-$ and $40$-dimensional problems, respectively. We can see that savings in the computational cost are significant for a reasonably accurate solution. We can achieve higher accuracy by increasing the stochastic dimension at a higher computational cost. 

\begin{figure}[t!]
    \centering
    \begin{subfigure}[t]{0.32\textwidth}
        \centering
        \includegraphics[height=1.2in]{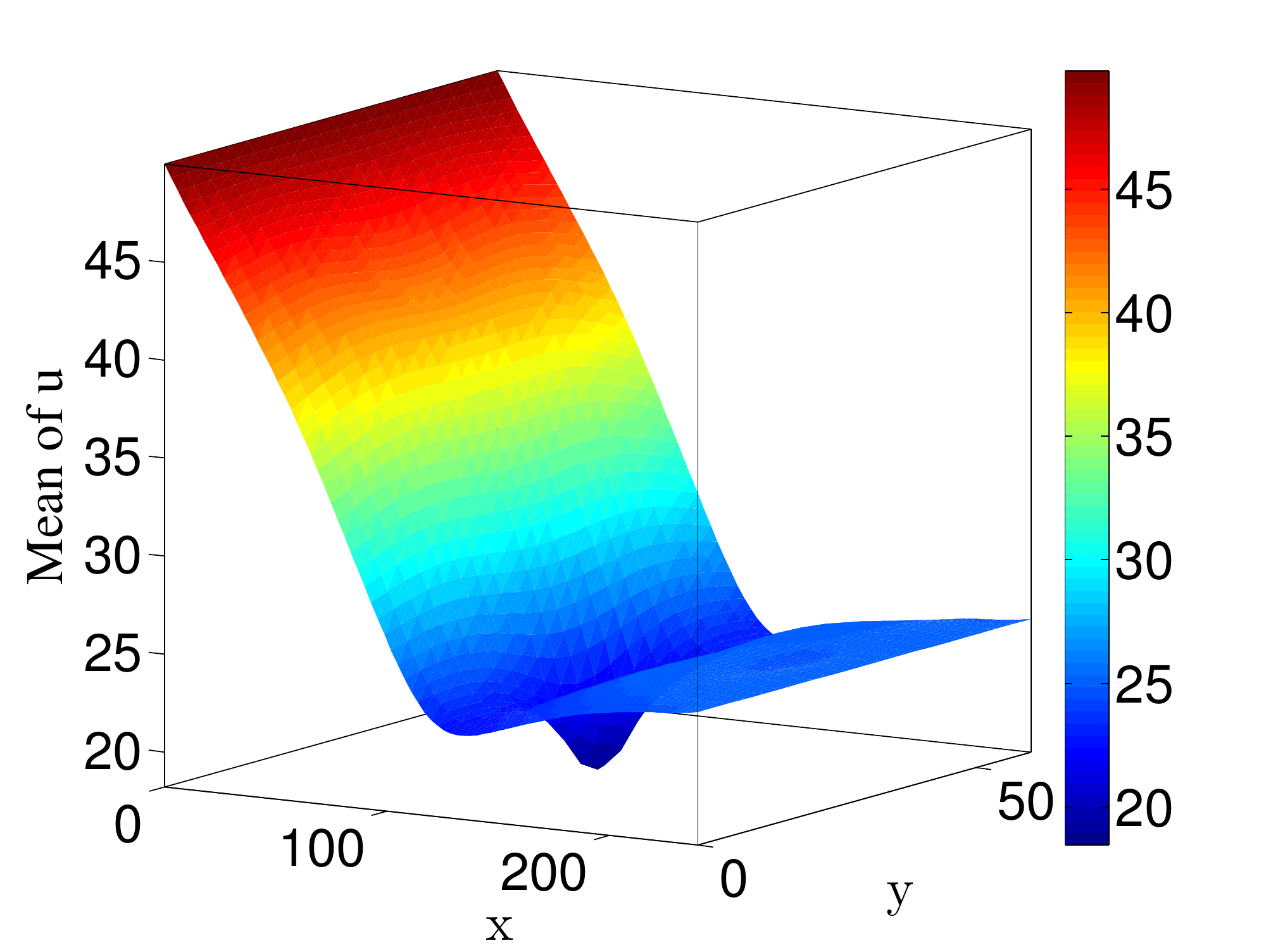}
        \caption{Mean, $\xi, d=10$} \label{RK:fig:u_mean_xid10_p3_gq5_src0_sink1_pce}
    \end{subfigure}        
    \begin{subfigure}[t]{0.32\textwidth}
        \centering
        \includegraphics[height=1.2in]{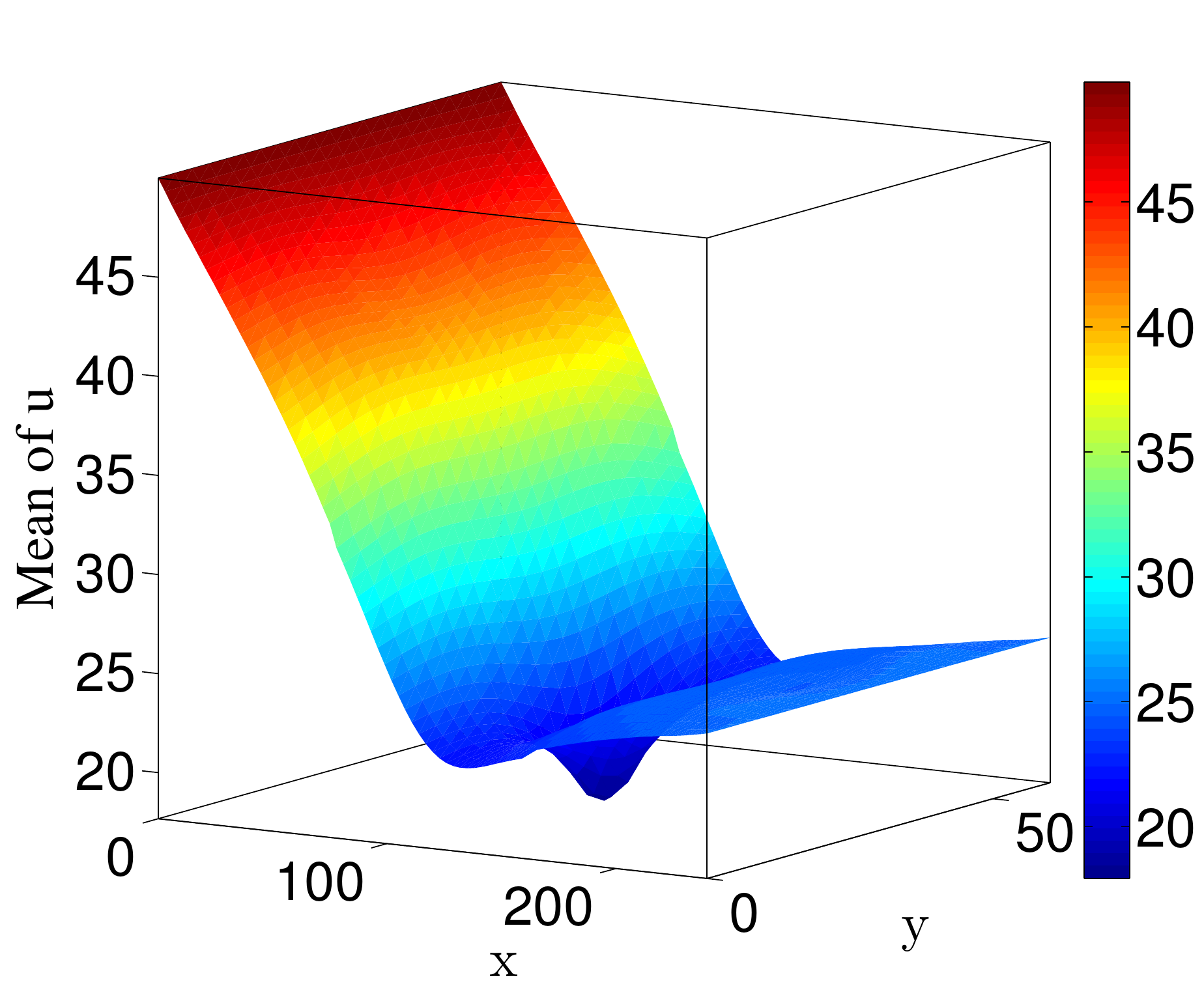}
        \caption{Mean, $\eta, r=5$} \label{RK:fig:u_mean_xid10_p3_gq5_src0_sink1_etad5_p3_3DOM}
    \end{subfigure}    
    \begin{subfigure}[t]{0.32\textwidth}
        \centering
        \includegraphics[height=1.2in]{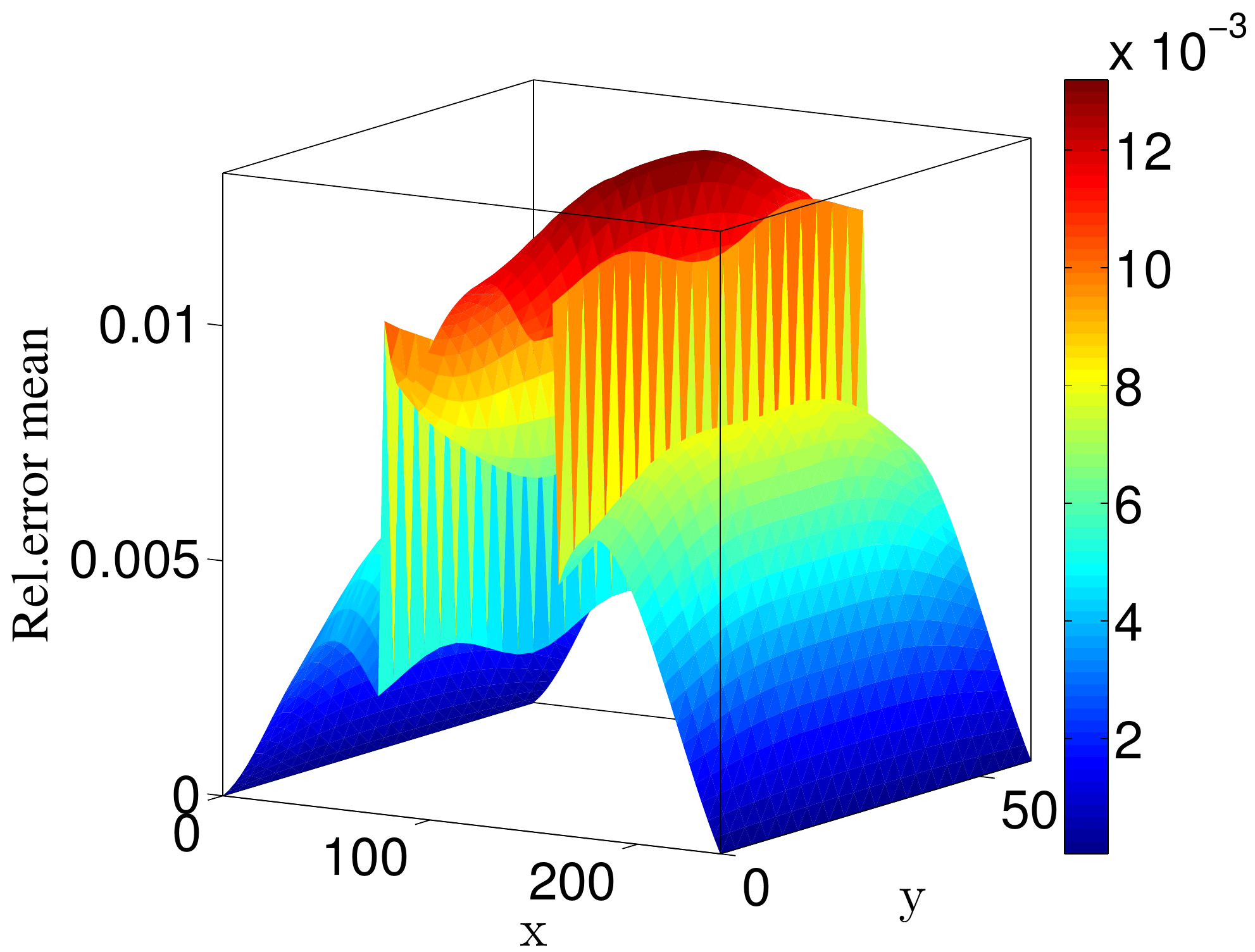}
        \caption{Error, $\eta, r=5$} \label{RK:fig:u_mean_xid10_p3_gq5_src0_sink1_etad5_p3_3DOM_rel_error}
    \end{subfigure}  
   
    \caption{Mean of the solution obtained by stochastic basis adaptation and domain decomposition (3 subdomains) with random variables $\eta$ and dimension, $r$ = 5, order, $p$ = 3, sparse-grid level, $l=5$ compared with the reference solution in $\xi$ of dimension, $d=10$.} \label{RK:fig:u_mean_xid10_etad5}
\end{figure}

\begin{figure}[t!]
    \centering
   \begin{subfigure}[t]{0.3\textwidth}
        \centering
        \includegraphics[height=1.2in]{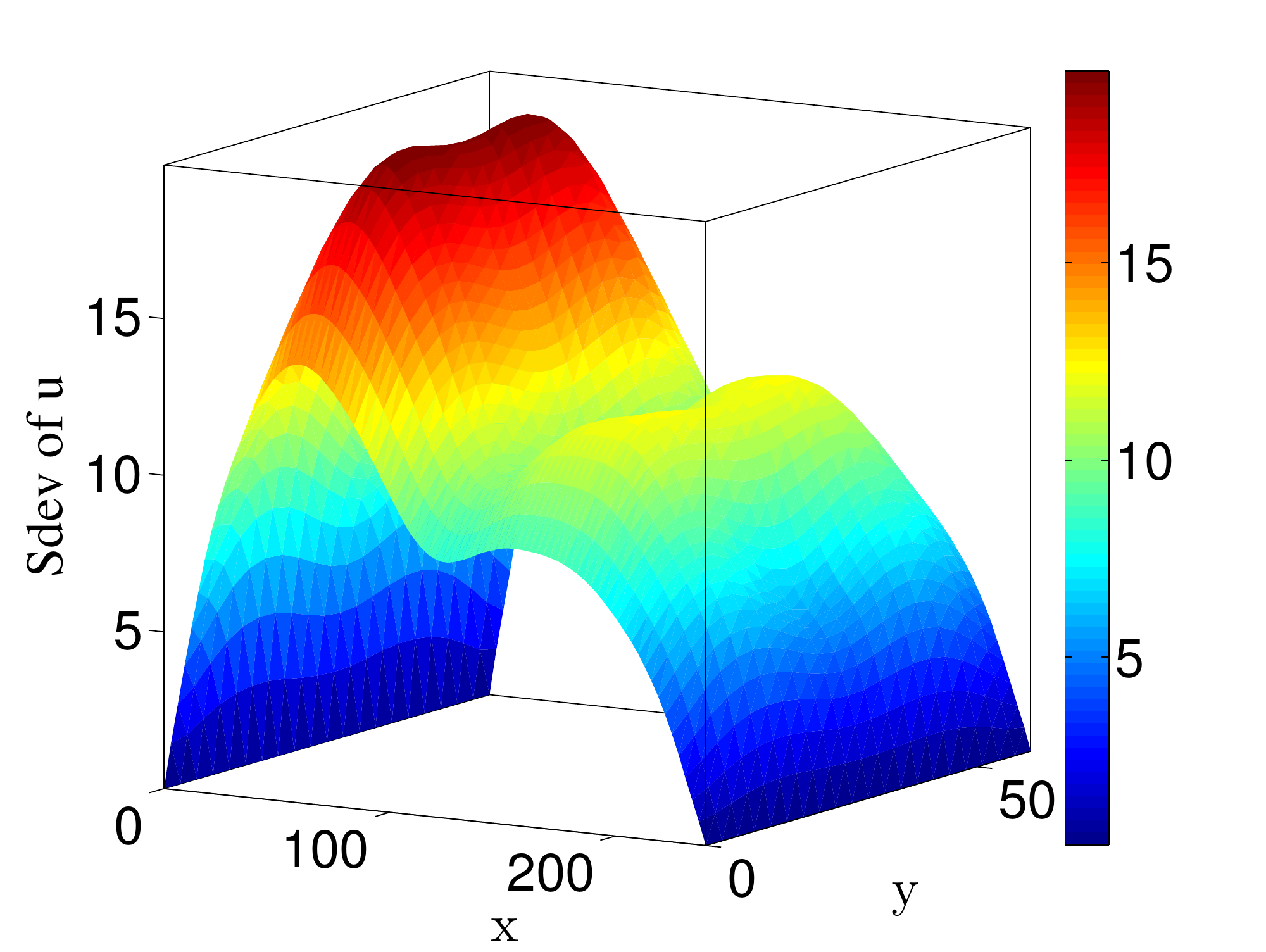}
        \caption{Sdev, $\xi, d=10$} \label{RK:fig:u_sdev_xid10_p3_gq5_src0_sink1_pce}
    \end{subfigure}        
    \begin{subfigure}[t]{0.3\textwidth}
        \centering
        \includegraphics[height=1.2in]{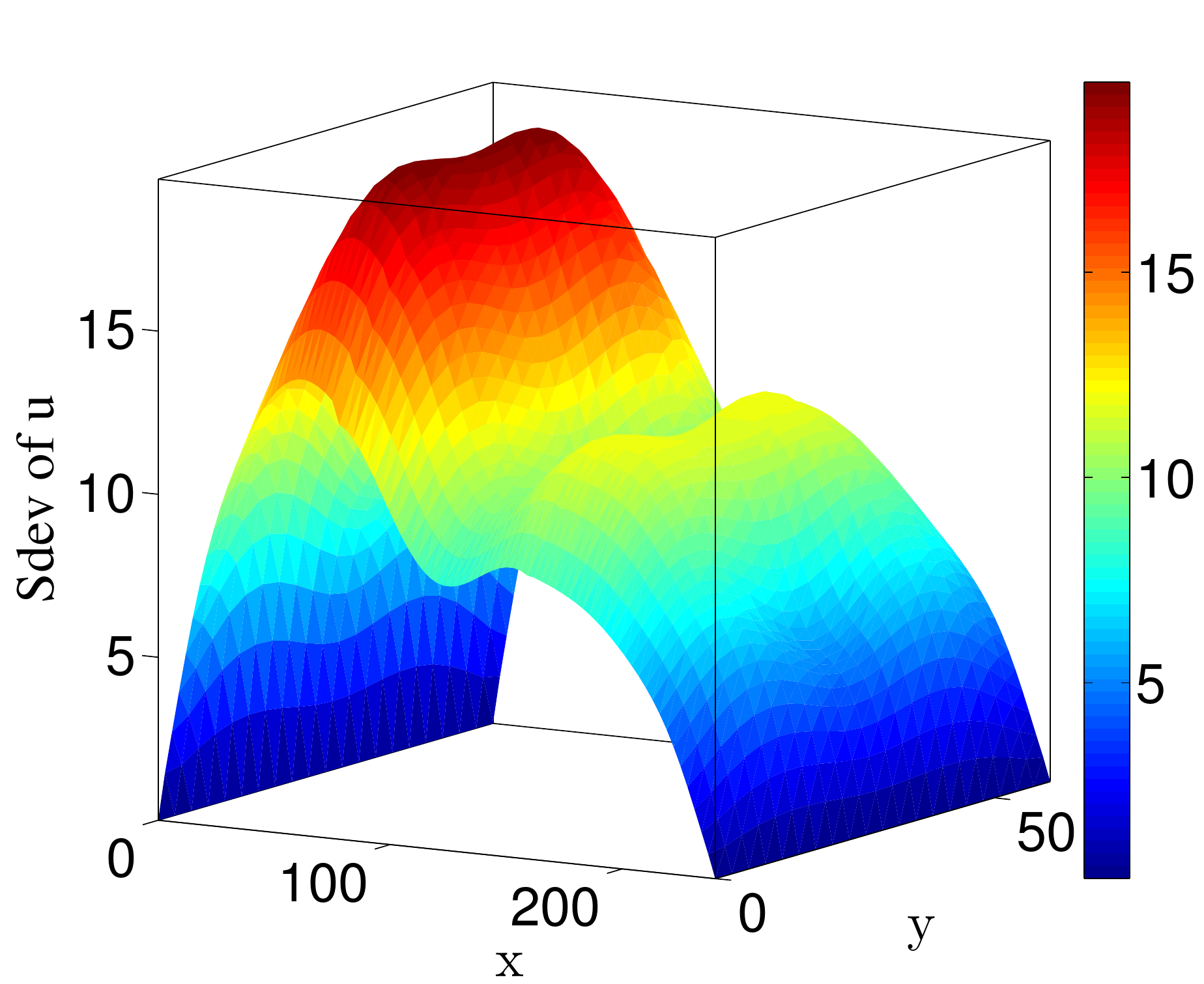}
        \caption{Sdev, $\eta, r=5$} \label{RK:fig:u_sdev_xid10_p3_gq5_src0_sink1_etad5_p3_3DOM}
    \end{subfigure}    
    \begin{subfigure}[t]{0.3\textwidth}
        \centering
        \includegraphics[height=1.2in]{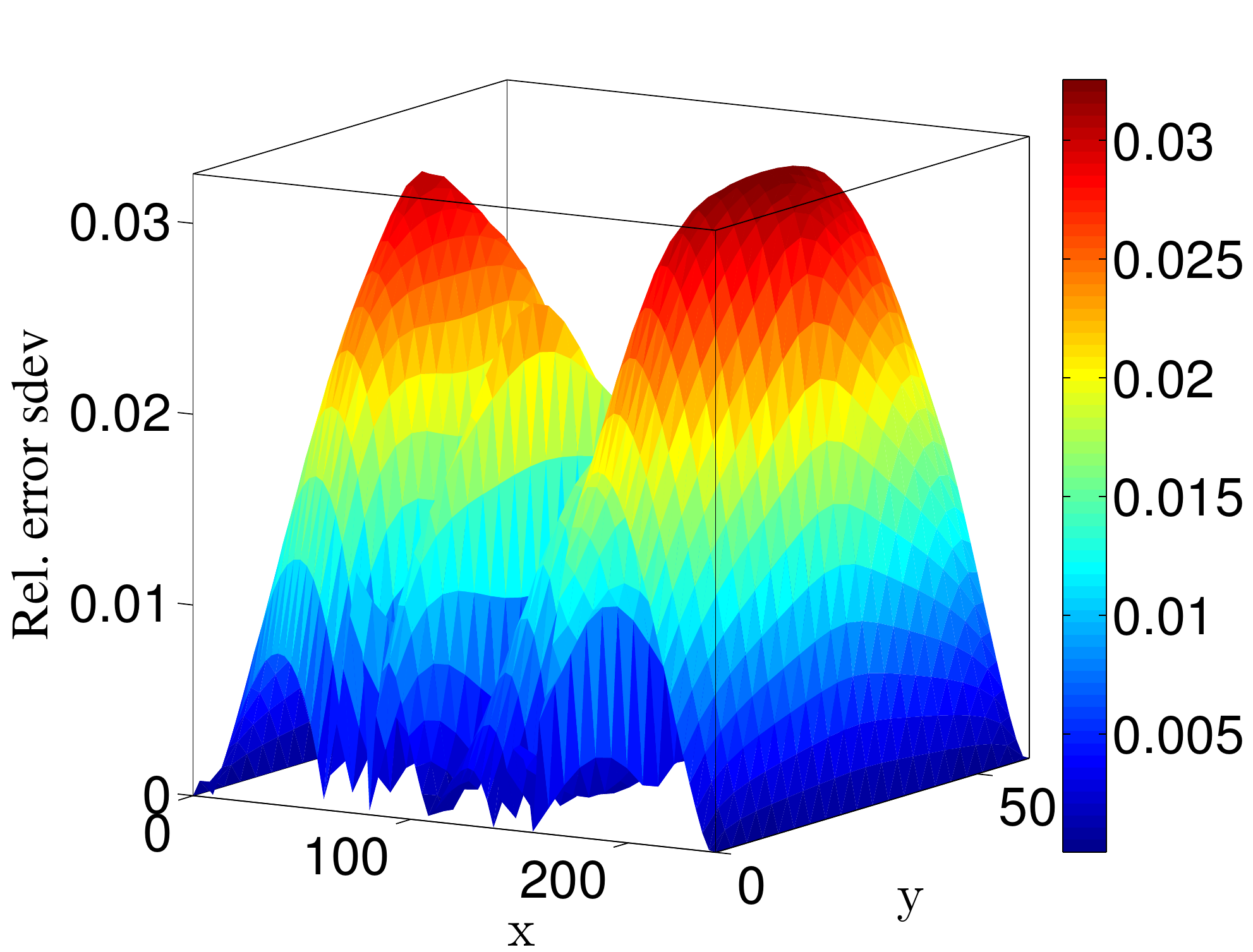}
        \caption{Error, $\eta, r=5$} \label{RK:fig:u_sdev_xid10_p3_gq5_src0_sink1_etad5_p3_3DOM_rel_error}
    \end{subfigure}  
   
    \caption{Standard deviation of the solution obtained by stochastic basis adaptation and domain decomposition (3 subdomains) with random variables $\eta$ and dimension, $r$ = 5, order, $p$ = 3, sparse-grid level, $l=5$ compared with the reference solution in $\xi$ of dimension, $d=10$.} \label{RK:fig:u_sdev_xid10_etad5}
\end{figure}

\begin{figure}[t!]
    \centering
   \begin{subfigure}[t]{0.3\textwidth}
        \centering
        \includegraphics[height=1.2in]{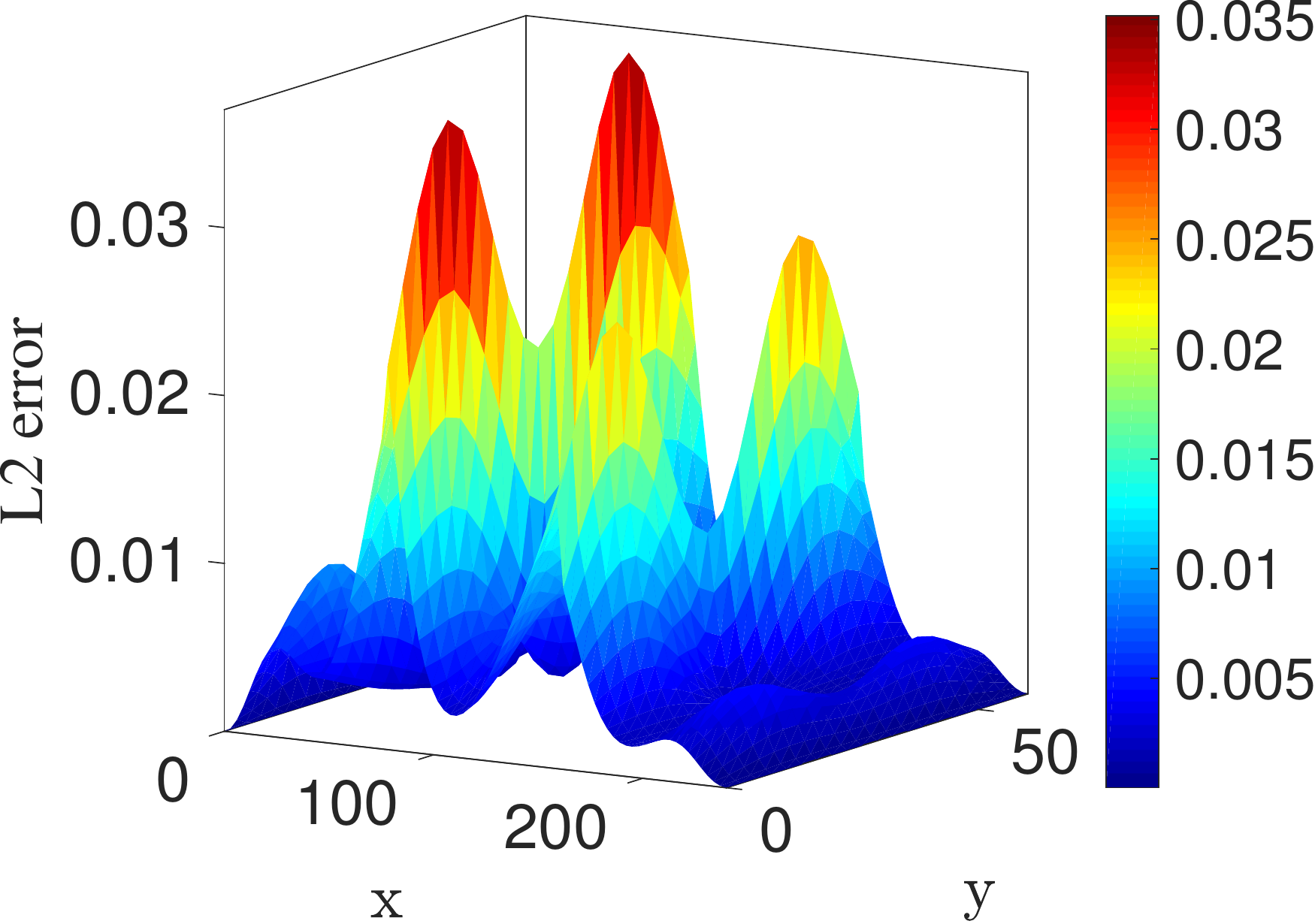}
        \caption{L2 error, $\eta, r=3$} \label{RK:fig:u_L2_xid10_p3_gq5_src0_sink1_etad3_p3_3DOM_rel_error}
    \end{subfigure}        
    \begin{subfigure}[t]{0.3\textwidth}
        \centering
        \includegraphics[height=1.2in]{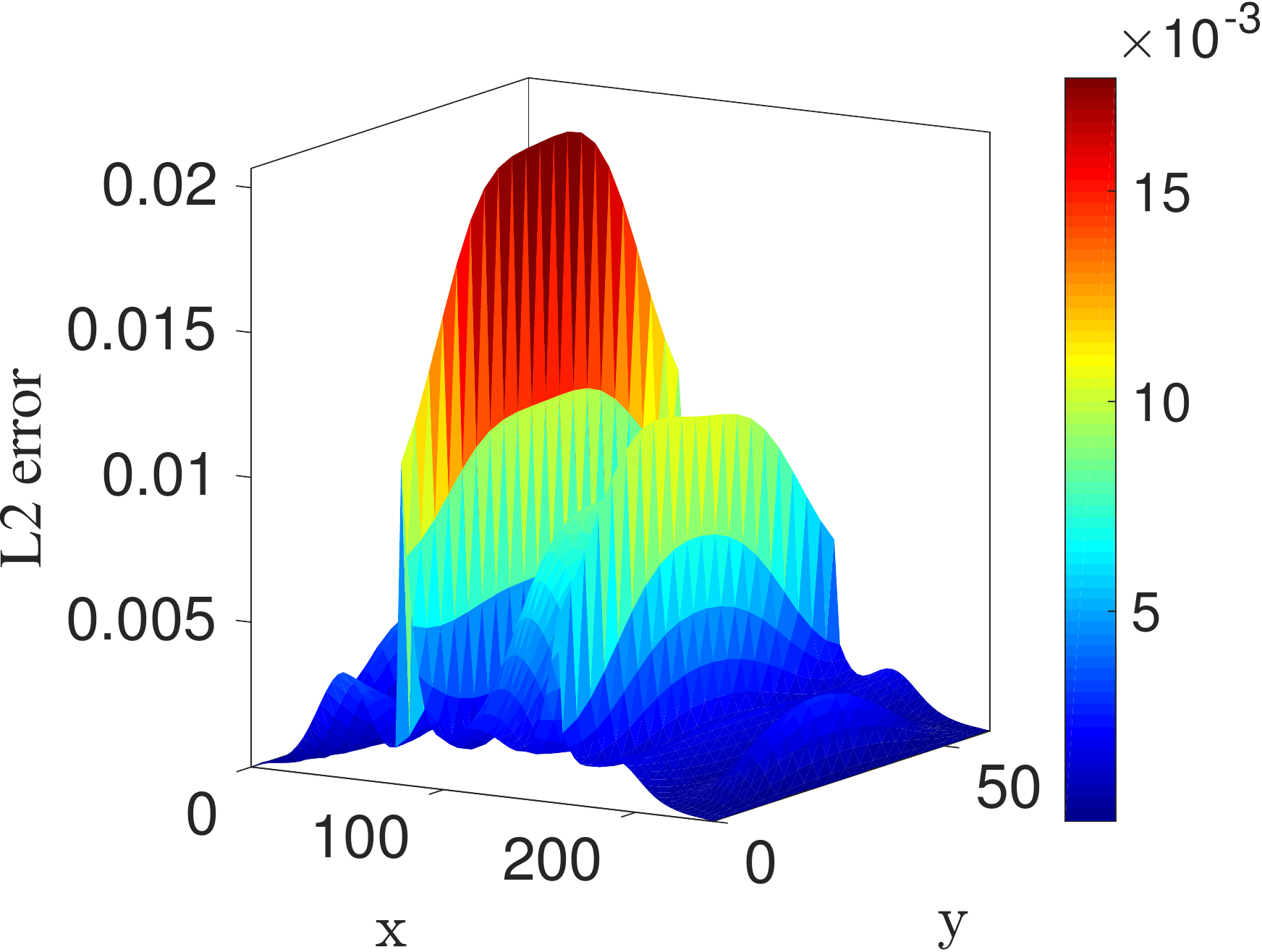}
        \caption{L2 error, $\eta, r=4$} \label{RK:fig:u_L2_xid10_p3_gq5_src0_sink1_etad4_p3_3DOM_rel_error}
    \end{subfigure}    
    \begin{subfigure}[t]{0.3\textwidth}
        \centering
        \includegraphics[height=1.2in]{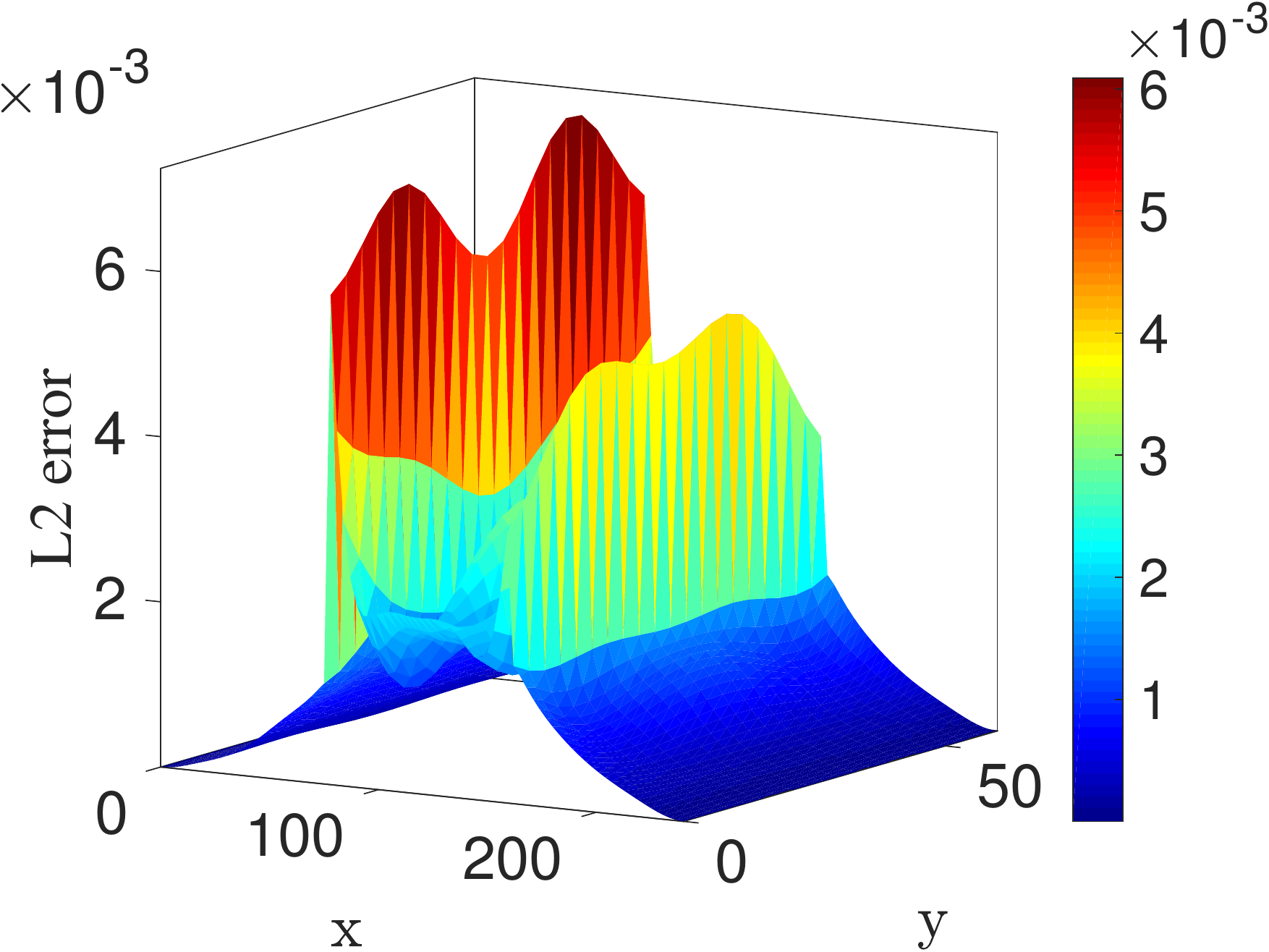}
        \caption{L2 error, $\eta, r=5$} \label{RK:fig:u_L2_xid10_p3_gq5_src0_sink1_etad5_p3_3DOM_rel_error}
    \end{subfigure}  
   
    \caption{L2 norm of the error ($\epsilon(x) = \mathbb{E}[(u(x,\boldsymbol{\xi})-u(x,\boldsymbol{\eta}))^2]$) in the solution obtained by stochastic basis adaptation and domain decomposition (3 subdomains) with random variables $\eta$ and dimension, $r$ = 3, 4 and 5, order, $p$ = 3, sparse-grid level, $l=5$ and the dimension of reference solution in $\xi$ is $d=10$.} \label{RK:fig:uu_L2_xid10_p3_gq5_src0_D3}
\end{figure}

\begin{figure}[t!]
    \centering
   \begin{subfigure}[t]{0.3\textwidth}
        \centering
        \includegraphics[height=1.2in]{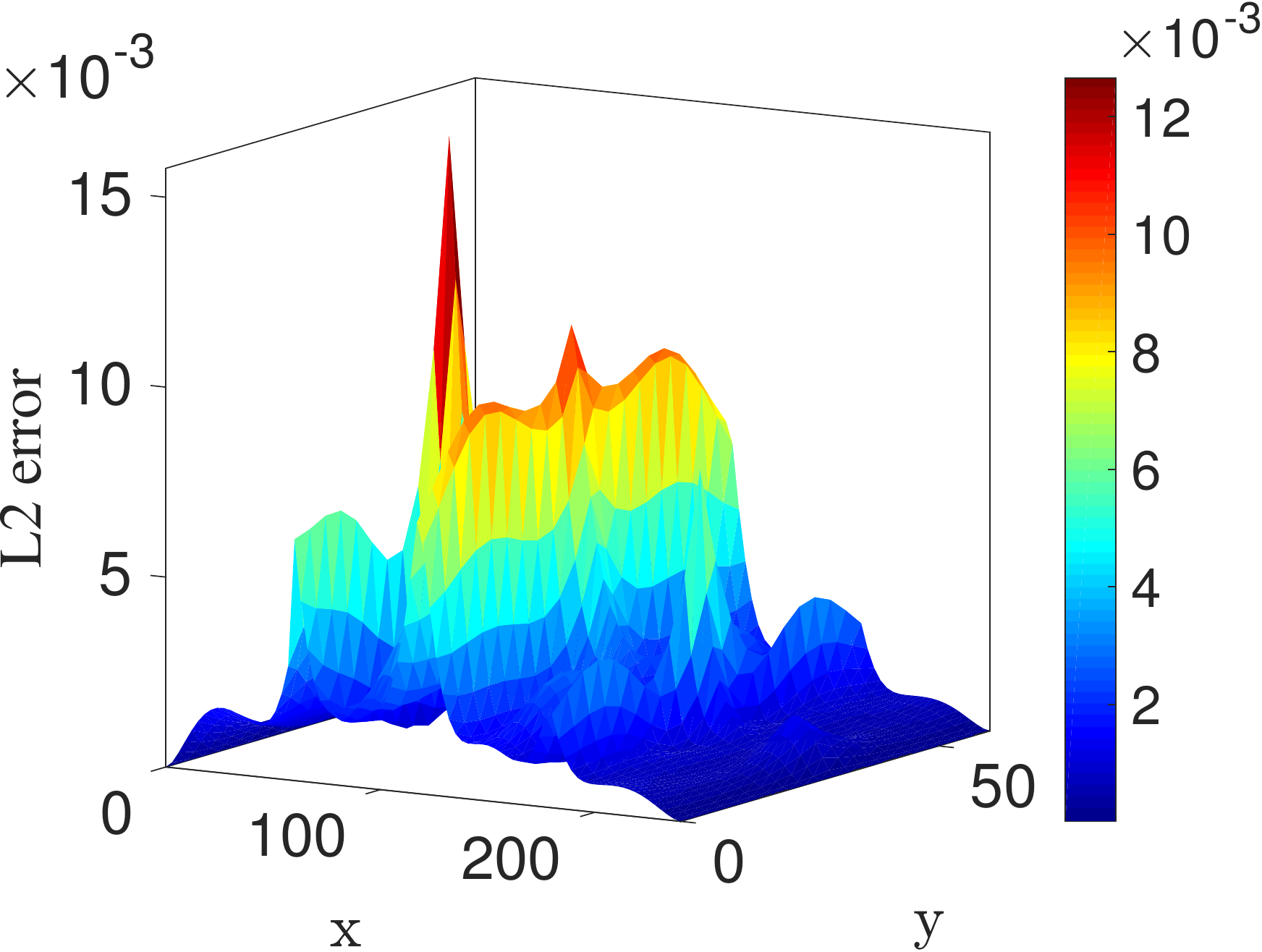}
        \caption{L2 error, $\eta, r=3$} \label{RK:fig:u_L2_xid10_p3_gq5_src0_sink1_etad3_p3_8DOM_rel_error}
    \end{subfigure}        
    \begin{subfigure}[t]{0.3\textwidth}
        \centering
        \includegraphics[height=1.2in]{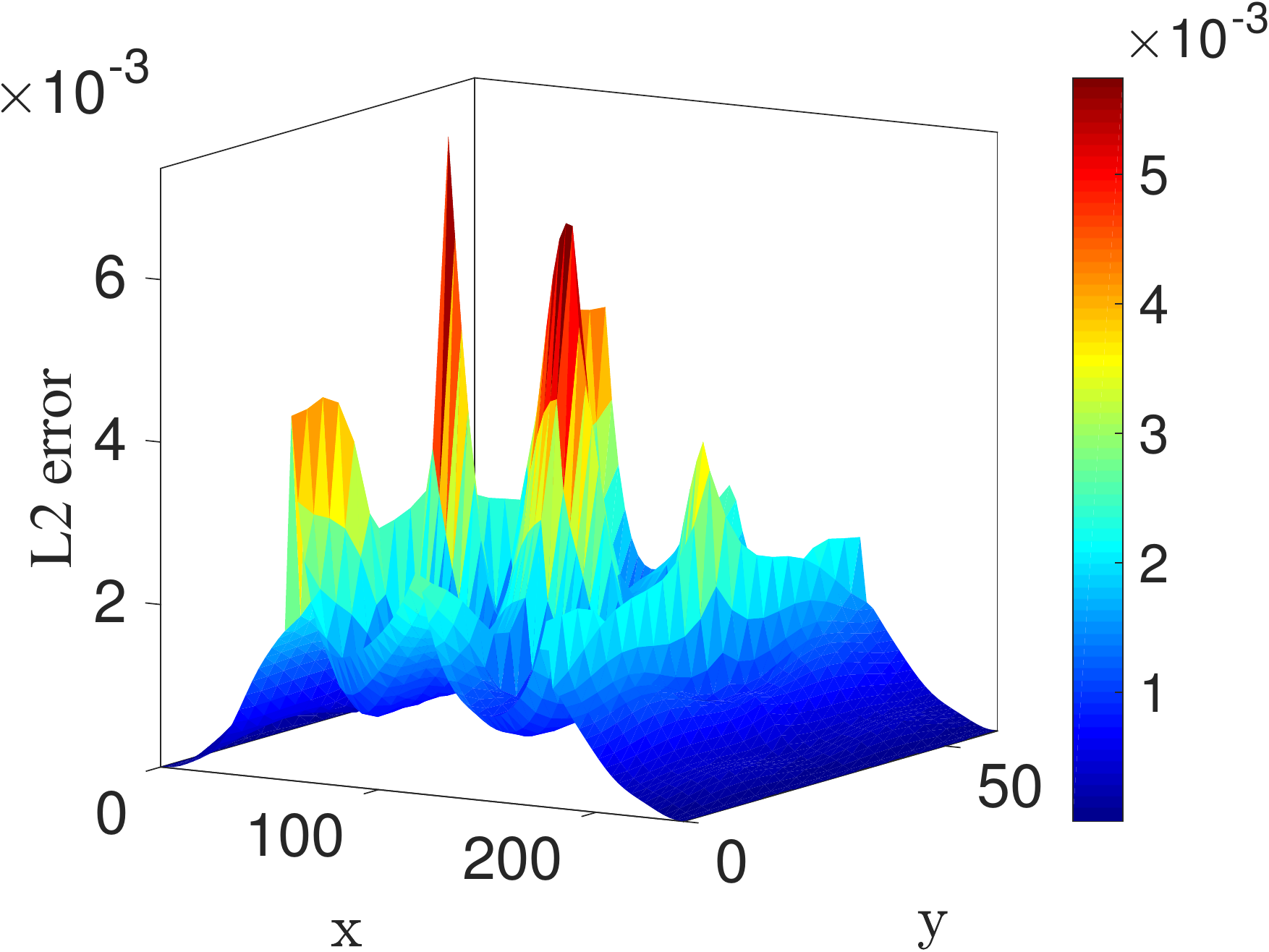}
        \caption{L2 error, $\eta, r=4$} \label{RK:fig:u_L2_xid10_p3_gq5_src0_sink1_etad4_p3_8DOM_rel_error}
    \end{subfigure}    
    \begin{subfigure}[t]{0.3\textwidth}
        \centering
        \includegraphics[height=1.2in]{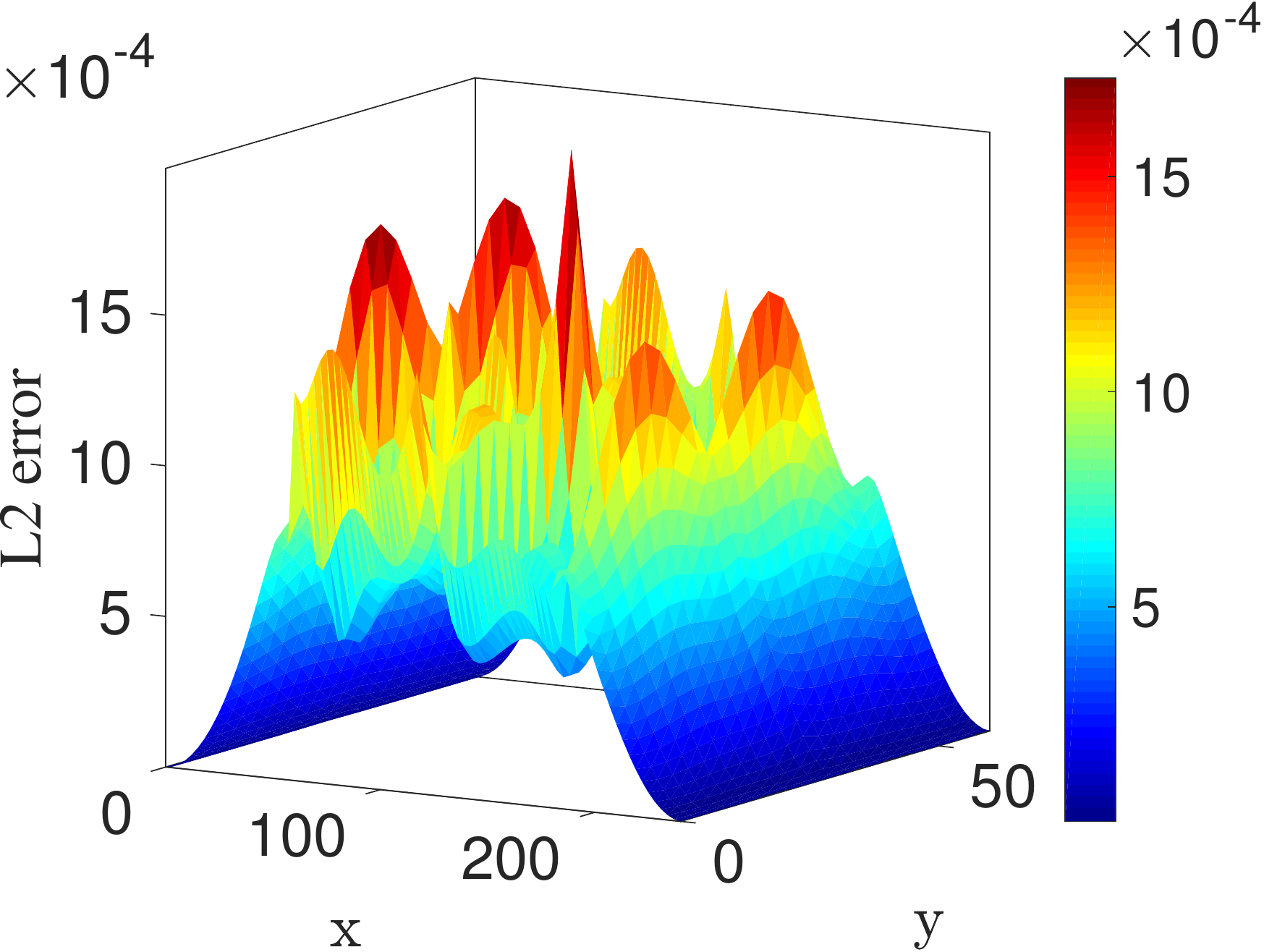}
        \caption{L2 error, $\eta, r=5$} \label{RK:fig:u_L2_xid10_p3_gq5_src0_sink1_etad5_p3_8DOM_rel_error}
    \end{subfigure}  
   
    \caption{L2 norm of the error ($\epsilon(x) = \mathbb{E}[(u(x,\boldsymbol{\xi})-u(x,\boldsymbol{\eta}))^2]$) in the solution obtained by stochastic basis adaptation and domain decomposition (8 subdomains) with random variables $\eta$ and dimension,  $r$ = 3, 4 and 5, order, $p$ = 3, sparse-grid level, $l=5$ and the dimension of reference solution in $\xi$ is $d=10$.} \label{RK:fig:uu_L2_xid10_p3_gq5_src0_D8}
\end{figure}

\begin{figure}[t!]
    \centering
   \begin{subfigure}[t]{0.3\textwidth}
        \centering
        \includegraphics[height=1.2in]{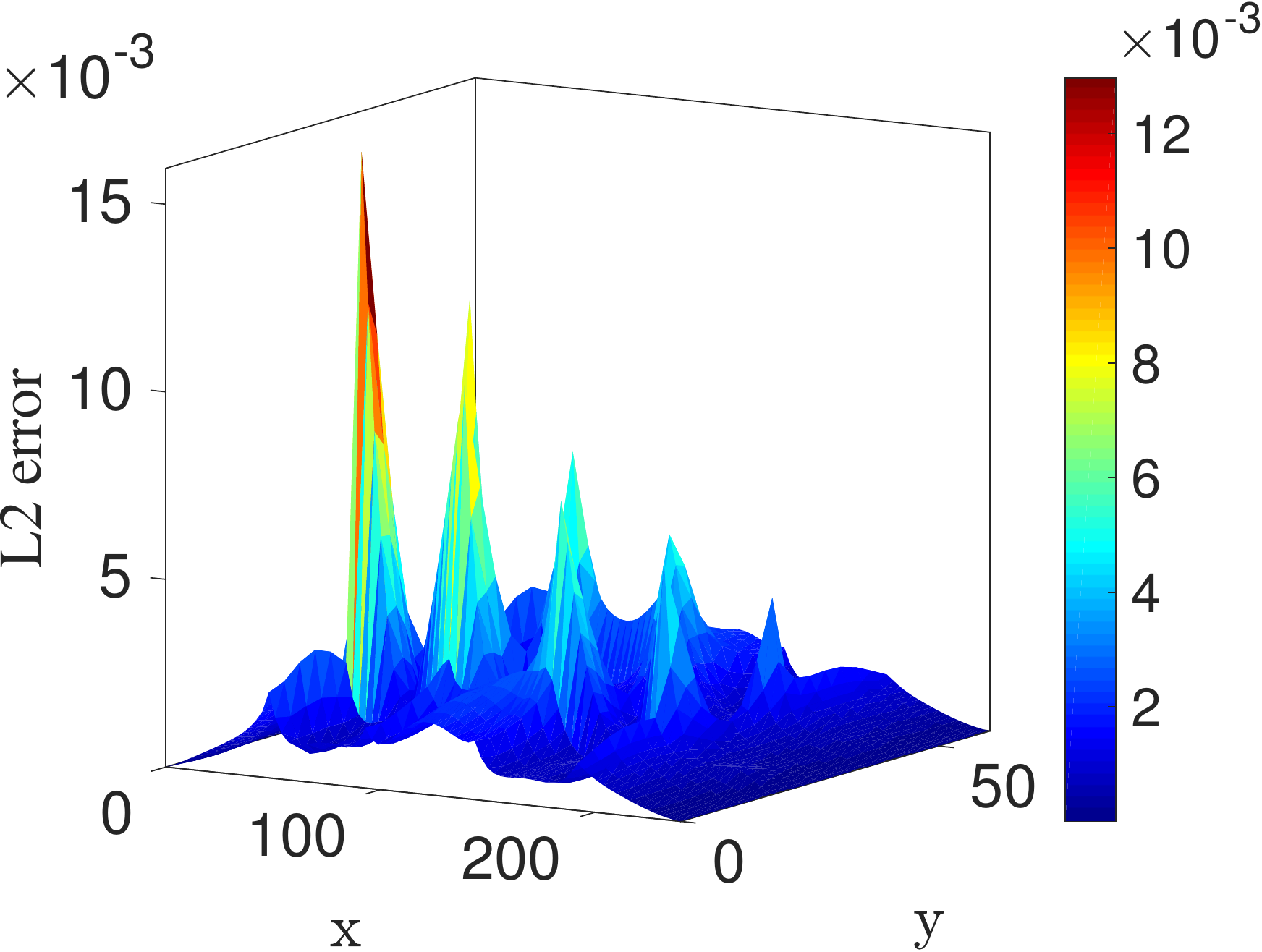}
        \caption{L2 error, $\eta, r=3$} \label{RK:fig:u_L2_xid10_p3_gq5_src0_sink1_etad3_p3_15DOM_rel_error}
    \end{subfigure}        
    \begin{subfigure}[t]{0.3\textwidth}
        \centering
        \includegraphics[height=1.2in]{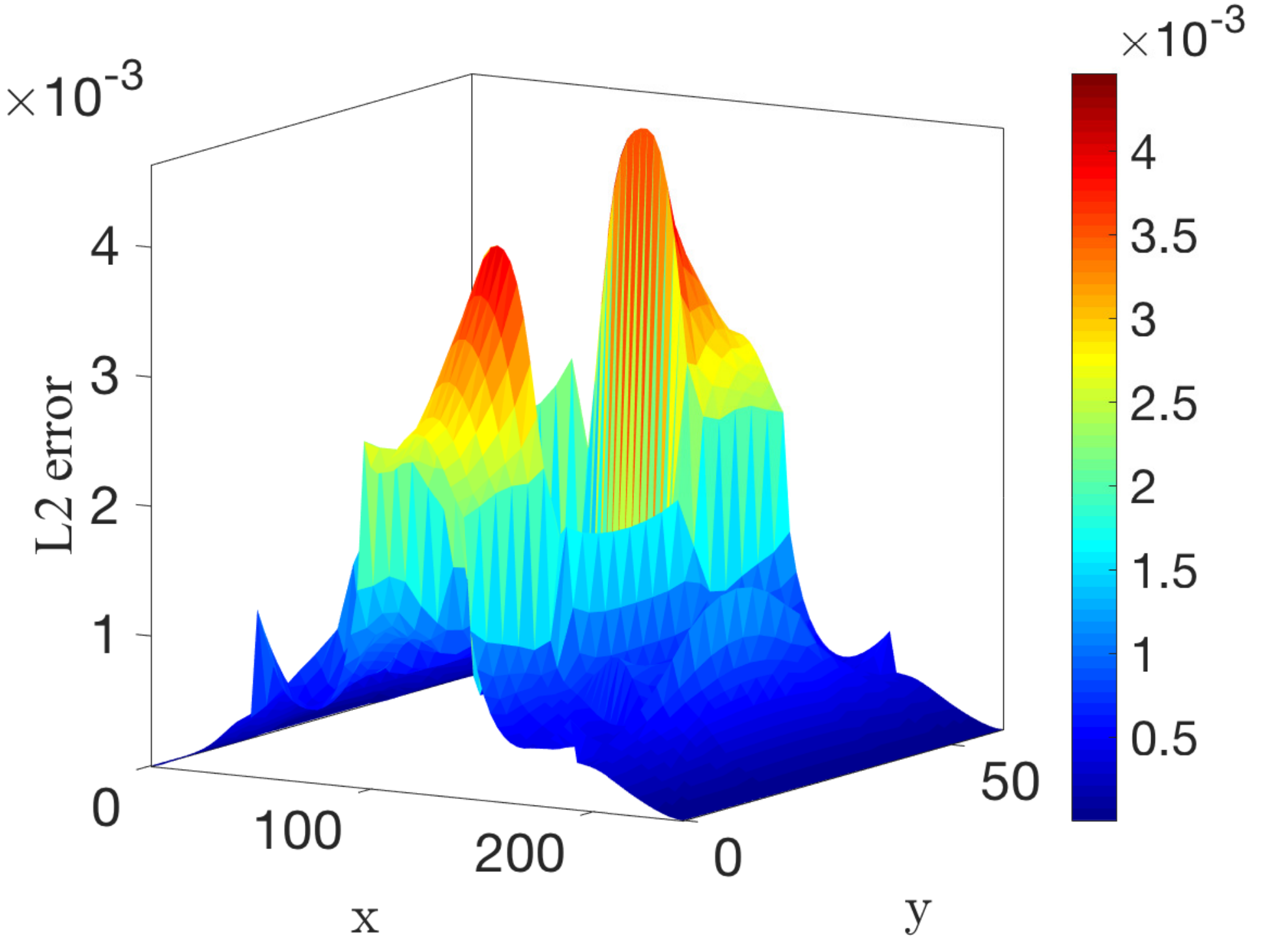}
        \caption{L2 error, $\eta, r=4$} \label{RK:fig:u_L2_xid10_p3_gq5_src0_sink1_etad4_p3_15DOM_rel_error}
    \end{subfigure}    
    \begin{subfigure}[t]{0.3\textwidth}
        \centering
        \includegraphics[height=1.2in]{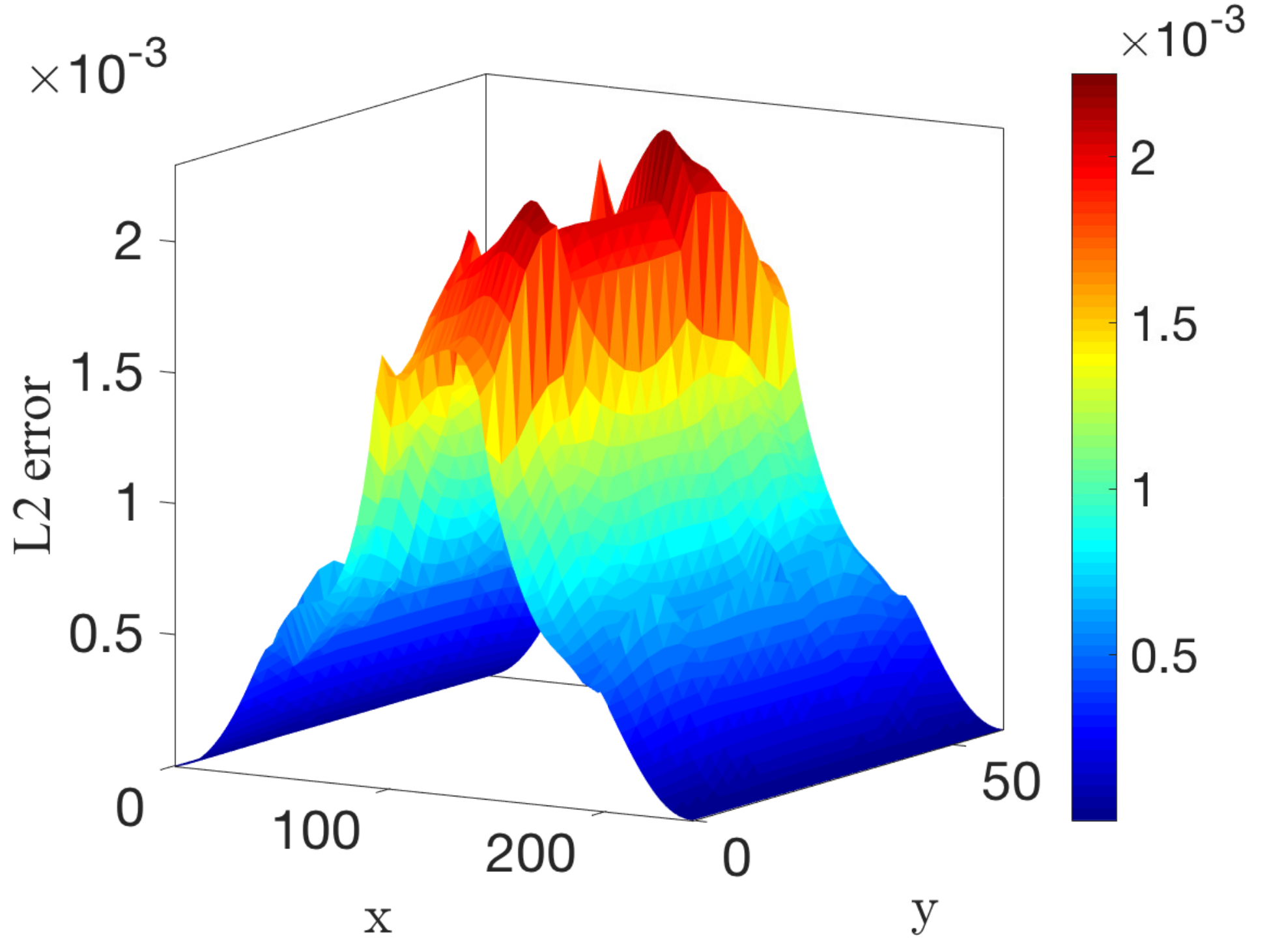}
        \caption{L2 error, $\eta, r=5$} \label{RK:fig:u_L2_xid10_p3_gq5_src0_sink1_etad5_p3_15DOM_rel_error}
    \end{subfigure}  
   
    \caption{L2 norm of the error ($\epsilon(x) = \mathbb{E}[(u(x,\boldsymbol{\xi})-u(x,\boldsymbol{\eta}))^2]$) in the solution obtained by stochastic basis adaptation and domain decomposition (15 subdomains) with random variables $\eta$ and dimension,  $r$ = 3, 4 and 5, order, $p$ = 3, sparse-grid level, $l=5$ and the dimension of reference solution in $\xi$ is $d=10$.} \label{RK:fig:uu_L2_xid10_p3_gq5_src0_D15}
\end{figure}

Fig. \ref{RK:fig:u_mean_xid10_etad5} shows the mean of the reference solution ($d=10$) and mean of the low-dimensional solution computed with $r =5$ with domain decomposition and basis adaptation, as well as the corresponding relative error. Fig. \ref{RK:fig:u_sdev_xid10_etad5} shows the corresponding plots for the standard deviation of the solution. There is good agreement between the reference and low-dimensional solutions. The mean and standard deviation correspond to the first and second moments of the solution. Another metric to compare the solution is to compute the expectation of the error between the two solutions, that is $\epsilon(x) = \mathbb{E}[(u(x,\boldsymbol{\xi})-u(x,\boldsymbol{\eta}))^2]$ which gives a measure of the error in all the statistics of the solutions, instead of just the first and second moments. We computed the relative error with respect to the solution $u(x,\boldsymbol{\xi})$ in Fig. \ref{RK:fig:uu_L2_xid10_p3_gq5_src0_D3}, \ref{RK:fig:uu_L2_xid10_p3_gq5_src0_D8} and \ref{RK:fig:uu_L2_xid10_p3_gq5_src0_D15} for subdomain number $D$ of 3, 8 and 15 respectively.

\begin{figure}[t!]
    \centering
    \begin{subfigure}[t]{0.45\textwidth}
        \centering
        \includegraphics[height=1.95in]{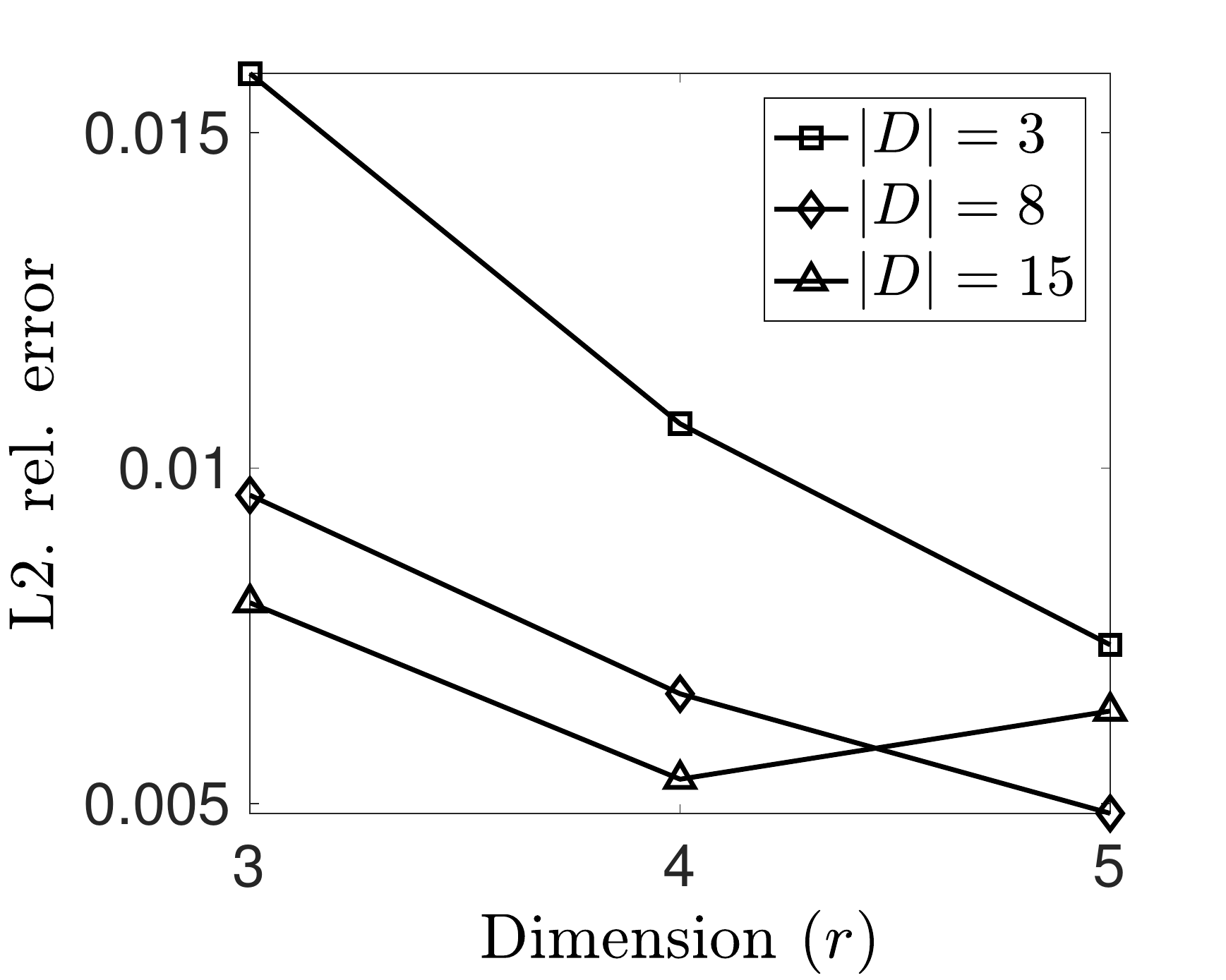}
        \caption{Mean} \label{RK:fig:mean_xid10_p3_gq5_src0_sink1_eta_p3_L2_sqrtn_rel_error_plot}
    \end{subfigure}        
    \begin{subfigure}[t]{0.45\textwidth}
        \centering
        \includegraphics[height=1.95in]{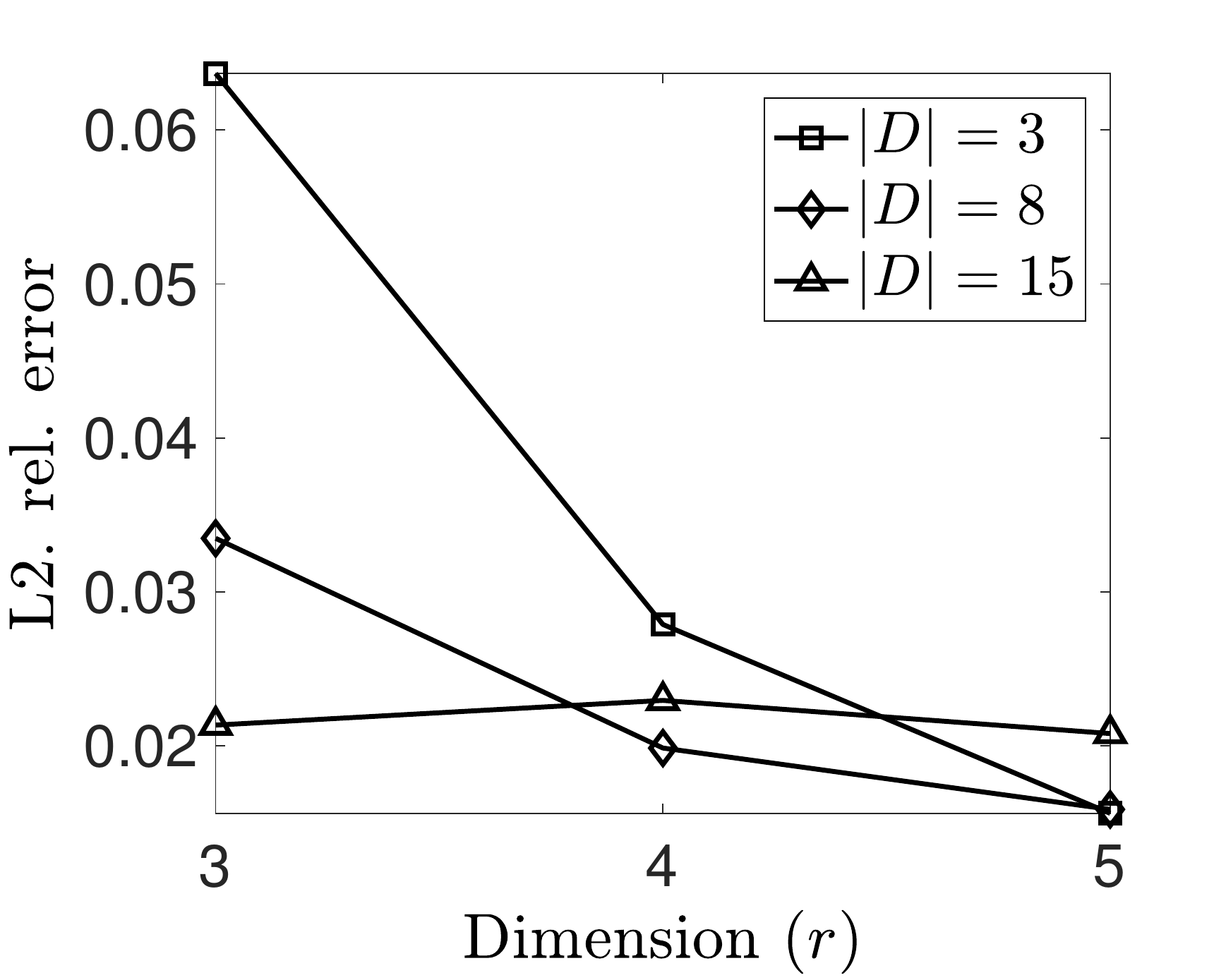}
        \caption{Standard deviation} \label{RK:fig:sdev_xid10_p3_gq5_src0_sink1_eta_p3_L2_sqrtn_rel_error_plot}
    \end{subfigure}    
     \begin{subfigure}[t]{0.45\textwidth}
        \centering
        \includegraphics[height=1.95in]{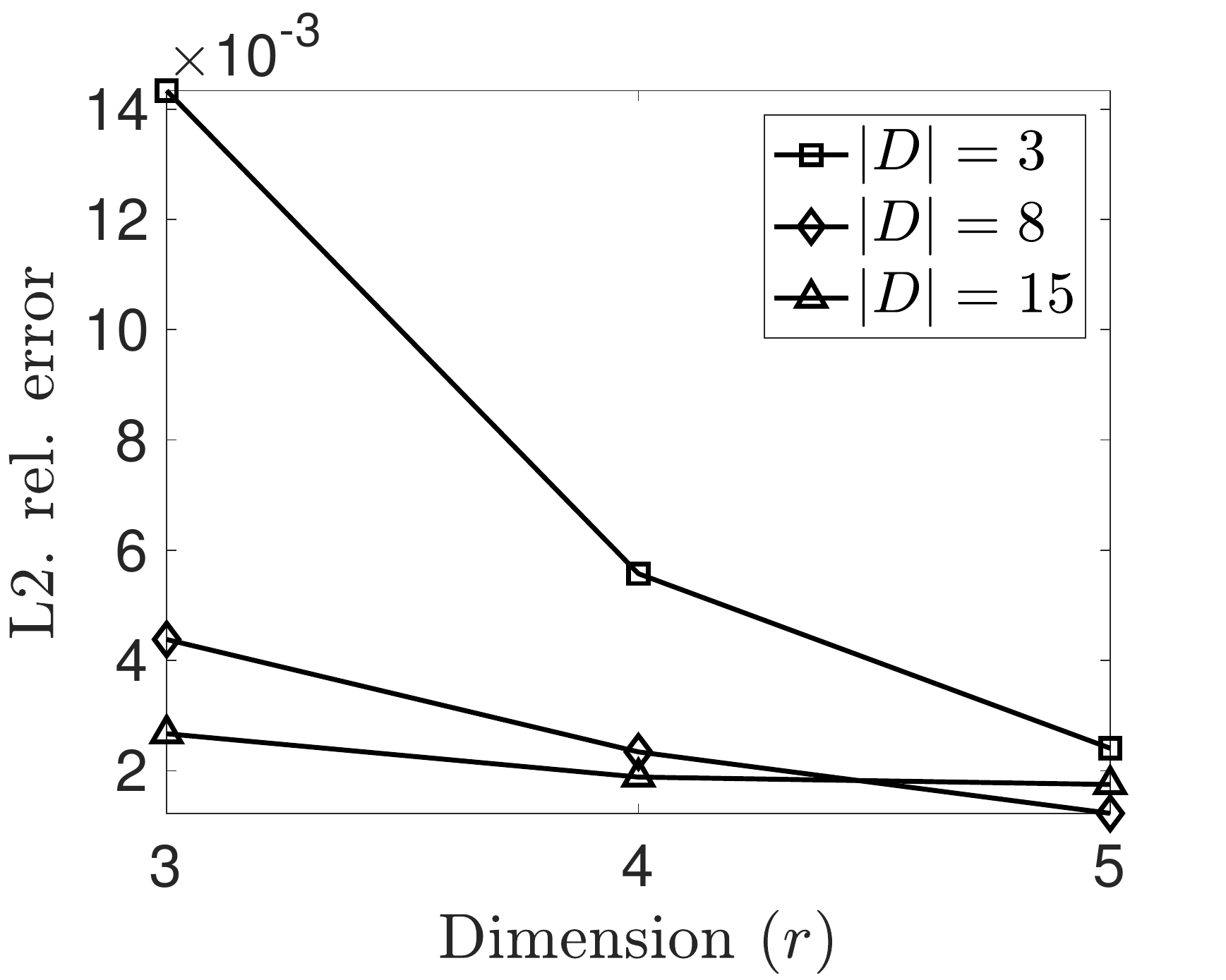}
        \caption{$\epsilon = \mathbb{E} [\|u(x,\boldsymbol{\xi})-u(x,\boldsymbol{\eta})\|^2_2]$} \label{RK:fig:u_L2_error_gauss_xid10_eta_d3-4-5_D3-8-15}
    \end{subfigure}   
    \caption{L2 relative error in mean and standard deviation of the solution obtained by stochastic basis adaptation and domain decomposition (3, 8 and 15 subdomains)  with random variables $\eta$ and dimension, $r$ = 3, 4 and 5, order, $p$ = 3 and sparse-grid level, $l=5$. The reference solution is computed in 10-dimensional random variables in $\xi$ using sparse-grid collocation method with level 5.} \label{RK:fig:u_mean_L2_error_plot}
\end{figure}

Fig. \ref{RK:fig:mean_xid10_p3_gq5_src0_sink1_eta_p3_L2_sqrtn_rel_error_plot} and Fig.  \ref{RK:fig:sdev_xid10_p3_gq5_src0_sink1_eta_p3_L2_sqrtn_rel_error_plot} show the respective relative error for the mean and standard deviation as a function of the number of dimensions $r$ in the reduced model. In Fig. \ref{RK:fig:u_L2_error_gauss_xid10_eta_d3-4-5_D3-8-15}, we show the relative error norm ($\epsilon = \mathbb{E} [\|u(x,\boldsymbol{\xi})-u(x,\boldsymbol{\eta})\|^2_2]$) as a function of the number of dimensions. In general, the relative error decreases with increasing $r$ and number of subdomains.

However, we observe that the error slightly increases for 15 subdomains and $r=5$. There are two sources of error: one due to the use of a reduced number of stochastic dimensions and the other due to the number of subdomains. Although an increase in the number of subdomains improves the accuracy of the local basis adaptation, there also can be an adverse effect on solution accuracy if there are too many subdomains. The increase in error is caused by the projection of components of far away subdomains on the local subdomain (see Step 5 in Algorithm 1). In the particular example where we observe the slight increase in error, it stems from using 15 subdomains for case with a relatively small $d$ ($d=10$). A more thorough examination of the trade-off between the number of reduced stochastic dimensions and the number of subdomains will be published elsewhere.

\begin{figure}[t!]
    \centering
    \begin{subfigure}[t]{0.32\textwidth}
        \centering
        \includegraphics[height=1.2in]{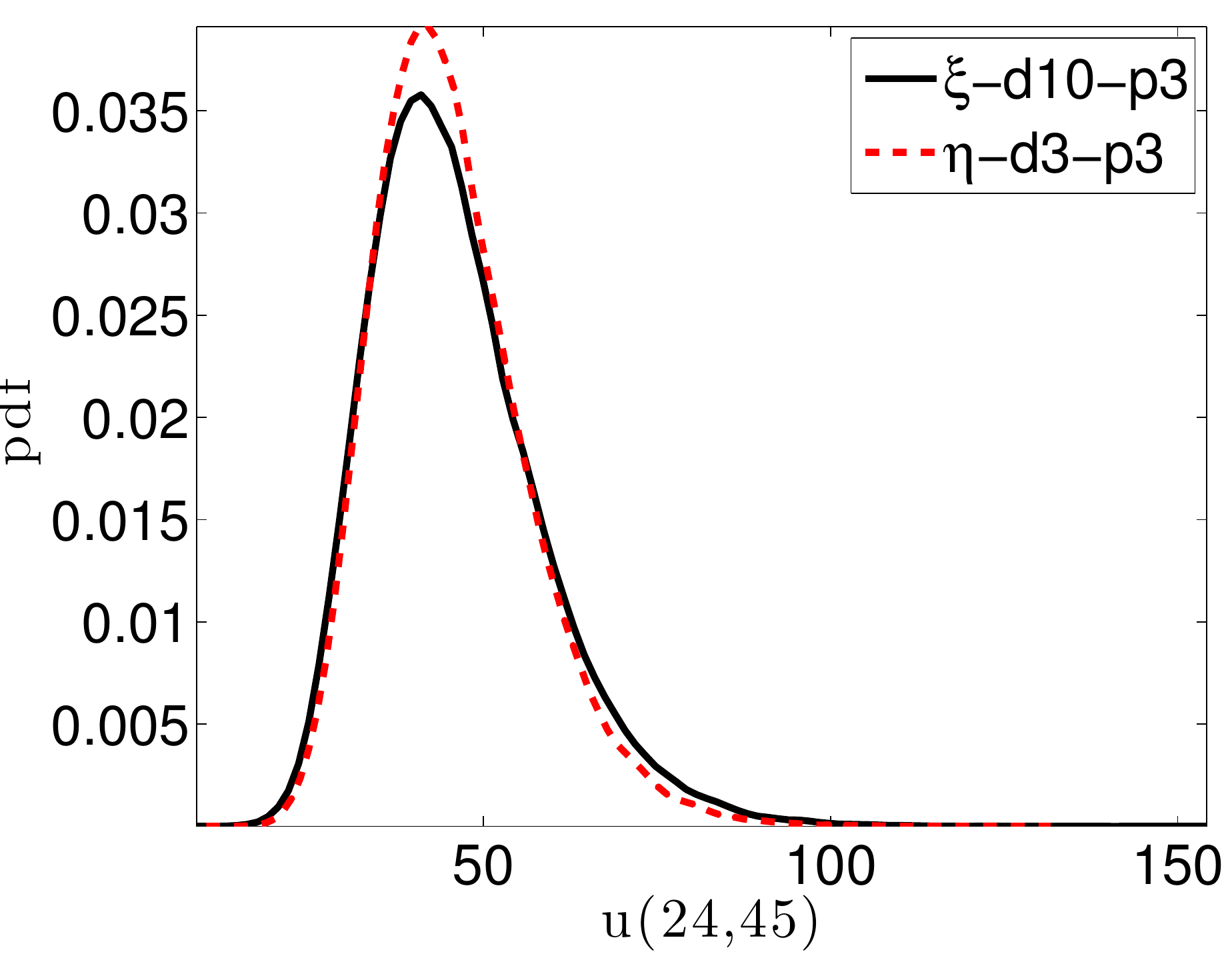}
        \caption{$\eta, d=3,  |D|=3$} \label{RK:fig:pdf_x24_y45_xid10_p3_gq5_src0_sink1_etad3_p3_3DOM}
    \end{subfigure}        
    \begin{subfigure}[t]{0.32\textwidth}
        \centering
        \includegraphics[height=1.2in]{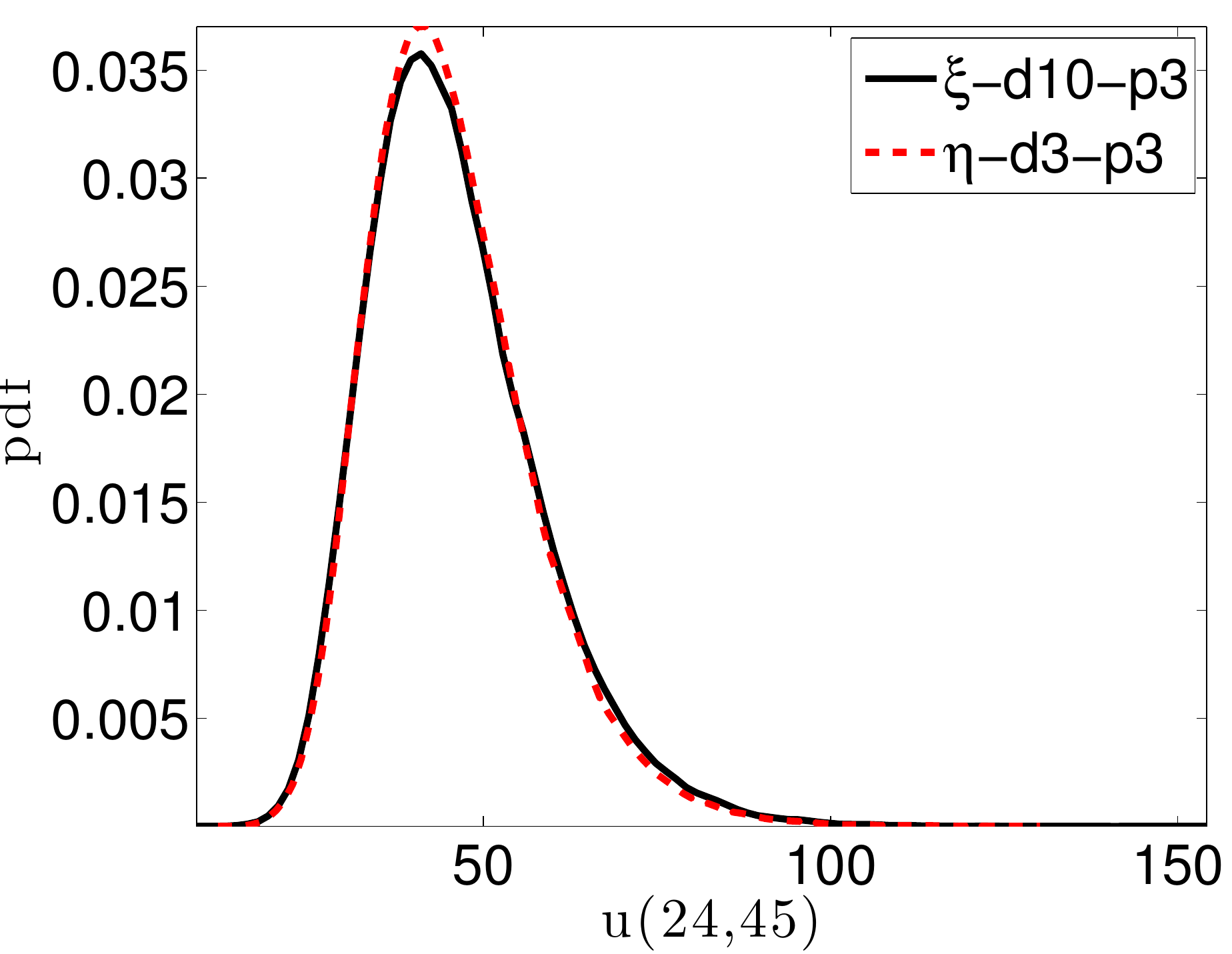}
        \caption{$\eta, d=3,  |D|=8$} \label{RK:fig:pdf_x24_y45_xid10_p3_gq5_src0_sink1_etad3_p3_8DOM}
    \end{subfigure}    
    \begin{subfigure}[t]{0.32\textwidth}
        \centering
        \includegraphics[height=1.2in]{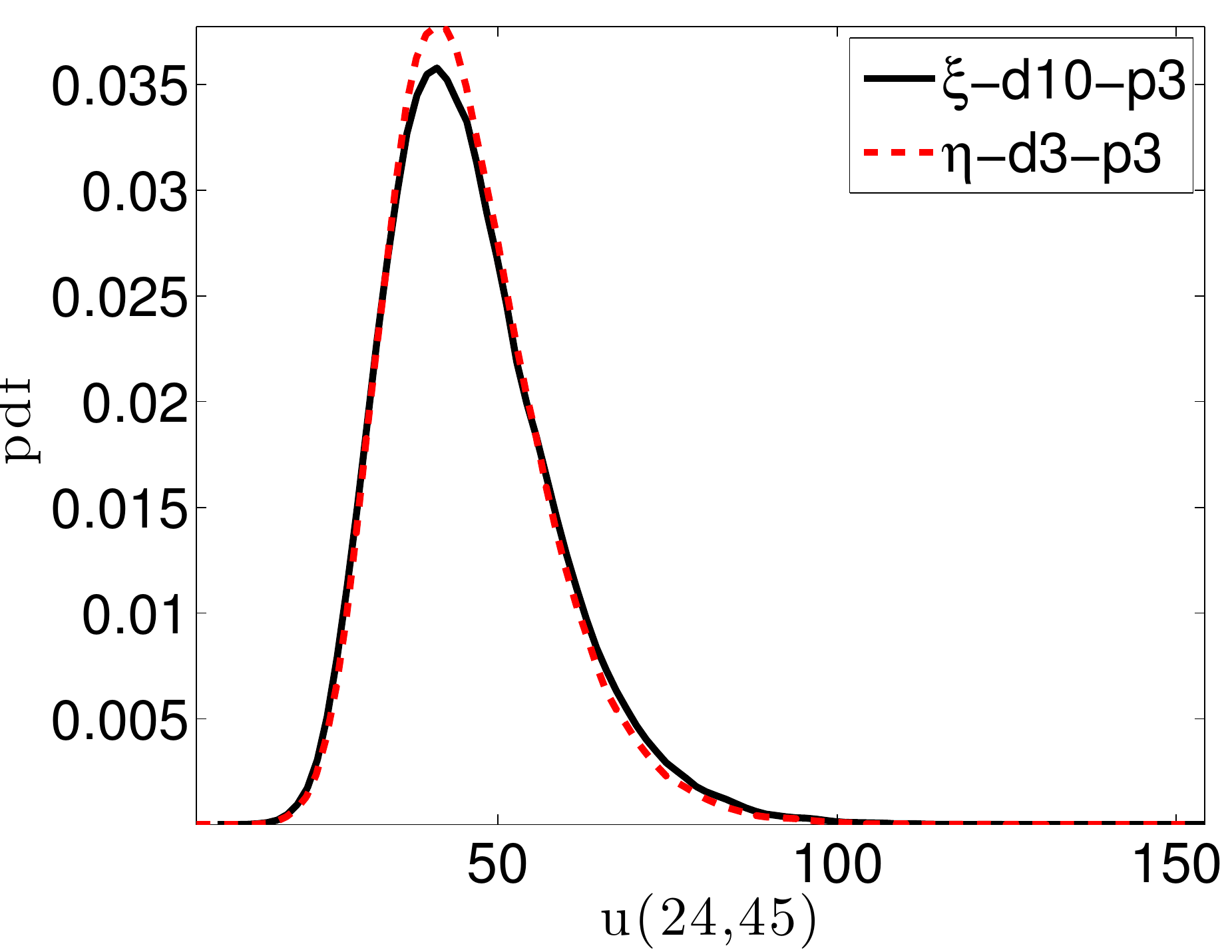}
        \caption{$\eta, d=3,  |D|=15$} \label{RK:fig:pdf_x24_y45_xid10_p3_gq5_src0_sink1_etad3_p3_15DOM}
    \end{subfigure}  
    \begin{subfigure}[t]{0.32\textwidth}
        \centering
        \includegraphics[height=1.2in]{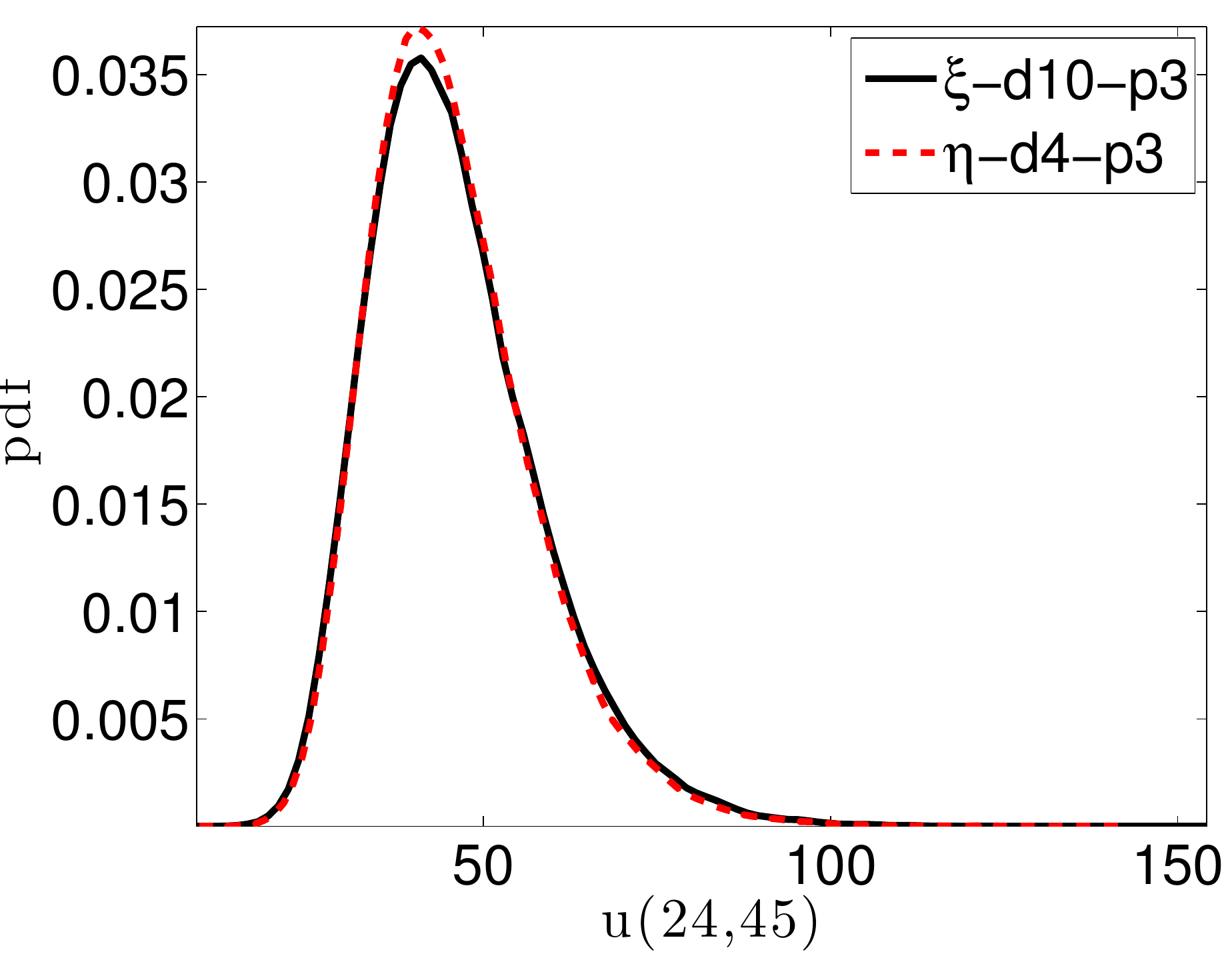}
        \caption{$\eta, d=4,  |D|=3$} \label{RK:fig:pdf_x24_y45_xid10_p3_gq5_src0_sink1_etad4_p3_3DOM}
    \end{subfigure}        
    \begin{subfigure}[t]{0.32\textwidth}
        \centering
        \includegraphics[height=1.2in]{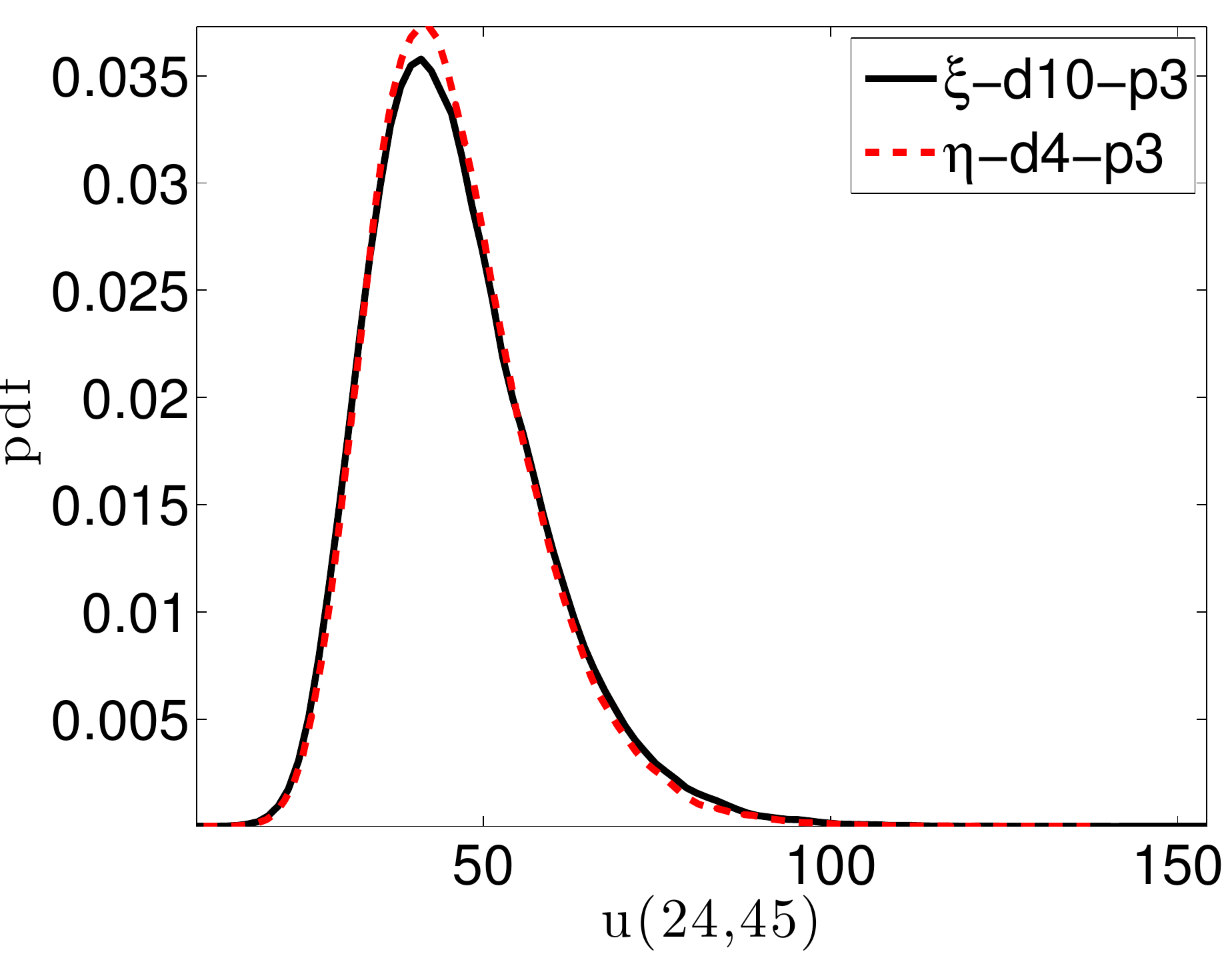}
        \caption{$\eta, d=4,  |D|=8$} \label{RK:fig:pdf_x24_y45_xid10_p3_gq5_src0_sink1_etad4_p3_8DOM}
    \end{subfigure}    
    \begin{subfigure}[t]{0.32\textwidth}
        \centering
        \includegraphics[height=1.2in]{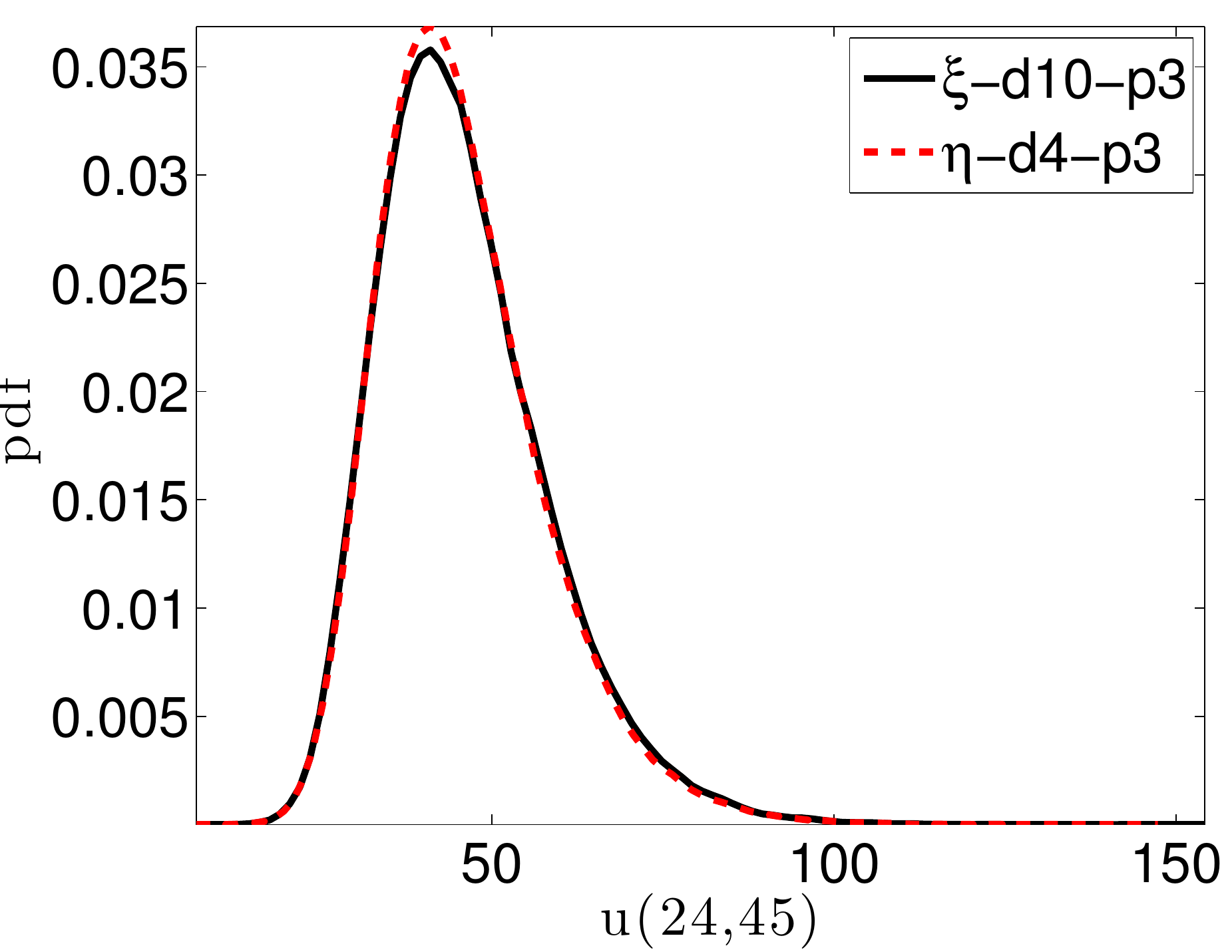}
        \caption{$\eta, d=4,  |D|=15$} \label{RK:fig:pdf_x24_y45_xid10_p3_gq5_src0_sink1_etad4_p3_15DOM}
    \end{subfigure}  
     \begin{subfigure}[t]{0.32\textwidth}
        \centering
        \includegraphics[height=1.2in]{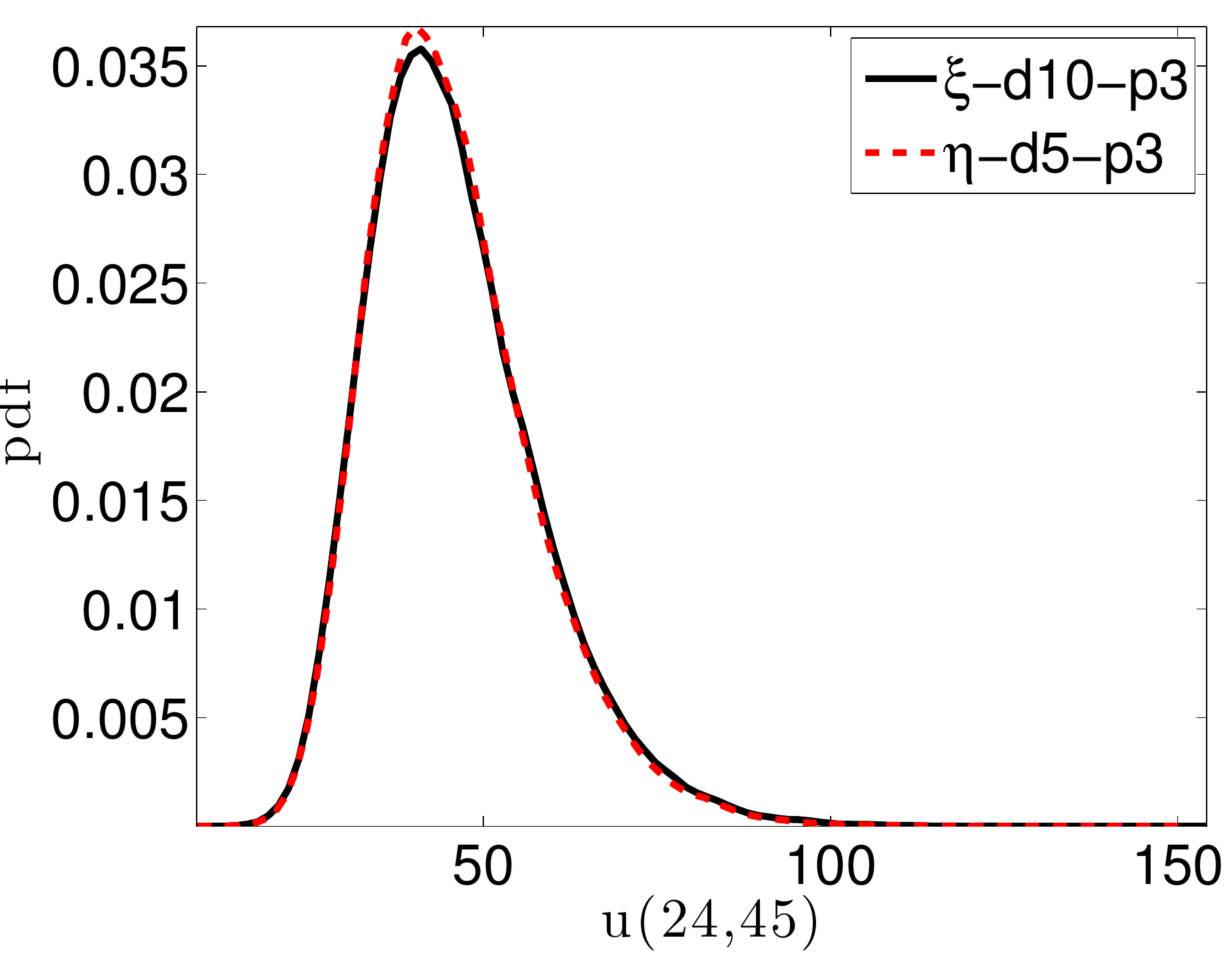}
        \caption{$\eta, d=5,  |D|=3$} \label{RK:fig:pdf_x24_y45_xid10_p3_gq5_src0_sink1_etad5_p3_3DOM}
    \end{subfigure} 
     \begin{subfigure}[t]{0.32\textwidth}
        \centering
        \includegraphics[height=1.2in]{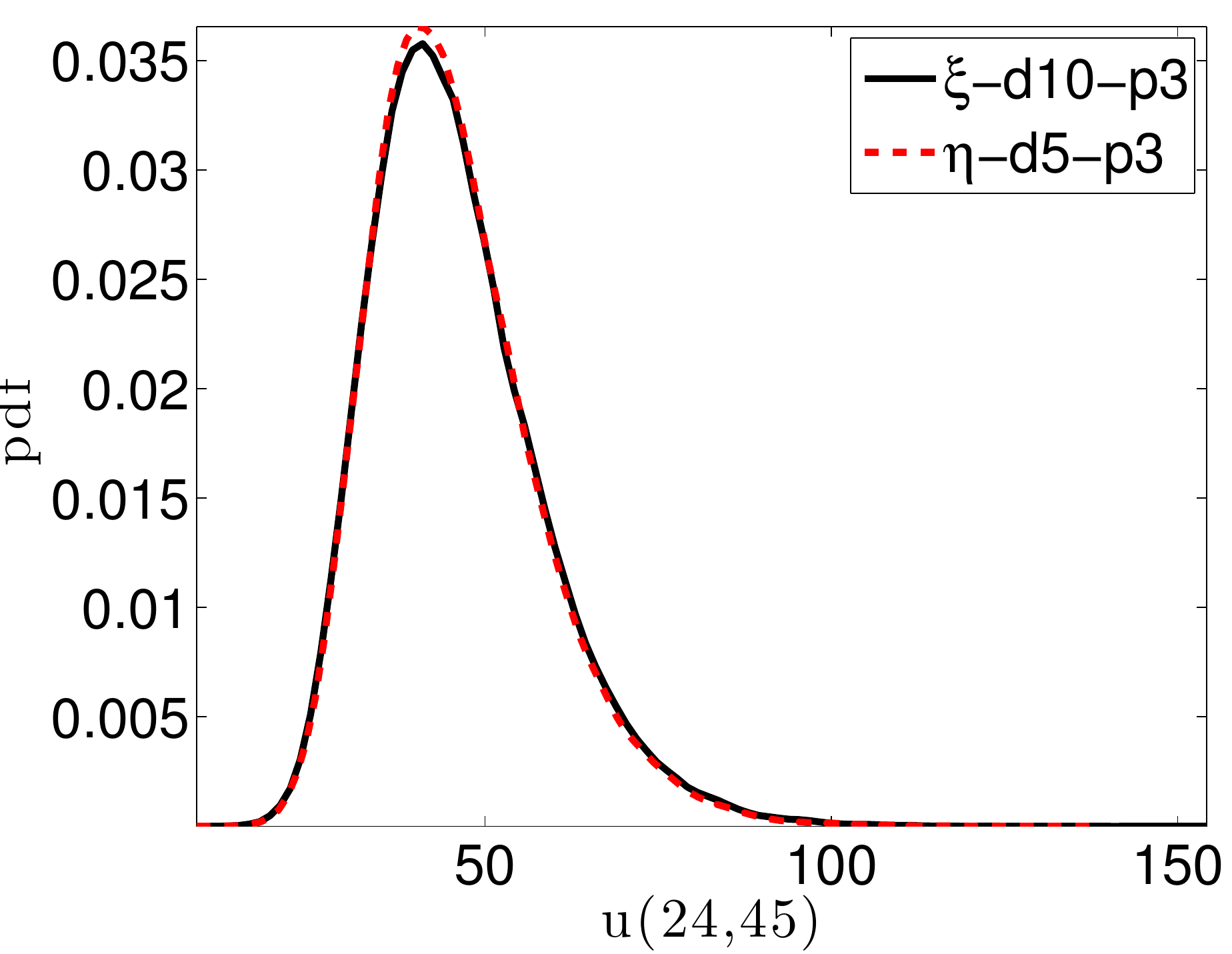}
        \caption{$\eta, d=5,  |D|=8$} \label{RK:fig:pdf_x24_y45_xid10_p3_gq5_src0_sink1_etad5_p3_8DOM}
    \end{subfigure} 
     \begin{subfigure}[t]{0.32\textwidth}
        \centering
        \includegraphics[height=1.2in]{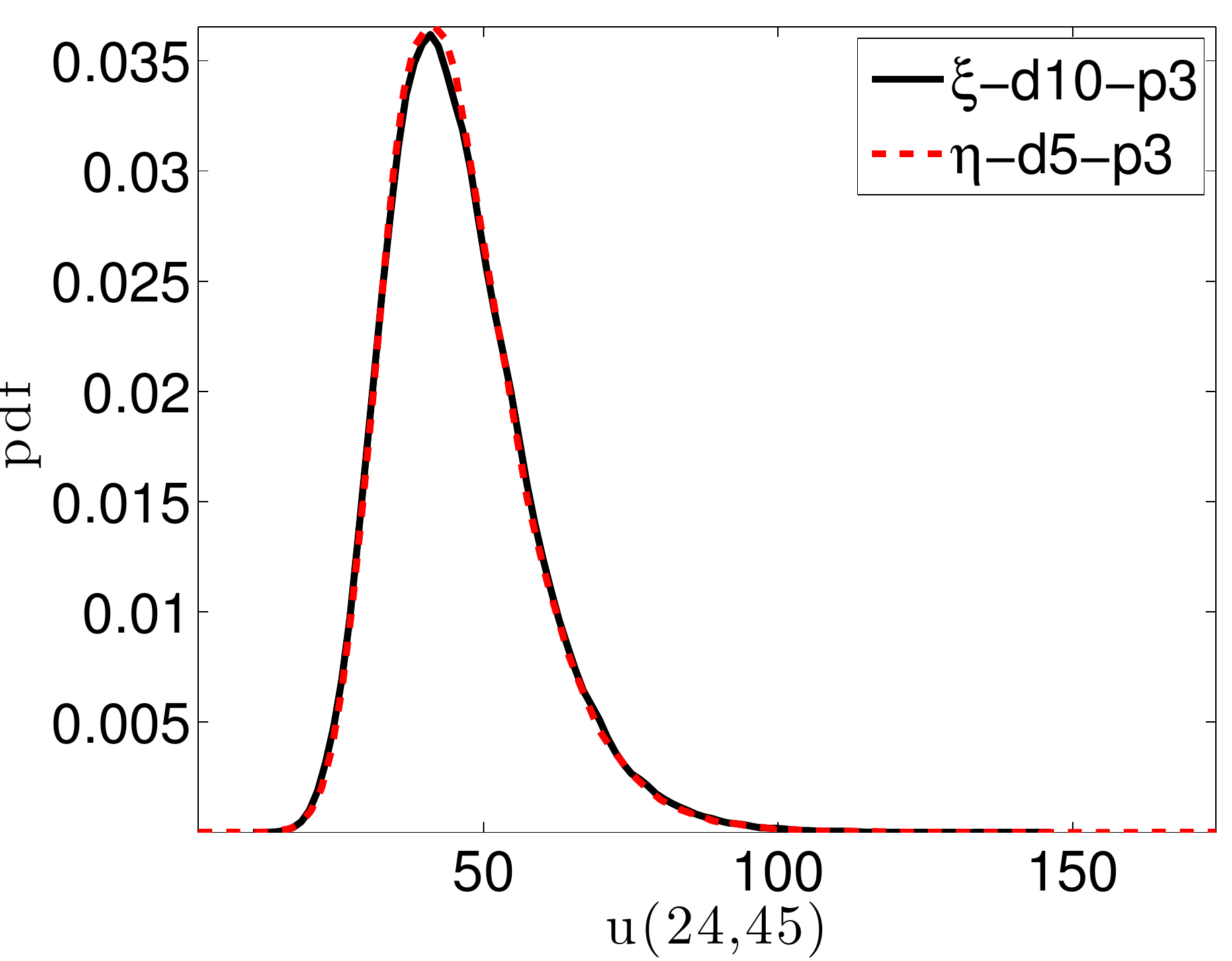}
        \caption{$\eta, d=5,  |D|=15$} \label{RK:fig:pdf_x24_y45_xid10_p3_gq5_src0_sink1_etad5_p3_15DOM}
    \end{subfigure} 
    \caption{Probability density of the solution at ($x=24,y=45$) obtained by stochastic basis adaptation and domain decomposition (3, 8, and 15 subdomains) with random variables $\eta$ and dimension,  $r$ = 3, 4 and 5, order, $p$ = 3 and sparse-grid level, $l=5,$ compared with that from reference solution in $\xi$ of dimension, $d=10$.} \label{RK:fig:pdf_x24_y45_xid10}
\end{figure}

Fig. \ref{RK:fig:pdf_x24_y45_xid10} shows the pdf computed at point $(24, 45)$ using all dimensions ($d=10$) and the low-dimensional representation with various combinations of reduced dimensions and number of subdomains. To compute the reference pdf, we generated 100000 samples. 
We can see good agreement between the reference solution and the solutions obtained with domain decomposition and basis adaptation.

\begin{figure}[t!]
    \centering
    \begin{subfigure}[t]{0.32\textwidth}
        \centering
        \includegraphics[height=1.2in]{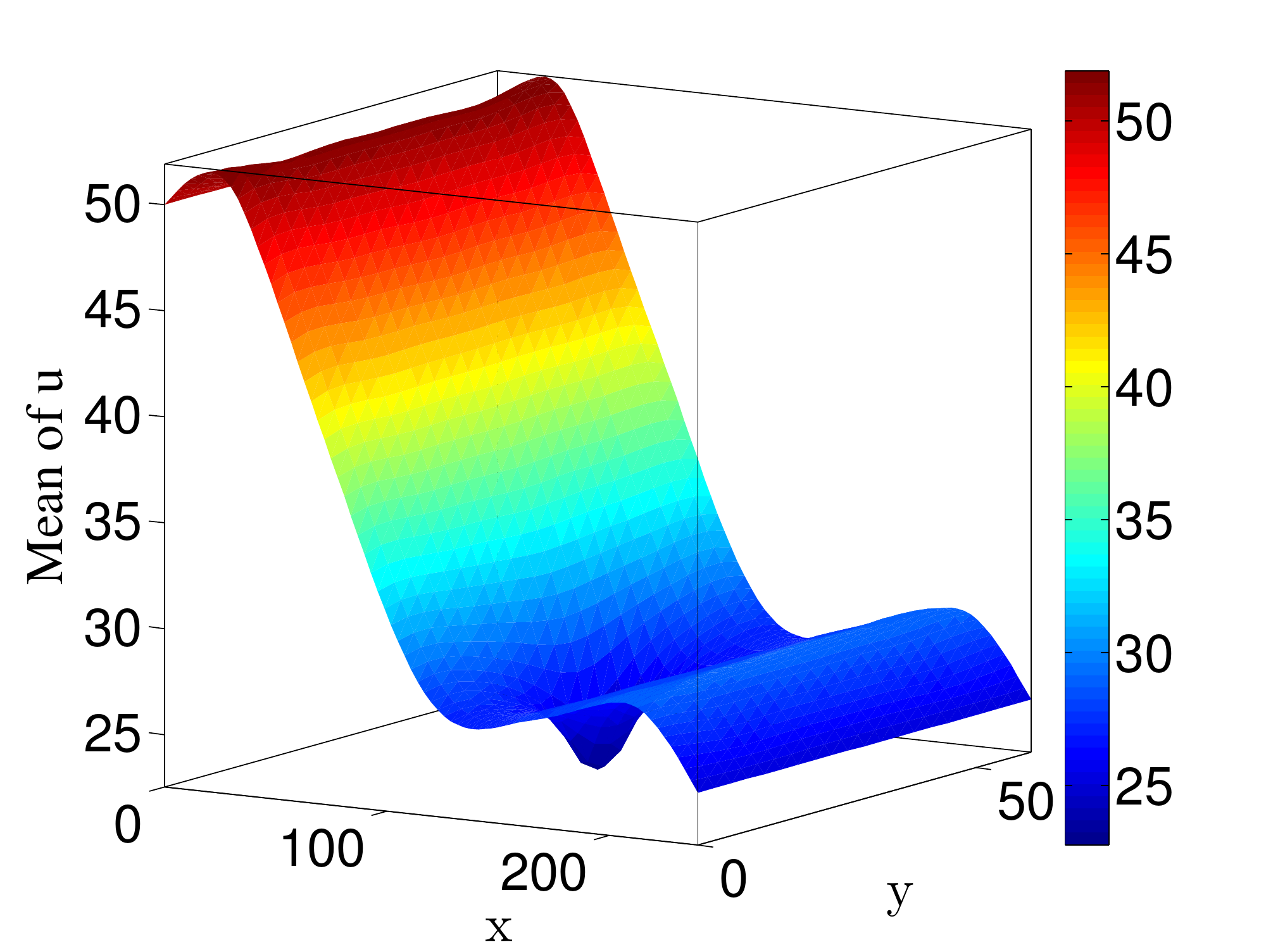}
        \caption{Mean, $\xi, d=40$} \label{RK:fig:u_mean_mc_xi_d40_src0_sink1_mc100000_bcs50_25}
    \end{subfigure}        
    \begin{subfigure}[t]{0.32\textwidth}
        \centering
        \includegraphics[height=1.2in]{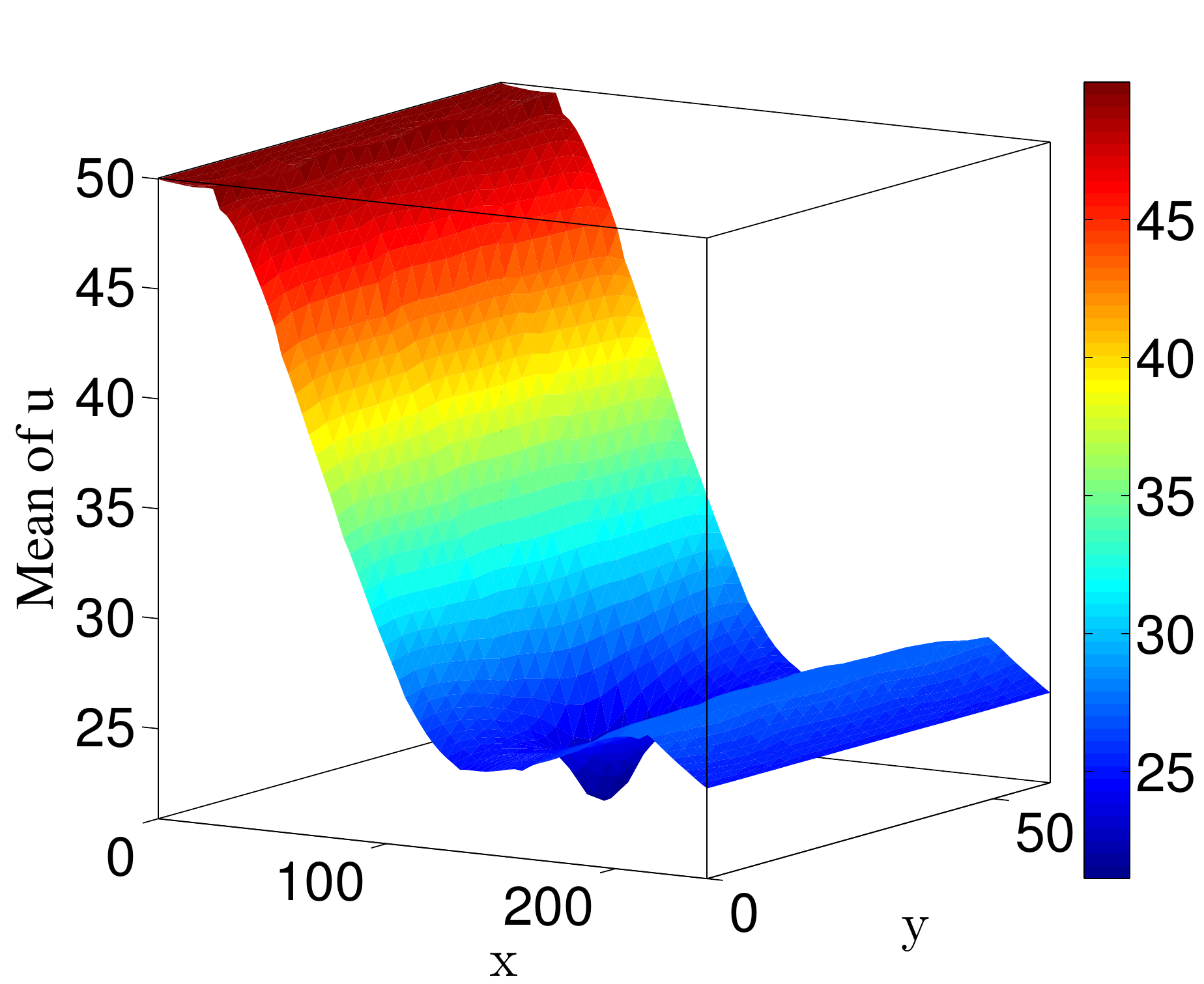}
        \caption{Mean, $\eta, r=4$} \label{RK:fig:u_mean_mc_xi_d40_src0_sink1_mc100000_bcs50_25_etad4_p3_27DOM}
    \end{subfigure}    
    \begin{subfigure}[t]{0.32\textwidth}
        \centering
        \includegraphics[height=1.2in]{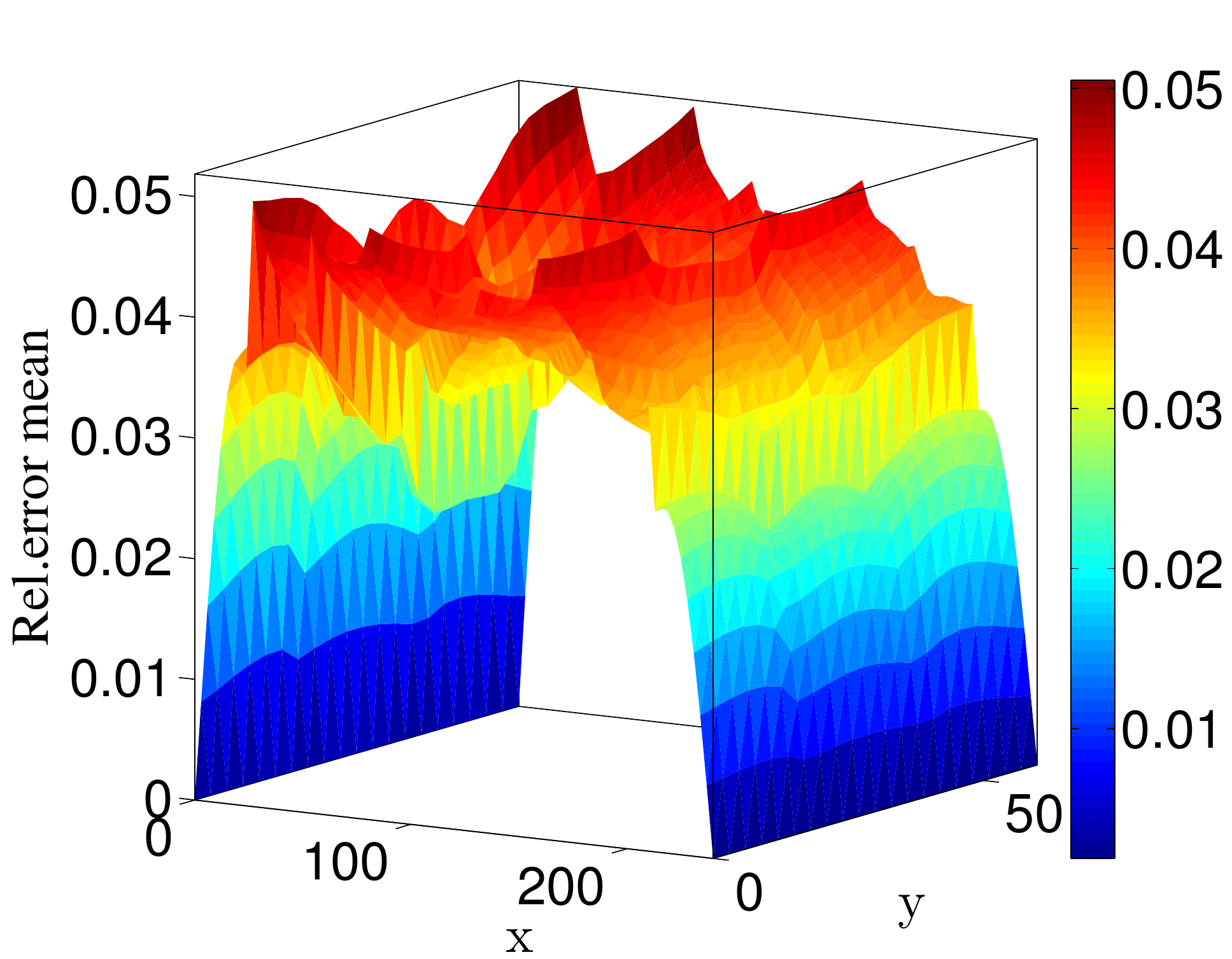}
        \caption{Error, $\eta, r=4$} \label{RK:fig:u_mean_mc_xi_d40_src0_sink1_mc100000_bcs50_25_etad4_p3_27DOM_rel_error}
    \end{subfigure}  
   
    \caption{Mean of the solution obtained by stochastic basis adaptation and domain decomposition (27 subdomains)  with random variables $\eta$ and dimension, $r$ = 4, order, $p$ = 3, sparse-grid level, $l=5$ compared with the reference solution in $\xi$ of dimension, $d=40$.} \label{RK:fig:mean_xid40_etad4_DOM27}
\end{figure}

\begin{figure}[t!]
    \centering
   \begin{subfigure}[t]{0.3\textwidth}
        \centering
        \includegraphics[height=1.2in]{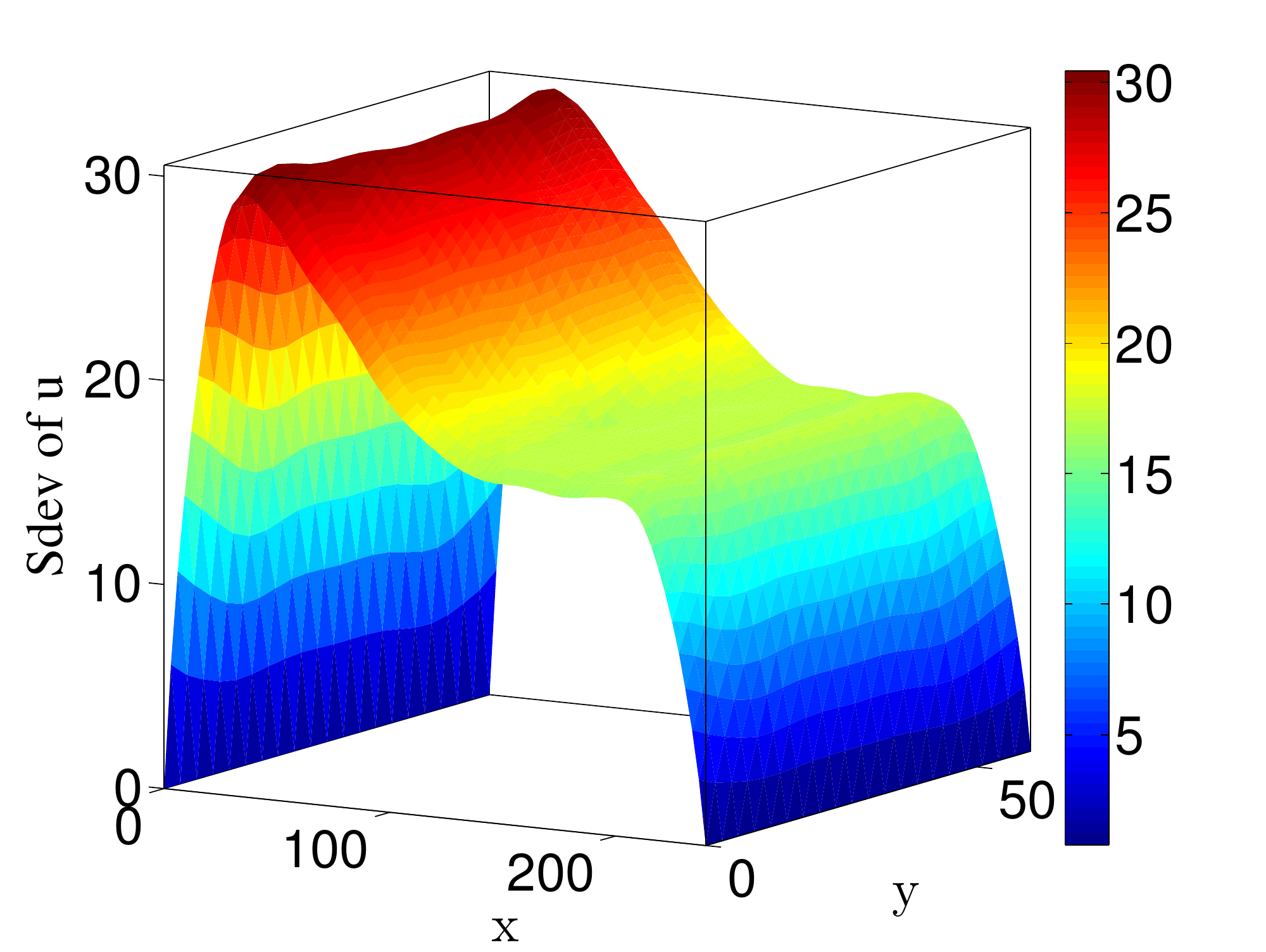}
        \caption{Sdev, $\xi, d=40$} \label{RK:fig:u_sdev_mc_xi_d40_src0_sink1_mc100000_bcs50_25}
    \end{subfigure}        
    \begin{subfigure}[t]{0.3\textwidth}
        \centering
        \includegraphics[height=1.2in]{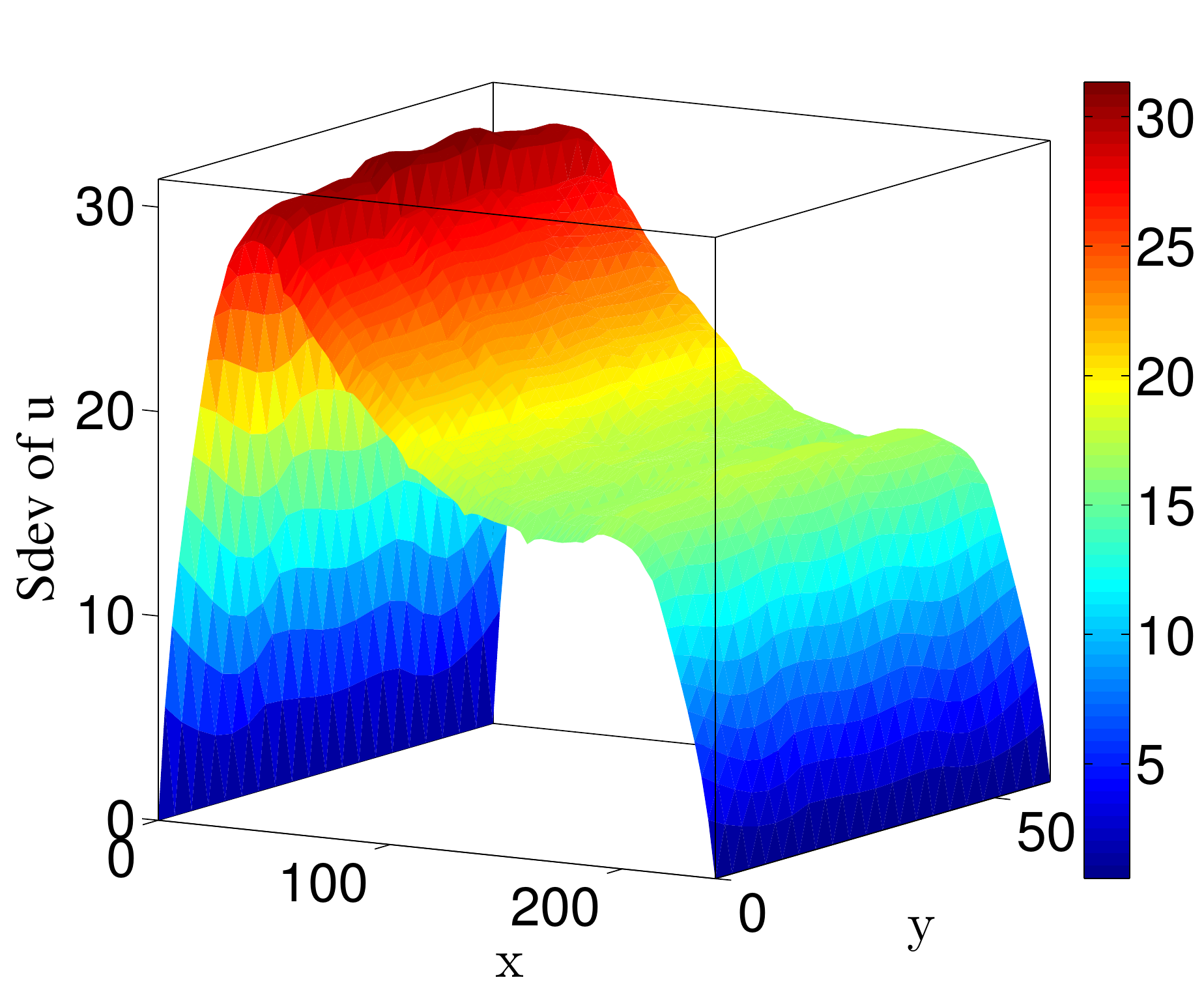}
        \caption{Sdev, $\eta, r=4$} \label{RK:fig:u_sdev_mc_xi_d40_src0_sink1_mc100000_bcs50_25_etad4_p3_27DOM}
    \end{subfigure}    
    \begin{subfigure}[t]{0.3\textwidth}
        \centering
        \includegraphics[height=1.2in]{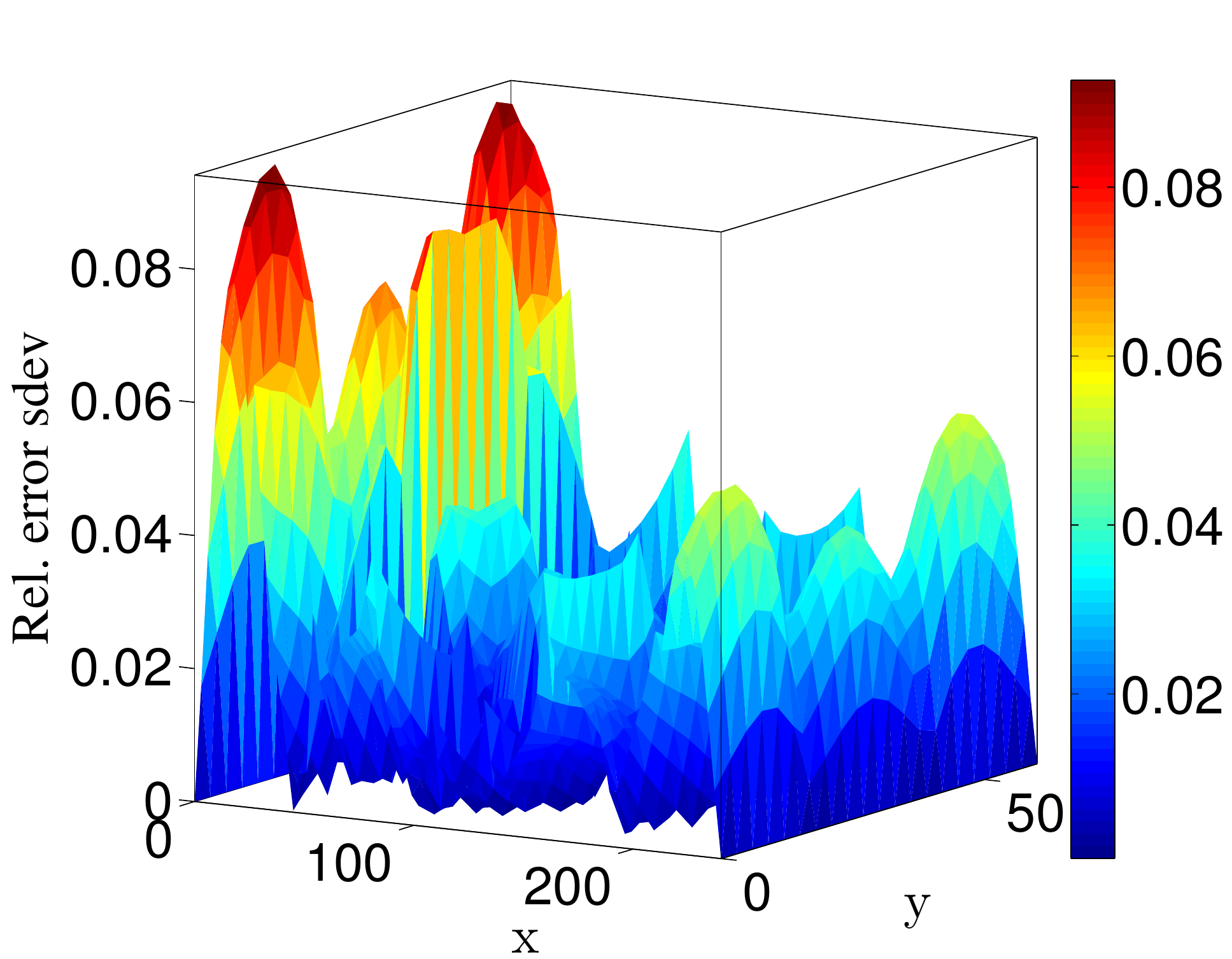}
        \caption{Error, $\eta, r=4$} \label{RK:fig:u_sdev_mc_xi_d40_src0_sink1_mc100000_bcs50_25_etad4_p3_27DOM_rel_error}
    \end{subfigure}  
   
    \caption{Standard deviation of the solution obtained by stochastic basis adaptation and domain decomposition (27 subdomains) with random variables $\eta$ and dimension, $r$ = 4, order, $p$ = 3, sparse-grid level, $l=5$ compared with the reference solution in $\xi$ of dimension, $d=40$.}\label{RK:fig:sdev_xid40_etad4_DOM27}
\end{figure}

\begin{figure}[t!]
    \centering
    \begin{subfigure}[t]{0.45\textwidth}
        \centering
        \includegraphics[height=1.95in]{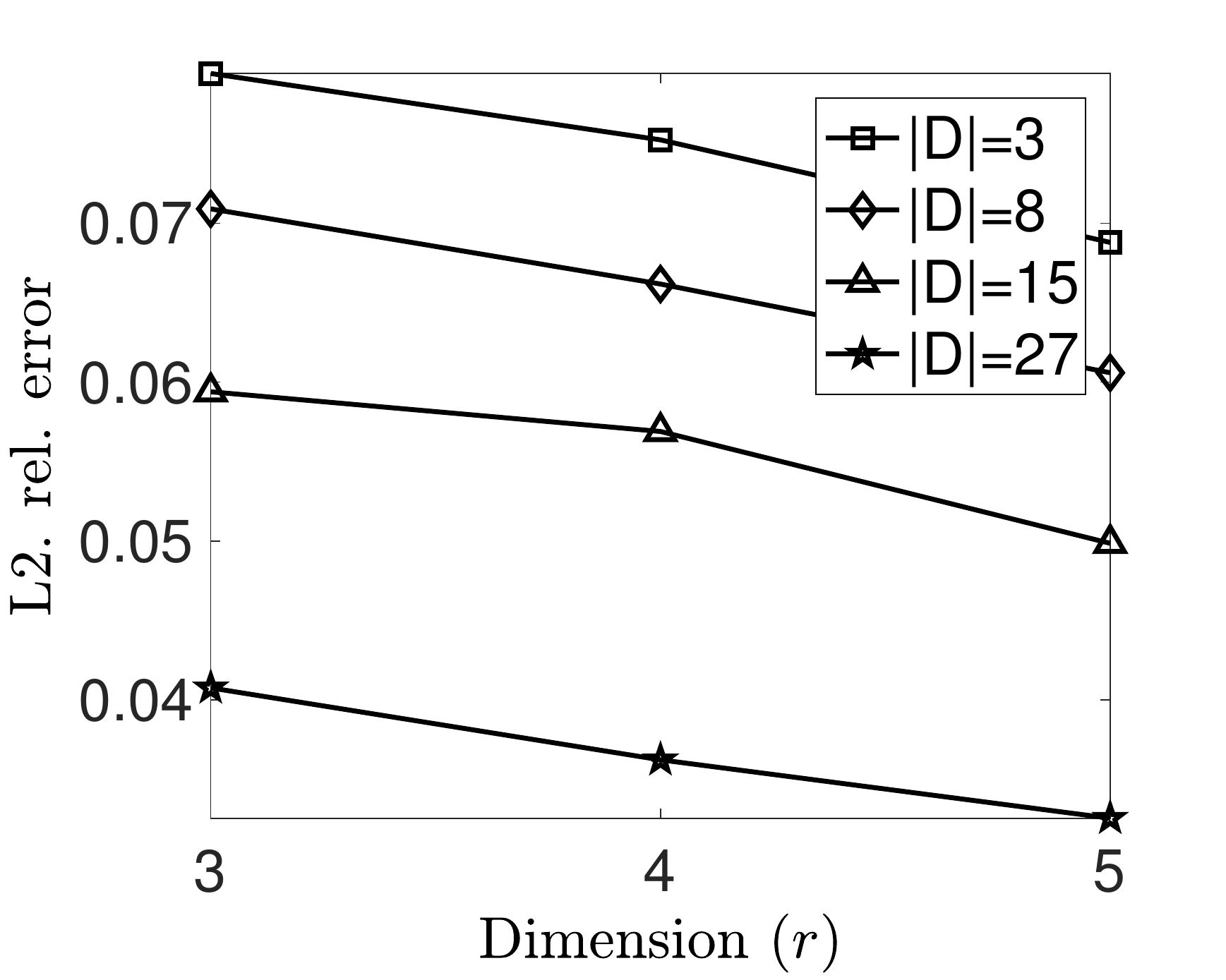}
        \caption{Mean} \label{RK:fig:mean_mc_xi_d40_src0_sink1_mc100000_bcs50_25_etad_p3_L2_sqrtn_rel_error_plot}
    \end{subfigure}        
    \begin{subfigure}[t]{0.45\textwidth}
        \centering
        \includegraphics[height=1.95in]{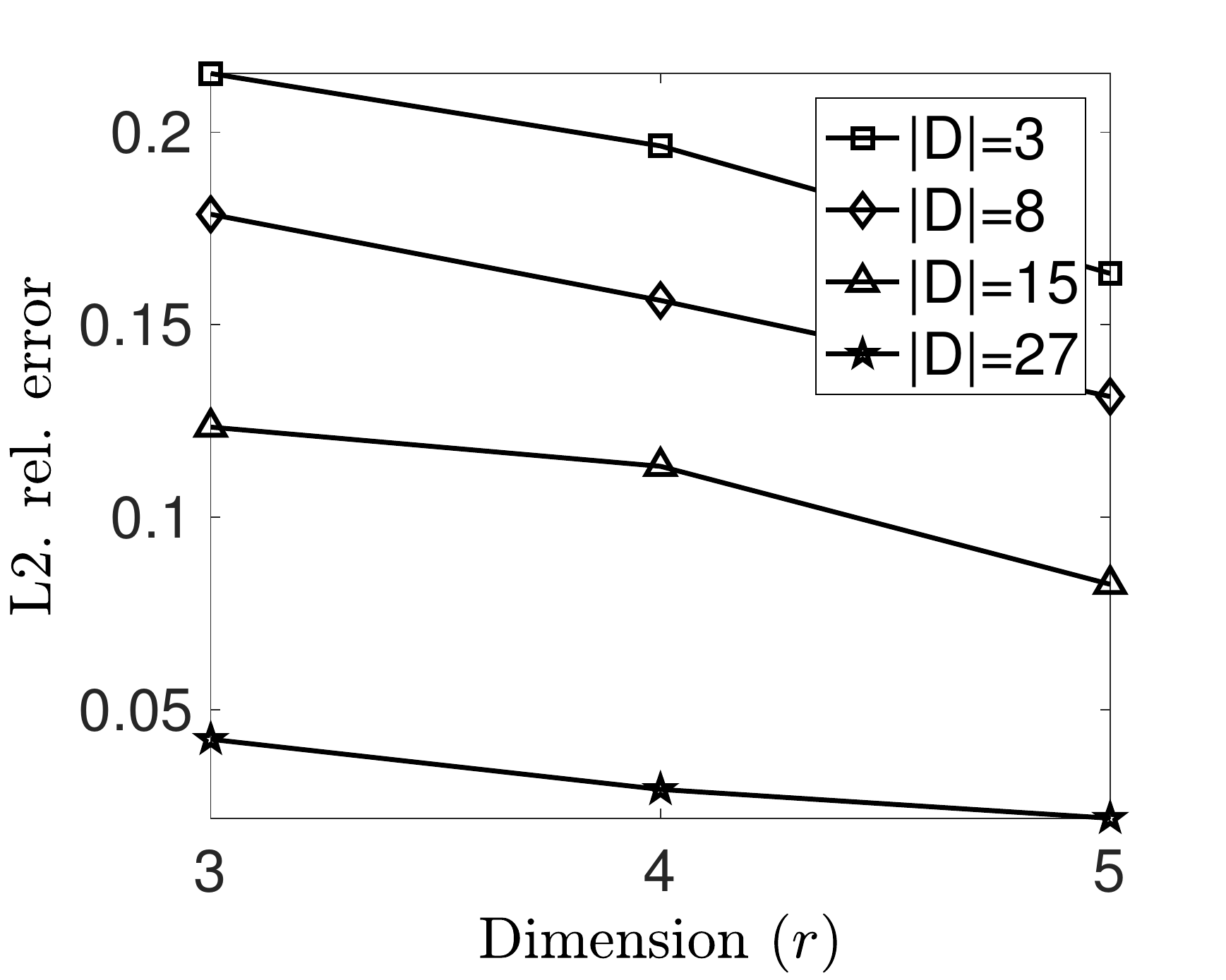}
        \caption{Standard deviation} \label{RK:fig:sdev_mc_xi_d40_src0_sink1_mc100000_bcs50_25_etad_p3_L2_sqrtn_rel_error_plot}
    \end{subfigure}    
    
    \caption{L2 relative error of the mean and standard deviation of the solution obtained by stochastic basis adaptation and domain decomposition (3, 8, 15, and 27 subdomains) with random variables $\eta$ and dimension, $r$ = 3, 4 and 5, order, $p$ = 3 and sparse-grid level, $l=5$. The reference solution is computed with $d=40$ random variables in $\xi$ using 100000 Monte Carlo simulations.} \label{RK:fig:u_mean_xid40_L2_error_plot}
\end{figure}

\begin{figure}[t!]
    \centering
    \begin{subfigure}[t]{0.45\textwidth}
        \centering
        \includegraphics[height=1.8in]{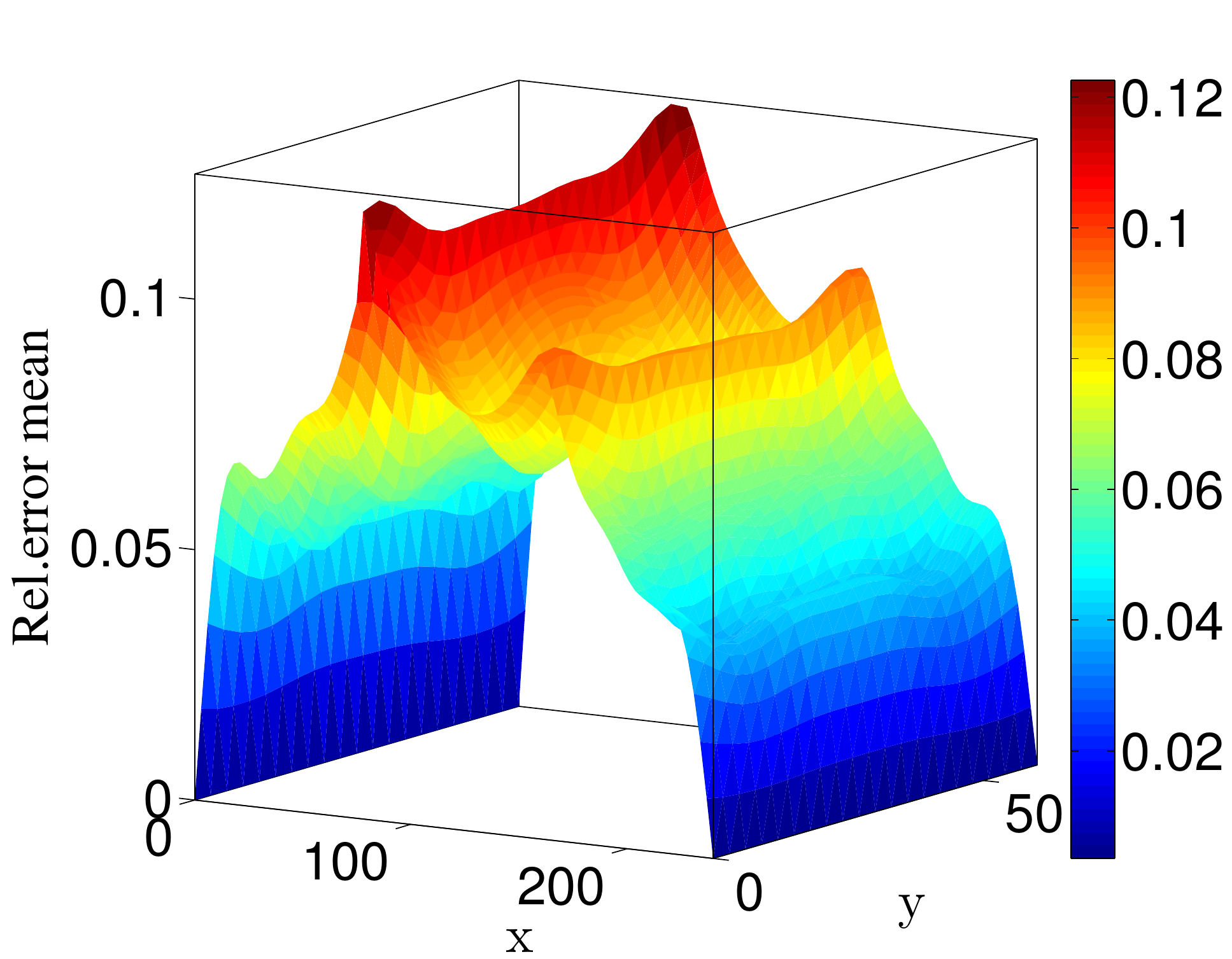}
        \caption{Error, $\xi, d=40; \eta, r=5,  |D|=3$} \label{RK:fig:u_mean_mc_xi_d40_src0_sink1_mc100000_bcs50_25_etad5_p3_3DOM_rel_error}
    \end{subfigure}        
    \begin{subfigure}[t]{0.45\textwidth}
        \centering
        \includegraphics[height=1.8in]{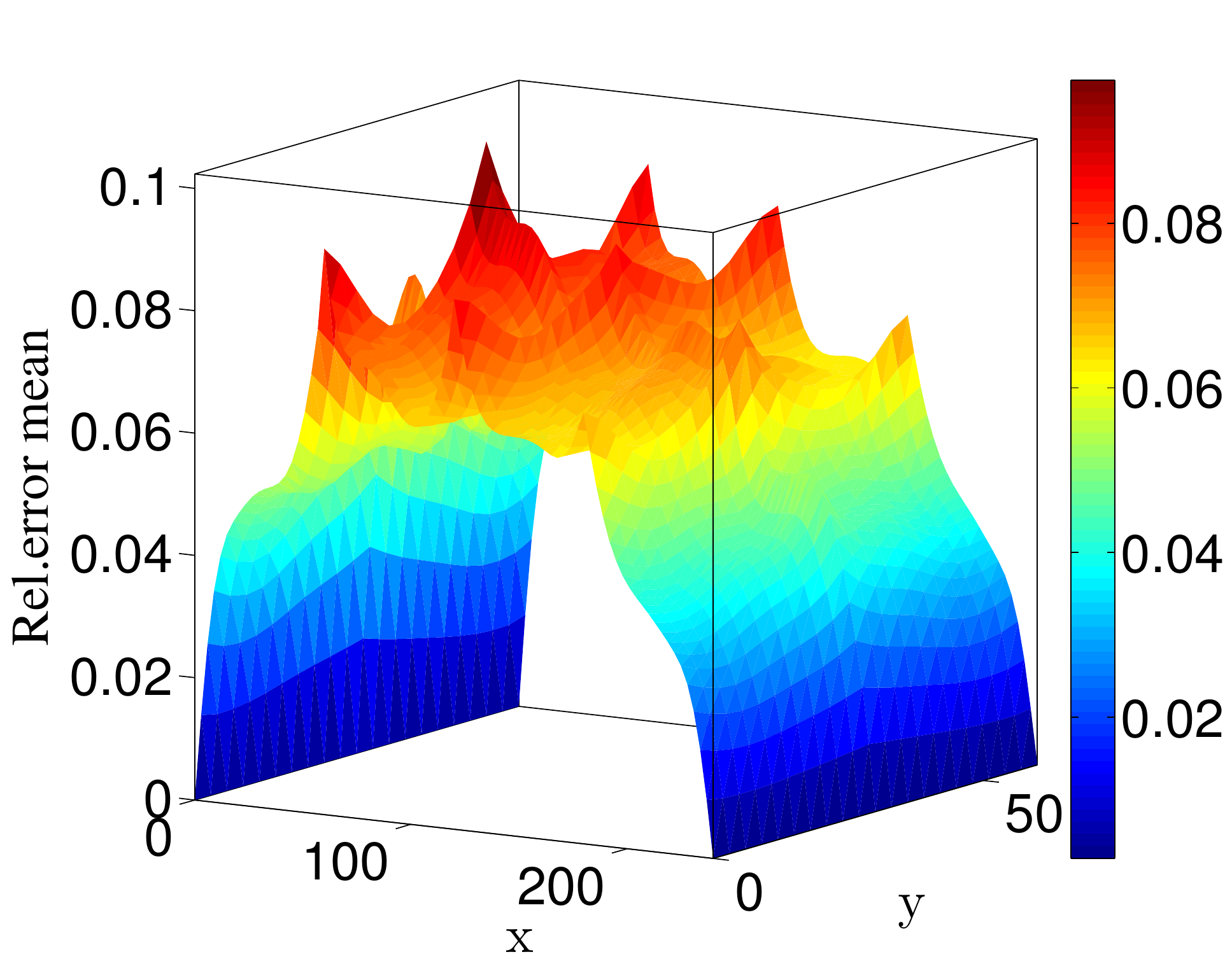}
        \caption{Error, $\xi, d=40; \eta, r=5,  |D|=8$} \label{RK:fig:u_mean_mc_xi_d40_src0_sink1_mc100000_bcs50_25_etad5_p3_8DOM_rel_error}
    \end{subfigure}    
    \begin{subfigure}[t]{0.45\textwidth}
        \centering
        \includegraphics[height=1.8in]{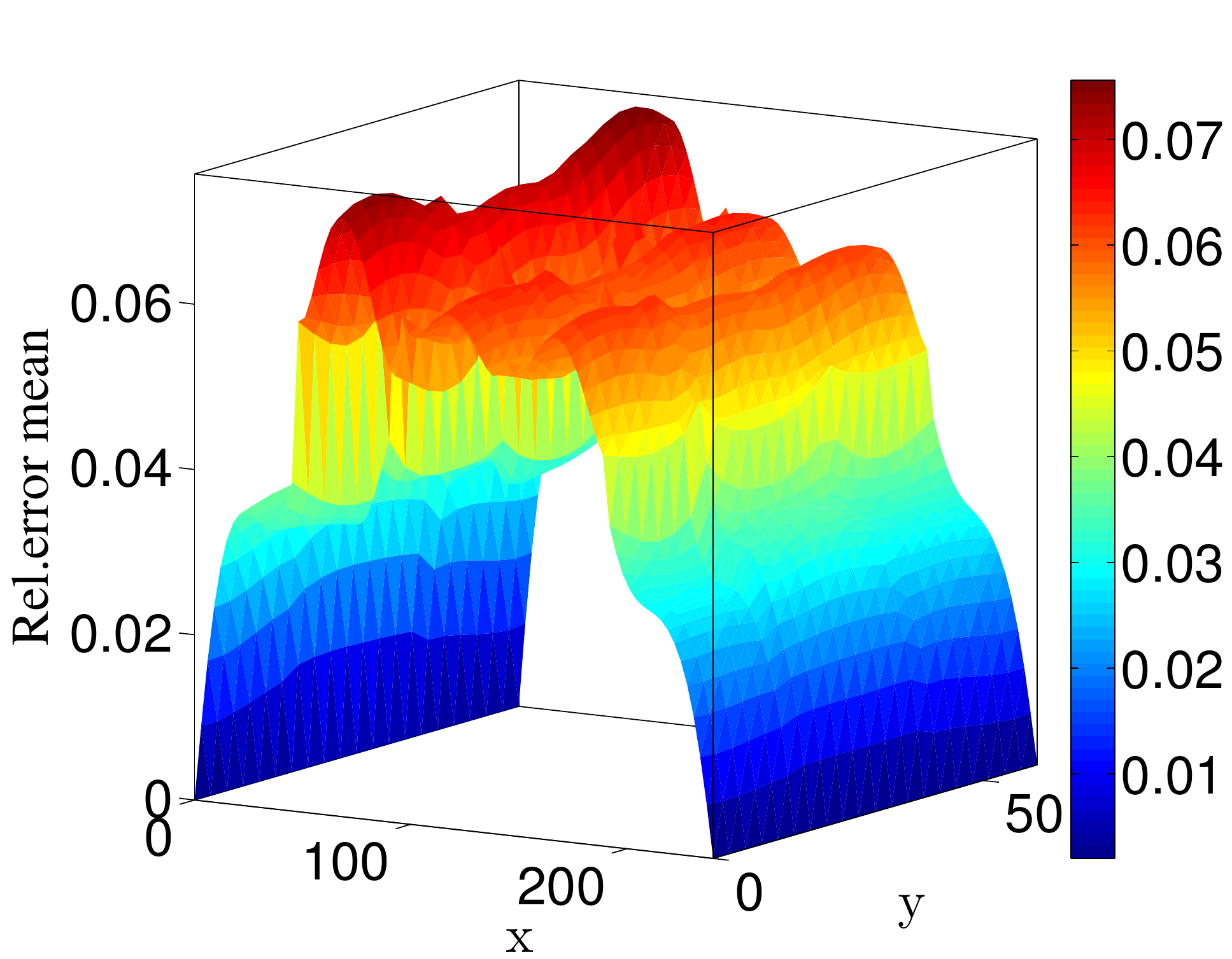}
        \caption{Error, $\xi, d=40; \eta, r=5,  |D|=15$} \label{RK:fig:u_mean_mc_xi_d40_src0_sink1_mc100000_bcs50_25_etad5_p3_15DOM_rel_error}
    \end{subfigure}  
     \begin{subfigure}[t]{0.45\textwidth}
        \centering
        \includegraphics[height=1.8in]{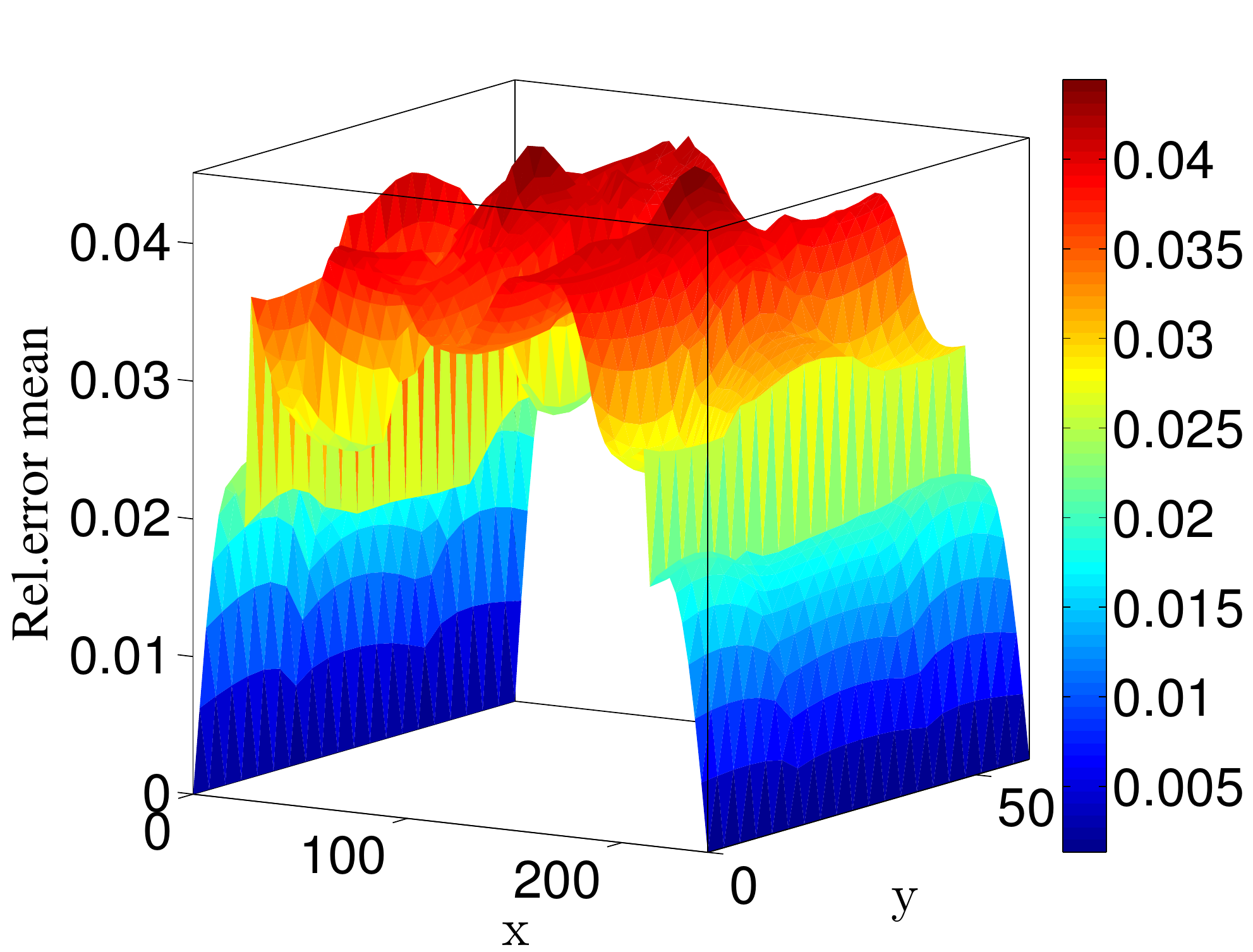}
        \caption{Error, $\xi, d=40; \eta, r=5,  |D|=27$} \label{RK:fig:u_mean_mc_xi_d40_src0_sink1_mc100000_bcs50_25_etad5_p3_27DOM_rel_error}
    \end{subfigure} 
    
    \caption{Relative error of the mean of the solution obtained by stochastic basis adaptation and domain decomposition (3, 8, 15, and 27 subdomains) with random variables $\eta$ and dimension, $r$ = 5, order, $p$ = 3 and and sparse-grid level, $l=5$. The reference solution is computed with $d=40$ random variables in $\xi$ using 100000 Monte Carlo simulations.} \label{RK:fig:u_mean_mc_xi_d40_src0_sink1_mc100000_bcs50_25_etad5_rel_error}
\end{figure}

\begin{figure}[t!]
    \centering
    \begin{subfigure}[t]{0.45\textwidth}
        \centering
       \includegraphics[height=1.8in]{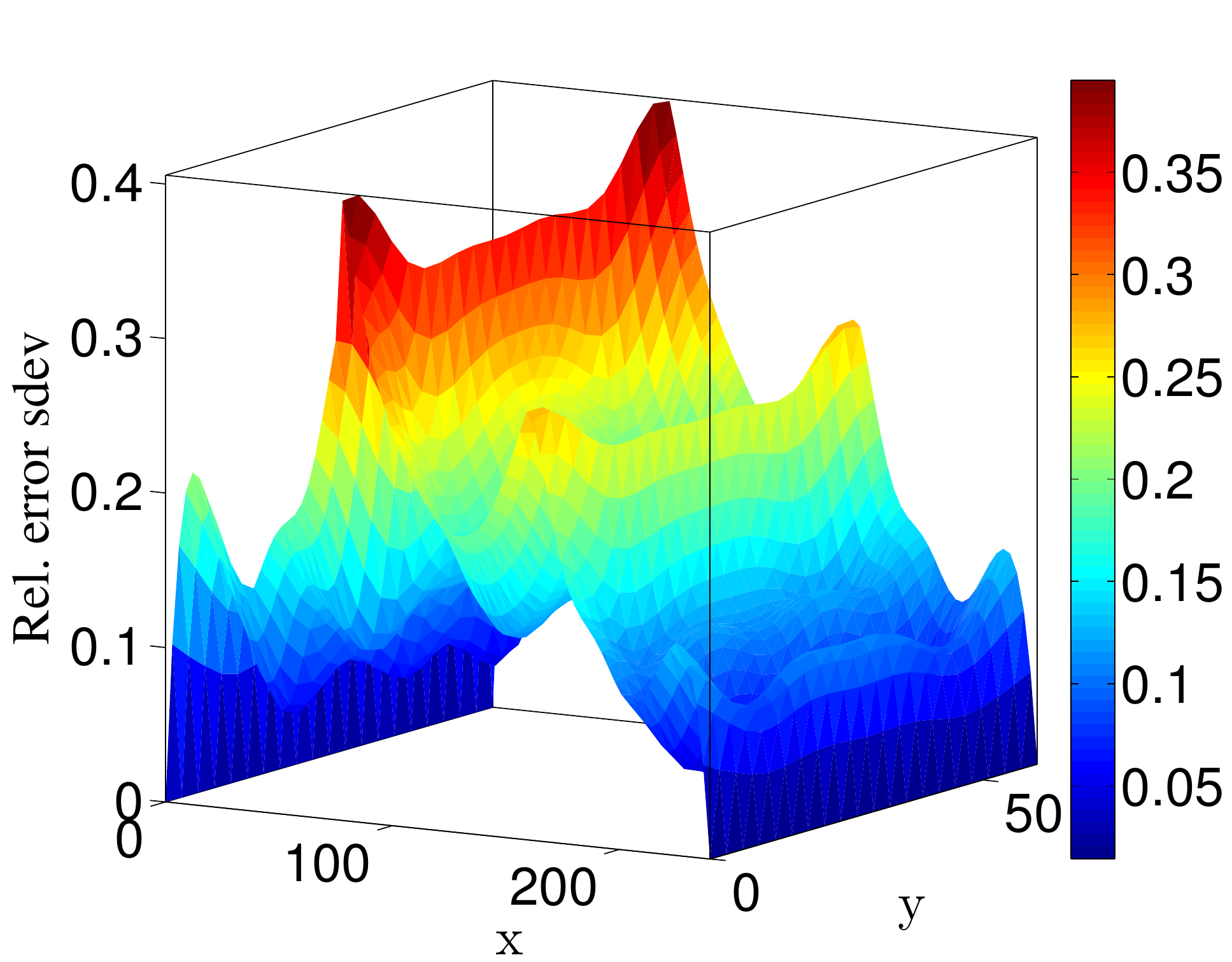}
        \caption{Error, $\xi, d=40; \eta, r=5,  |D|=3$} \label{RK:fig:u_sdev_mc_xi_d40_src0_sink1_mc100000_bcs50_25_etad5_p3_3DOM_rel_error}
    \end{subfigure}        
    \begin{subfigure}[t]{0.45\textwidth}
        \centering
        \includegraphics[height=1.8in]{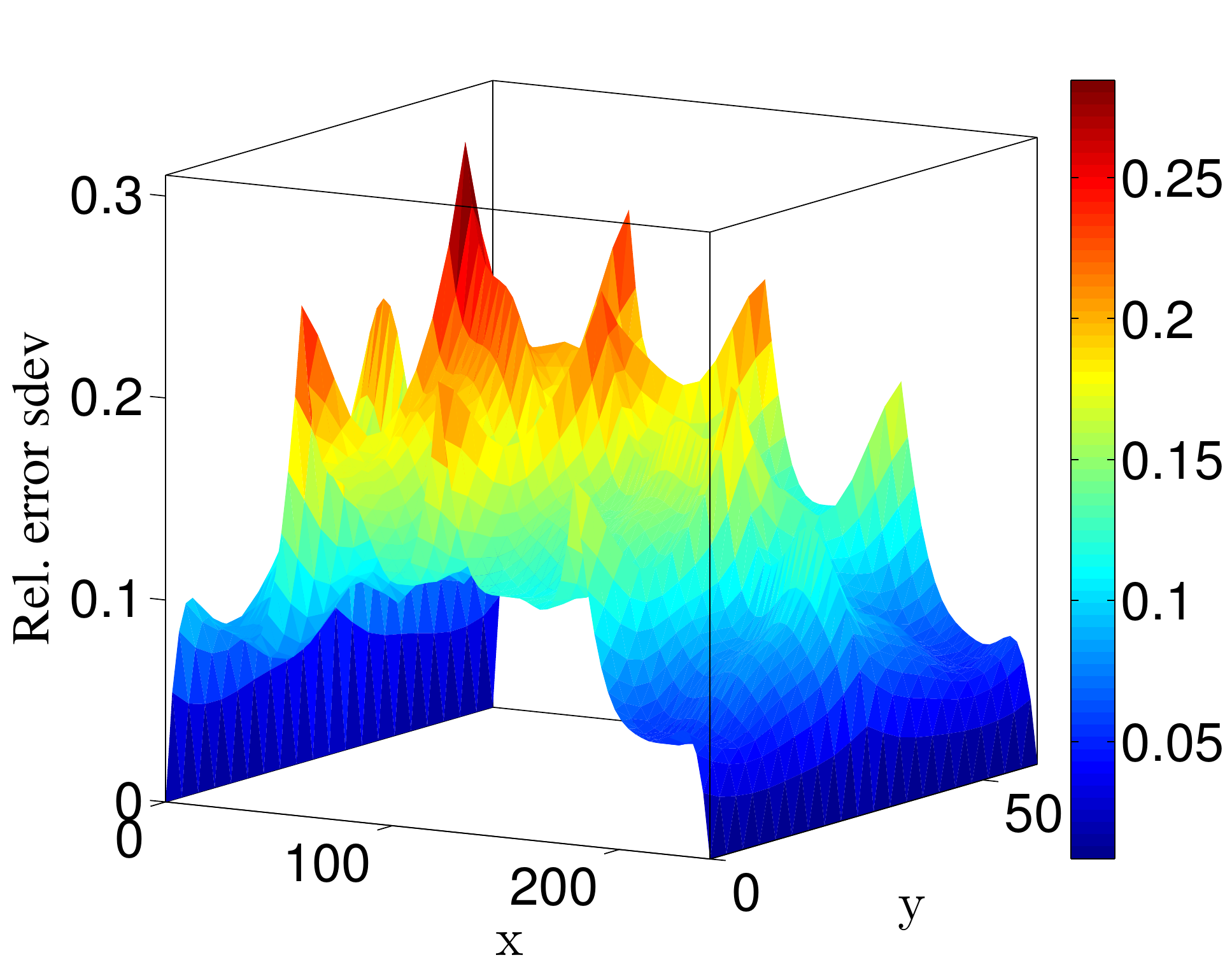}
        \caption{Error, $\xi, d=40; \eta, r=5,  |D|=8$} \label{RK:fig:u_sdev_mc_xi_d40_src0_sink1_mc100000_bcs50_25_etad5_p3_8DOM_rel_error}
    \end{subfigure}    
    \begin{subfigure}[t]{0.45\textwidth}
        \centering
        \includegraphics[height=1.8in]{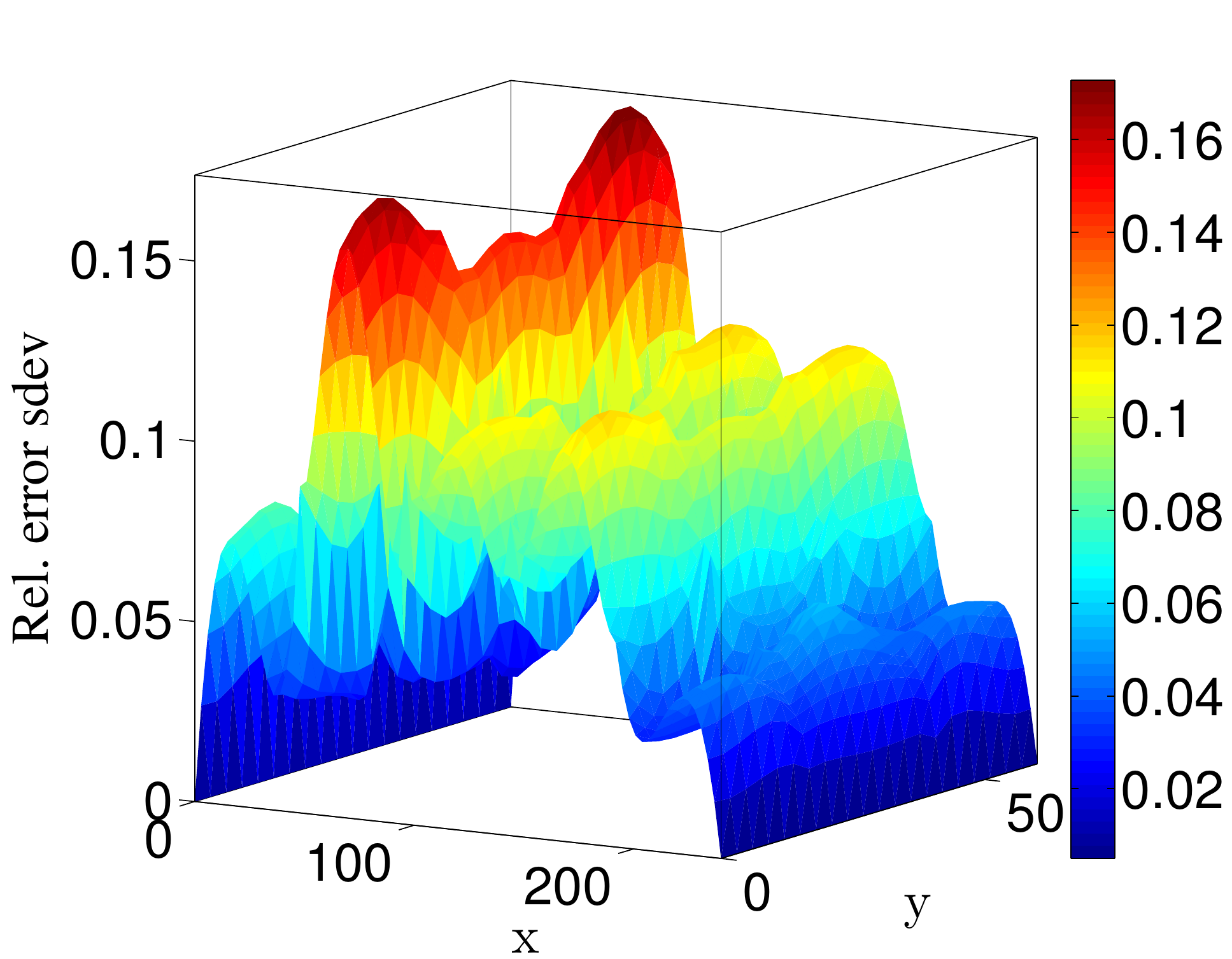}
        \caption{Error, $\xi, d=40; \eta, r=5,  |D|=15$} \label{RK:fig:u_sdev_mc_xi_d40_src0_sink1_mc100000_bcs50_25_etad5_p3_15DOM_rel_error}
    \end{subfigure}  
     \begin{subfigure}[t]{0.45\textwidth}
        \centering
        \includegraphics[height=1.8in]{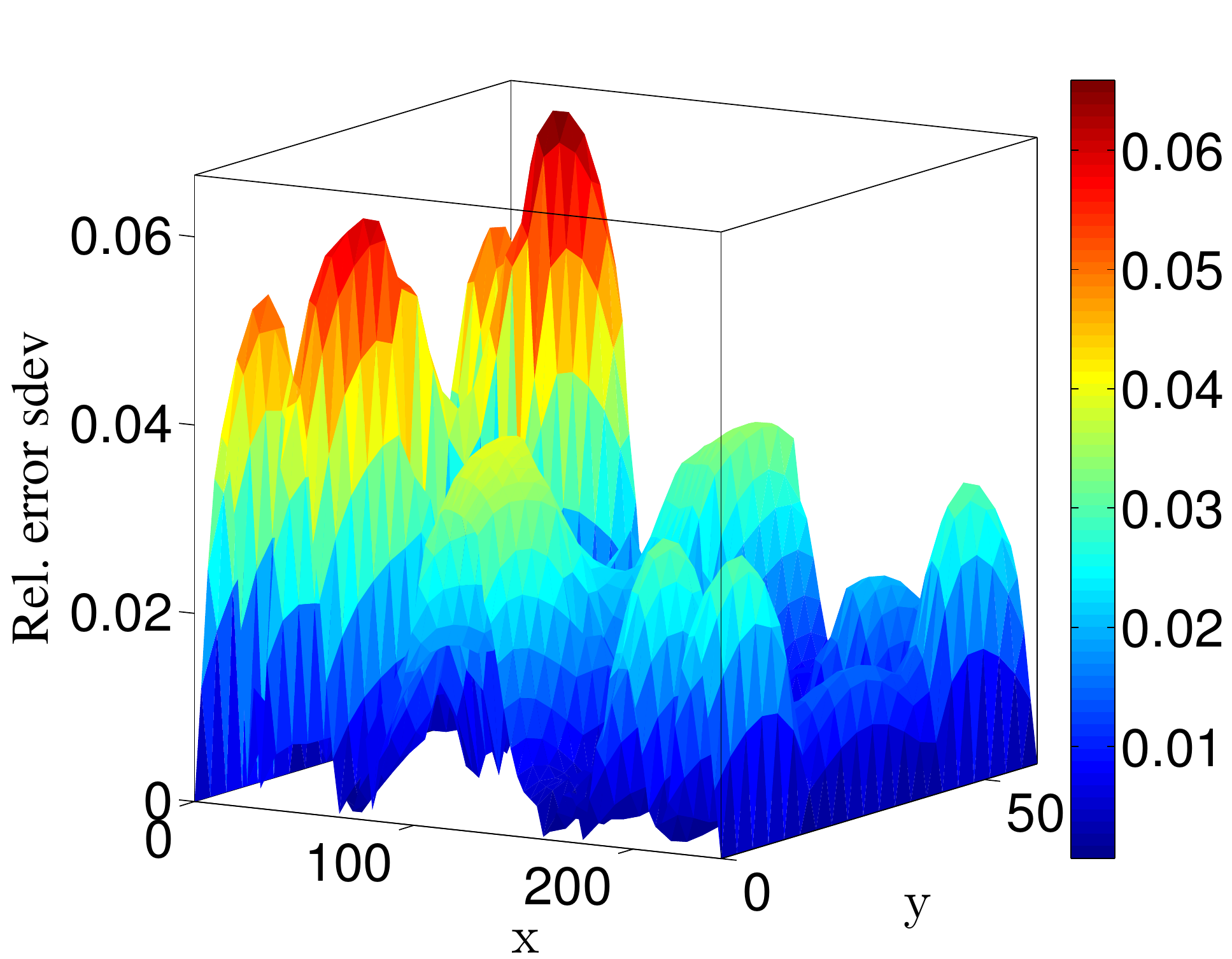}
        \caption{Error, $\xi, d=40; \eta, r=5,  |D|=27$} \label{RK:fig:u_sdev_mc_xi_d40_src0_sink1_mc100000_bcs50_25_etad5_p3_27DOM_rel_error}
    \end{subfigure} 
    
    \caption{Relative error of the standard deviation of the solution obtained by stochastic basis adaptation and domain decomposition (3, 8, 15, and 27 subdomains) with random variables $\eta$ and dimension, $r$ = 5, order, $p$ = 3 and sparse-grid level, $l=5$. The reference solution is computed with $d=40$ random variables in $\xi$ using 100000 Monte Carlo simulations.} \label{RK:fig:u_sdev_mc_xi_d40_src0_sink1_mc100000_bcs50_25_etad5_rel_error}
\end{figure}

For the $40$-dimensional case, we plot the mean from the reference solution with $100000$ Monte Carlo simulations (see Fig. \ref{RK:fig:mean_xid40_etad4_DOM27}) and from the reduced solution with stochastic dimension $4$ and number of subdomains $27.$ We also plot the relative error for the mean. Fig. \ref{RK:fig:sdev_xid40_etad4_DOM27} depicts the corresponding plots for the standard deviation.  Fig. \ref{RK:fig:mean_mc_xi_d40_src0_sink1_mc100000_bcs50_25_etad_p3_L2_sqrtn_rel_error_plot} and Fig. \ref{RK:fig:sdev_mc_xi_d40_src0_sink1_mc100000_bcs50_25_etad_p3_L2_sqrtn_rel_error_plot} show the respective relative error for the mean and standard deviation for various combinations of $r$ and number of subdomains. These plots show that accuracy increases with increasing $r$.
Fig. \ref{RK:fig:u_mean_mc_xi_d40_src0_sink1_mc100000_bcs50_25_etad5_rel_error} and Fig. \ref{RK:fig:u_sdev_mc_xi_d40_src0_sink1_mc100000_bcs50_25_etad5_rel_error} show the relative error plots for the mean and standard deviation, respectively, for reduced dimension $r=5$ and various numbers of subdomains. These plots show improvement in accuracy as the number of domains increases. 

\clearpage 

\begin{table}[h]
\caption{L2 relative error of the mean and standard deviation and cost ratio for the solution obtained by stochastic basis adaptation and domain decomposition (3, 8, and 15 subdomains) with random variables $\eta$ and dimension, $r$ = 3, 4, and 5, order, $p$ = 3 and sparse-grid level, $l=5$. The reference solution is computed with $d=10$ random variables in $\xi$ using the sparse-grid collocation method with level 5.}
\label{RK:tab:10d}
\begin{center}
\begin{tabular}{|c|c|c|c|c|c|c|c|c|c|}
 \hline
   & \multicolumn{3}{|c|}{r=3} &
 \multicolumn{3}{|c|}{r=4} &\multicolumn{3}{|c|}{r=5} \\ \cline{2-10} 
   $N_D$ & $\mu_e$(\%)  & $\sigma_e$(\%)   & CR & $\mu_e$(\%)   & $\sigma_e$ (\%)  & CR & $\mu_e$ (\%)  & $\sigma_e$(\%)   & CR  \\ \hline
  3  &1.59  & 6.37 & 237 & 1.07 &2.79 & 151   & 0.74  &  1.56  & 91   \\ \hline
  8  & 0.96  & 3.35 & 367  & 0.66  &1.99 & 316  & 0.49  &  1.59  & 248  \\ \hline
  15  & 0.80  & 2.13 & 304  & 0.54  & 2.29 & 218 & 0.64  &  2.08 & 141  \\ 
 \hline
\end{tabular}
\end{center}
\end{table}

\begin{table}[h]
\caption{L2 relative error of the mean and standard deviation and cost ratio for the solution obtained by stochastic basis adaptation and domain decomposition (3, 8, 15, and 27 subdomains)  with random variables $\eta$ and dimension, $r$ = 3, 4, and 5, order, $p$ = 3 and sparse-grid level, $l=5$. The reference solution is computed with $d=40$ random variables in $\xi$ using 100000 Monte Carlo simulations.}
\label{RK:tab:40d}
\begin{center}
\begin{tabular}{|c|c|c|c|c|c|c|c|c|c|}
 \hline
   & \multicolumn{3}{|c|}{r=3} &
 \multicolumn{3}{|c|}{r=4} &\multicolumn{3}{|c|}{r=5} \\ \cline{2-10} 
   $N_D$ & $\mu_e$(\%)  & $\sigma_e$(\%)   & CR & $\mu_e$(\%)   & $\sigma_e$ (\%)  & CR & $\mu_e$ (\%)  & $\sigma_e$(\%)   & CR  \\ \hline
  3  &7.94  &21.5 & 1032 & 7.52 &19.6 & 847   & 6.88  &  16.3  & 640   \\ \hline
  8  & 7.09  &17.9 & 1194  & 6.62  &15.6 & 1140  & 6.06  &  13.1  & 1049  \\ \hline
  15  & 5.94  & 12.3 & 1126  & 5.69  & 11.3 & 998 & 4.98  &  8.26 & 818  \\ \hline
  27  & 4.07  &  4.23  & 915 & 3.62  &  2.93  & 673   & 3.25  &  2.18  & 450  \\  
 \hline
\end{tabular}
\end{center}
\end{table}

\clearpage 

\section{Conclusions}\label{conclusions}
We have presented a novel approach  for time-independent PDEs with random parameters. In this approach, we decomposed the spatial domain into a set of non-overlapping subdomains and used stochastic basis adaptation methods to compute a low-dimensional representation of the solution in each subdomain. For each adapted basis, the solution was computed locally, subject to the interface compatibility conditions between disjoint subdomains. Specifically, we used a non-iterative version of the N-N algorithm to find the solution in each subdomain. This involved two steps: 1) solution at the interface of the subdomains and 2) solution in the interior of each subdomain. The interior solution of each subdomain was computed independently, which makes our method highly parallelizable.

We presented several numerical examples including a one-dimensional non-linear diffusion equations and two-dimensional linear diffusion equation. We provided results comparing the mean, standard deviation, and pdfs for solutions computed with both full and reduced dimensional representations, and estimates of the  computational cost in terms of flops for full and reduced dimensional solutions. Our numerical experiments demonstrate that the low-dimensional solutions agree well with high-dimensional solutions at a significantly smaller computational cost.

\section{Acknowledgments}
This research was supported by the U.S. Department of Energy, Office of Science, Office of Advanced 
Scientific Computing Research as part of the ``Uncertainty Quantification For Complex Systems Described by Stochastic Partial Differential Equations'' project.  Pacific Northwest National Laboratory is operated by Battelle for the DOE under Contract DE-AC05-76RL01830. We would like to thank the anonymous reviewers for their insightful comments which helped us improve the manuscript significantly.

 \bibliographystyle{siamplain} 
 
 \bibliography{StochDD_BA_NN_Schur}

\begin{thebibliography}{10}

\bibitem{Babuka2002}
{\sc I.~Babu{\v{s}}ka and P.~Chatzipantelidis}, {\em On solving elliptic
  stochastic partial differential equations}, Computer Methods in Applied
  Mechanics and Engineering, 191 (2002), pp.~4093--4122.

\bibitem{Babuka2010}
{\sc I.~Babu{\v{s}}ka, F.~Nobile, and R.~Tempone}, {\em A stochastic
  collocation method for elliptic partial differential equations with random
  input data}, {SIAM} Review, 52 (2010), pp.~317--355.

\bibitem{Berkooz1993539}
{\sc G.~Berkooz, P.~Holmes, and J.~Lumley}, {\em The proper orthogonal
  decomposition in the analysis of turbulent flows}, Annual Review of Fluid
  Mechanics, 25 (1993), pp.~539--575.

\bibitem{Cameron1947}
{\sc R.~H. Cameron and W.~T. Martin}, {\em The orthogonal development of
  non-linear functionals in series of {F}ourier-{H}ermite functionals}, The
  Annals of Mathematics, 48 (1947), p.~385.

\bibitem{Chen2015}
{\sc Y.~Chen, J.~Jakeman, C.~Gittelson, and D.~Xiu}, {\em Local polynomial
  chaos expansion for linear differential equations with high dimensional
  random inputs}, {SIAM} Journal on Scientific Computing, 37 (2015),
  pp.~A79--A102.

\bibitem{Christensen1999}
{\sc E.~A. Christensen., M.~Br{\o}ns, and J.~N. S{\o}rensen}, {\em Evaluation
  of proper orthogonal decomposition--based decomposition techniques applied to
  parameter-dependent nonturbulent flows}, {SIAM} Journal on Scientific
  Computing, 21 (1999), pp.~1419--1434.

\bibitem{RK:Doostan2007}
{\sc A.~Doostan, R.~G. Ghanem, and J.~Red-Horse}, {\em Stochastic model
  reduction for chaos representations}, Computer Methods in Applied Mechanics
  and Engineering, 196 (2007), pp.~3951--3966.

\bibitem{RK:Doostan2011}
{\sc A.~Doostan and H.~Owhadi}, {\em A non-adapted sparse approximation of pdes
  with stochastic inputs}, Journal of Computational Physics, 230 (2011),
  pp.~3015--3034.

\bibitem{Ghanem1999}
{\sc R.~Ghanem}, {\em The nonlinear gaussian spectrum of log-normal stochastic
  processes and variables}, Journal of Applied Mechanics, 66 (1999), p.~964.

\bibitem{Ghanem2015}
{\sc R.~Ghanem and C.~Soize}, {\em Remarks on stochastic properties of
  materials through finite deformations}, International Journal for Multiscale
  Computational Engineering, 13 (2015), pp.~367--374.

\bibitem{RK:Ghanem1991}
{\sc R.~Ghanem and P.~Spanos}, {\em Stochastic Finite Elements: A Spectral
  Approach}, Springer-Verlag, 1991.

\bibitem{RK:Kirby1992}
{\sc M.~Kirby}, {\em Minimal dynamical systems from {PDE}s using {S}obolev
  eigenfunctions}, Physica D: Nonlinear Phenomena, 57 (1992), pp.~466--475.

\bibitem{RK:Levy1999}
{\sc A.~Levy and J.~Rubinstein}, {\em {H}ilbert--space {K}arhunen--{L}o\'eve
  transform with application to image analysis}, Journal of The Optical Society
  of America A, 16 (1999), pp.~28--35.

\bibitem{Lin2009AWR}
{\sc G.~Lin and A.~Tartakovsky}, {\em An efficient, high-order probabilistic
  collocation method on sparse grids for three-dimensional flow and solute
  transport in randomly heterogeneous porous media}, Advances in Water
  Resources, 32 (2009), pp.~712--722.

\bibitem{Lin2010JCP}
{\sc G.~Lin, A.~Tartakovsky, and D.~Tartakovsky}, {\em Uncertainty
  quantification via random domain decomposition and probabilistic collocation
  on sparse grids}, Journal of Computational Physics, 229 (2010),
  pp.~6995--7012.

\bibitem{Lin2010JSC}
{\sc G.~Lin and A.~M. Tartakovsky}, {\em Numerical studies of three-dimensional
  stochastic {D}arcy's equation and stochastic advection-diffusion-dispersion
  equation}, Journal of Scientific Computing, 43 (2010), pp.~92--117.

\bibitem{RK:Loeve1977}
{\sc M.~Lo\'eve}, {\em Probability Theory}, Springer-Verlag, 1977.

\bibitem{RK:Nobile2008}
{\sc F.~Nobile, R.~Tempone, and C.~G. Webster}, {\em A sparse grid stochastic
  collocation method for partial differential equations with random input
  data}, SIAM Journal on Numerical Analysis, 46 (2008), pp.~2309--2345.

\bibitem{Nouy2007}
{\sc A.~Nouy}, {\em A generalized spectral decomposition technique to solve a
  class of linear stochastic partial differential equations}, Computer Methods
  in Applied Mechanics and Engineering, 196 (2007), pp.~4521--4537.

\bibitem{Pan1995}
{\sc L.~Pan and P.~J. Wierenga}, {\em A transformed pressure head-based
  approach to solve richards{\textquotesingle} equation for variably saturated
  soils}, Water Resources Research, 31 (1995), pp.~925--931.

\bibitem{Pranesh2016}
{\sc S.~Pranesh and D.~Ghosh}, {\em Addressing the curse of dimensionality in
  {SSFEM} using the dependence of eigenvalues in {KL} expansion on domain
  size}, Computer Methods in Applied Mechanics and Engineering, 311 (2016),
  pp.~457--475.

\bibitem{RK:Silverman1996}
{\sc B.~W. Silverman}, {\em Smoothed functional principal components analysis
  by choice of norm}, The Annals of Statistics, 24 (1996), pp.~1--24.

\bibitem{RK:Smolyak1963}
{\sc S.~Smolyak}, {\em Quadrature and interpolation formulas for tensor
  products of certain classes of functions}, Doklady {A}kademii {N}auk {SSSR},
  4 (1963), pp.~240--243.

\bibitem{Soize2015}
{\sc C.~Soize}, {\em Polynomial chaos expansion of a multimodal random vector},
  {SIAM}/{ASA} Journal on Uncertainty Quantification, 3 (2015), pp.~34--60.

\bibitem{Soize2016}
{\sc C.~Soize and R.~Ghanem}, {\em Data-driven probability concentration and
  sampling on manifold}, Journal of Computational Physics, 321 (2016),
  pp.~242--258.

\bibitem{Tartakovsky2008}
{\sc A.~M. Tartakovsky, D.~Bolster, and D.~M. Tartakovsky}, {\em
  Hydrogeophysical approach for identification of layered structures of the
  vadose zone from electrical resistivity data}, Vadose Zone Journal, 7 (2008),
  p.~1253.

\bibitem{Tipireddy2013}
{\sc R.~Tipireddy}, {\em {Algorithms for stochastic Galerkin projections:
  Solvers, basis adaptation and multiscale modeling and reduction}}, {T}heses,
  {University of Southern California}, Aug. 2013.

\bibitem{RK:Tipireddy2014}
{\sc R.~Tipireddy and R.~Ghanem}, {\em Basis adaptation in homogeneous chaos
  spaces}, Journal of Computational Physics, 259 (2014), pp.~304--317.

\bibitem{RK:Tipireddy2016}
{\sc R.~Tipireddy, P.~Stinis, and A.~Tartakovsky}, {\em Basis adaptation and
  domain decomposition for steady partial differential equations with random
  coefficients}, arXiv:1607.08280,  (2016).

\bibitem{RK:Toselli2005}
{\sc A.~Toselli and O.~B. Widlund}, {\em Domain Decomposition
  Methods--Algorithms and Theory}, Springer-Verlag, 2005.

\bibitem{RK:Tsilifis2016}
{\sc P.~Tsilifis and R.~Ghanem}, {\em Reduced {Wiener Chaos} representation of
  random fields via basis adaptation and projection}, arXiv:1603.04803v3,
  (2016).

\bibitem{Venturi2013JCP}
{\sc D.~Venturi, D.~Tartakovsky, A.~Tartakovsky, and G.~Karniadakis}, {\em
  Exact {PDF} equations and closure approximations for advective-reactive
  transport}, Journal of Computational Physics, 243 (2013), pp.~323--343.

\bibitem{RK:Xiu2002}
{\sc D.~Xiu and G.~E. Karniadakis}, {\em The {W}iener--{A}skey polynomial chaos
  for stochastic differential equations}, SIAM Journal on Scientific Computing,
  24 (2002), pp.~619--644.

\bibitem{Xiu2004662}
{\sc D.~Xiu and D.~Tartakovsky}, {\em A two-scale nonperturbative approach to
  uncertainty analysis of diffusion in random composites}, Multiscale Modeling
  Simulation, 2 (2004), pp.~662--674.

\end{thebibliography}

\end{document}